\documentclass{amsart}
\usepackage{amsfonts,amssymb,amsmath}

  \def\F{\mathbb F}
\def\Fq{{\mathbb F}_q}

  \def\a{\alpha}

  \def\Ga{\Gamma}
\def\d{\delta}
   \def\D{\Delta}
  \def\e{\epsilon}
  \def\l{\lambda}

  \def\la{\langle}
  \def\ra{\rangle}
  \def\s{\sigma}

  \def\o{\omega}
    \def\Om{\Omega}
  \def\no{\noindent}
  
  \def\pf{\noindent {\bf Proof $\;$ }}
  
\def\vp{\varphi}
\def\M{{\rm M}}

\def\ve{\varepsilon}
\newcommand{\End}{{\mathrm {End}}}
\newcommand{\Syl}{{\mathrm {Syl}}}

\newcommand{\fs}{{f^*}}
\newcommand{\Hom}{{\mathrm {Hom}}}
\newcommand{\AAA}{{\sf A}}
\newcommand{\SSS}{{\sf S}}

\newcommand{\FQ}{\mathbb{F}_{q}}

\newcommand{\Ker}{{\rm Ker}}

\newcommand{\Irr}{{\rm Irr}}

  \def\hal{\unskip\nobreak\hfil\penalty50\hskip10pt\hbox{}\nobreak
  \hfill\vrule height 5pt width 6pt depth 1pt\par\vskip 2mm}

\newtheorem{theorem}{Theorem}
  \newtheorem{thm}{Theorem}[section]
  \newtheorem{prop}[thm]{Proposition}
  
  \newtheorem{lemma}[thm]{Lemma}
  \newtheorem{cor}[thm]{Corollary}
 \newtheorem{coroll}[theorem]{Corollary}
  \newtheorem{remark}[thm]{Remark}

\begin{document}

\title{Arithmetic results on orbits of linear groups}

\author[M. Giudici]{Michael Giudici}
\email{michael.giudici@uwa.edu.au}

\author[M.~W. Liebeck]{Martin W. Liebeck}
\email{m.liebeck@imperial.ac.uk}

\author[C.~E. Praeger]{Cheryl E. Praeger}
\email{Cheryl.Praeger@uwa.edu.au}

\author[J. Saxl]{Jan Saxl}
\email{saxl@dpmms.cam.ac.uk}

\author[P.~H. Tiep]{Pham Huu Tiep}
\email{tiep@math.arizona.edu}

\address[Michael Giudici and Cheryl E. Praeger]{Centre for the Mathematics of Symmetry and Computation\\
School of Mathematics and Statistics\\ 
University of Western Australia, 35 Stirling Highway\\ 
Crawley, Western Australia 6009}

\address[Martin W. Liebeck]{Department of Mathematics, Imperial College, London SW7 2BZ, England}
\address[Jan Saxl]{DPMMS, CMS, University of Cambridge, Wilberforce Road, Cambridge CB3 0WB, England}
\address[Pham Huu Tiep]{Department of Mathematics\\University of Arizona\\617 N. Santa Rita Ave.\\Tucson, AZ 85721-0089}

\begin{abstract}
Let $p$ be a prime and $G$ a subgroup of $GL_d(p)$. We define $G$ to be $p$-{\it exceptional} if it has order divisible by $p$, but all its orbits on vectors have size coprime to $p$. We obtain a classification of $p$-exceptional linear groups. This has consequences for a well known conjecture in representation theory, and also for a 
longstanding question concerning $\frac{1}{2}$-transitive linear groups (i.e. those having all orbits on nonzero vectors of equal length), classifying those of order divisible by $p$.
\end{abstract}

\maketitle


\section{Introduction}

The study of orbits of linear groups acting on finite vector spaces has  a long history. Zassenhaus \cite{Z} investigated linear groups for which all orbits on nonzero vectors are regular, classifying the insoluble examples, i.e. the insoluble Frobenius complements. If one merely assumes that there is at least one regular orbit, there are many examples and the investigation and classification of these is a lively area of current research. For example, if $p$ is the characteristic and $G$ is a quasisimple irreducible $p'$-group, there is almost always a regular orbit, the exceptions being classified in \cite{goodwin, KP}; this played a major role in the solution of the $k(GV)$-problem \cite{kgv}. In a different direction, linear groups acting transitively on the set of nonzero vectors were determined by Hering \cite{her}, leading to the classification of 2-transitive permutation groups of affine type. Results on groups with few orbits, or a long orbit,  or orbits with coprime lengths, can be found in \cite{dol, liebaff, malle}. A much weaker assumption than transitivity is that of $\frac{1}{2}$-{\it transitivity} -- namely, that all orbits on nonzero vectors have the same size. The soluble linear groups with this property were classified by Passman \cite{passman1, passman2}.

In this paper we study linear groups with the following property.

\vspace{4mm}
\no {\bf Definition } Let $V=V_d(p)$ be a vector space of dimension $d$ over $\F_p$ with $p$ prime, and let $G \le GL_d(p) =  GL(V)$. We say that $G$ is $p$-{\it exceptional} if $p$ divides $|G|$ and $G$ has no orbits
on $V$ of size divisible by $p$.

\vspace{4mm}
Note that if $d=ab$ for positive integers $a,b$, and $q = p^a$, then $\Ga L_b(q) \le GL_d(p)$, so the above definition also applies to subgroups of  $\Ga L_b(q)$.

If $G \le GL_d(p)$ has a regular orbit on vectors, then $G$ is not $p$-exceptional. On the other hand, if $G$ is transitive (or $\frac{1}{2}$-transitive) on nonzero vectors and has order divisible by $p$, then $G$ is $p$-exceptional. 

We shall obtain a classification of all $p$-exceptional linear groups, up to some undecided questions in the imprimitive case. We also give applications to  $\frac{1}{2}$-transitive groups, and to a conjecture in representation theory.

We begin with our result for primitive linear groups. In the statement, by the {\it deleted permutation module} over $\F_p$ ($p$ prime) for a symmetric group $\SSS_c$, we mean the irreducible $\F_p\SSS_c$-module $S/S\cap T$, where
$S = \{(a_1,\ldots ,a_c) : a_i \in \F_p, \sum a_i = 0\}$ and $T = \{(a,\ldots ,a) : a \in \F_p\}$, and $\SSS_c$ acts by permuting coordinates in the obvious way. Denote by $V^\sharp$ the set of nonzero vectors in a vector space $V$.

\begin{theorem}\label{prim}
Let $G$ be an irreducible 
$p$-exceptional subgroup of $GL_d(p) = GL(V)$, and suppose 
$G$ acts primitively on $V$. Then one of the following holds: 
\begin{enumerate}
\item[{\rm (i)}] $G$ is transitive on $V^\sharp$
(a list can be found in \cite[Appendix 1]{liebaff});
\item[{\rm (ii)}]  $G \le \Gamma L_1(q)$ ($q=p^d$), determined in Lemma $\ref{lem:1dim}$;
\item[{\rm (iii)}] $G$ is one of the following:
\begin{enumerate}
\item[{\rm (a)}] $G = \AAA_c$, $\SSS_c$ with $c = 2^r-2$ or $2^r-1$, with $V$ the deleted
permutation module over $\F_2$, of dimension $c-2$ or $c-1$ respectively (see Lemma $\ref{delete}$);
\item[{\rm (b)}] $SL_2(5) \trianglelefteq G < \Ga L_2(9)< GL_4(3)$, orbit sizes $1,40,40$;
\item[{\rm (c)}] $L_2(11) \trianglelefteq G < GL_5(3)$, orbit sizes $1,22,110,110$;
\item[{\rm (d)}] $M_{11} \trianglelefteq G < GL_5(3)$, orbit sizes $1,22,220$;
\item[{\rm (e)}] $M_{23} = G < GL_{11}(2)$, orbit sizes $1,23,253,1771.$
\end{enumerate}
\end{enumerate}
\end{theorem}

For the imprimitive case we first require a result on permutation groups.
For a prime $p$, we say a subgroup $K \le \SSS_n$
is $p$-{\it concealed} if it has order divisible by $p$,
and all its orbits on the power set of $\{1,\ldots ,n\}$ have size
coprime to $p$. The following result is an extension of \cite{CNS, Ser}, which classify primitive 
groups having no regular orbit on the power set.

\begin{theorem}\label{p-conc}
Let $H$ be a primitive subgroup of $\SSS_n$ of order divisible by a prime $p$. 
Then $H$ is $p$-concealed if and only if 
one of the following holds:

{\rm (i)} $\AAA_n \trianglelefteq H \leq \SSS_n$, and $n = ap^s-1$ 
with $s \geq 1$, $a \leq p-1$ and $(a,s) \neq (1,1)$; also 
$H \neq \AAA_3$ if $(n,p) = (3,2)$;

{\rm (ii)} $(n,p) = (8,3)$, and $H = AGL_3(2) = 2^3:SL_3(2)$ or 
$H = A\Gamma L_1(8) = 2^3:7:3$;

{\rm (iii)} $(n,p) = (5,2)$ and $H = D_{10}$, 
a dihedral group of order $10$. 
\end{theorem}

Theorem \ref{p-conc} will be proved in Section \ref{c2sec}.

Here is our main result on imprimitive $p$-exceptional linear groups.

\begin{theorem}\label{imprim}
Suppose $G \le GL_d(p) = GL(V)$ is irreducible, $p$-exceptional and 
imprimitive, and also $G = O^{p'}(G)$.
Let $V = V_1 \oplus \ldots \oplus V_n$ ($n>1$)
be any imprimitivity decomposition for $G$. Then 
$G_{V_1}$ is transitive on $V_1^\sharp$,
and $G$ induces a primitive $p$-concealed subgroup of $\SSS_n$ on 
$\{V_1,\ldots ,V_n\}$.
\end{theorem}

There is a partial converse: if $X \le GL(V_1)$ is transitive on 
$V_1^\sharp$ and $H \le \SSS_n$ is primitive and $p$-concealed, then the
full wreath product $X \wr H$ acting on $V = V_1^n$ is $p$-exceptional (see Lemma \ref{wreath}).

The following is a general structure theorem for irreducible $p$-exceptional groups. 

\begin{theorem}\label{gen}
Let $G \le GL_d(p) = GL(V)$ be an irreducible $p$-exceptional group, and let 
$G_0 = O^{p'}(G)$. Write $V \downarrow G_0 = V_1\oplus \cdots \oplus V_t$
with $V_i$ irreducible $G_0$-modules. 
Then $G_0^{V_1}$ is either a primitive $p$-exceptional group (given by
Theorem \ref{prim}), or an imprimitive $p$-exceptional group (given by 
Theorem \ref{imprim}). Moreover,  the $V_i$ are pairwise non-isomorphic $G_0$-modules, and 
$G$ acts on $\{V_1,\ldots ,V_t\}$ as a transitive $p'$-group.
\end{theorem}

Again, there is a partial converse (Lemma \ref{wreath}): the full wreath product
of a $p$-exceptional group and a transitive $p'$-group is $p$-exceptional.

The next result has important applications in the modular 
representation theory of finite groups. Recall that, if $p$ is any prime and 
$B$ is a Brauer $p$-block of any finite group $G$ with defect group $P$,
then the {\it Brauer height zero conjecture} asserts that
all irreducible complex characters in $B$ have height
zero if and only if $P$ is abelian. One of the significant results of the
representation theory of finite groups in the 1980's was to prove  that if
$G$ is $p$-soluble and $\lambda \in \Irr(Z)$ is
an irreducible complex character of a normal
subgroup $Z \lhd G$   such that $\chi(1)/\lambda(1)$ is
not divisible by $p$ for all $\chi \in \Irr(G)$
lying over $\lambda$, then $G/Z$ has abelian Sylow
$p$-subgroups. This theorem, established by D. Gluck and T. Wolf in \cite{GW1, GW2}
led to a proof of the Brauer height zero conjecture
for $p$-soluble groups. As shown by 
very recent results on the Brauer height zero conjecture,
in particular, the proof \cite{NT1} of the conjecture in the case 
$p=2$ and $P \in \Syl_p(G)$, as well as reduction theorems for the conjecture
\cite{Mu} and \cite{NS}, one of the main obstacles to proving
the conjecture in full generality is
to obtain a proof of the Gluck-Wolf theorem for
arbitrary finite groups. This has now been achieved in \cite{NT2}, which uses
Theorem \ref{tiepnav} in a crucial way. 

\begin{theorem}\label{tiepnav} 
Let $G$ be a non-identity finite group
and let $p$ be an odd prime. Assume that $G = O^{p'}(G) = O^{p}(G)$ and
$G$ has abelian Sylow $p$-subgroups. Suppose that
$V$ is a finite-dimensional, faithful, irreducible $\F_pG$-module 
such that every orbit of $G$ on $V$ has length coprime to $p$. Then
one of the following holds:

{\rm (i)} $G = SL_{2}(q)$ and $|V| = q^2$ for some $q = p^a$;

{\rm (ii)} $G$ acts transitively on the $n$ summands of a decomposition
$V = \bigoplus^{n}_{i=1}V_i$, where $p < n < p^2$,  
$n \equiv -1 \hbox{ mod }p$. 
Furthermore, $G_{V_1}$ acts transitively on $V_1^\sharp$, 
and the action of $G$
on $\{V_1,...,V_n\}$ induces either $\AAA_n$, or the affine group 
$2^3:SL_3(2)$ for $(n,p) = (8,3)$;

{\rm (iii)} $(G,|V|) = (SL_2(5),3^4)$, $(2^{1+4}_{-}\cdot \AAA_5,3^4)$,
$(L_2(11),3^5)$, $(M_{11}, 3^5)$, or $(SL_2(13),3^6)$.
\end{theorem}

Here is a further consequence concerning $\frac{1}{2}$-transitive linear groups.

\begin{theorem}\label{threehalf}
Let $G\le GL_d(p)$ ($p$ prime) be a $\frac{1}{2}$-transitive linear group, and suppose $p$ divides $|G|$.
Then one of the following holds:

{\rm (i)} $G$ is transitive on the set of nonzero vectors (given by \cite[Appendix]{liebaff});

{\rm (ii)} $G \le \Ga L_1(p^d)$;

{\rm (iii)} $SL_2(5)\trianglelefteq G \le \Ga L_2(9) < GL_4(3)$ and $G$ has two orbits on nonzero vectors of size $40$.
\end{theorem}

Concerning part (ii) of the theorem, some examples of $\frac{1}{2}$-transitive subgroups of $\Ga L_1(p^d)$ are given in Lemma \ref{lem:1dim}.

 Recall that a finite transitive permutation group is said to be $\frac{3}{2}$-{\it transitive} if all nontrivial orbits of a point stabiliser have the same size, this size being greater than 1. By \cite[Theorem 10.4]{W}, such groups are either Frobenius groups or primitive. Steps towards the classification of the primitive examples were taken in \cite{bam}, where it was shown that they must be either almost simple or affine, and the former were classified. The next result deals with the modular affine case. The non-modular case will be the subject of a future paper.

\begin{coroll}\label{threecor}
If $G \le AGL_d(p)$ is a $\frac{3}{2}$-transitive affine permutation group with point-stabiliser $G_0$ of order divisible by $p$, then one of the following holds:

{\rm (i)} $G$ is $2$-transitive;

{\rm (ii)} $G \le A\Ga L_1(p^d)$;

{\rm (iii)} $SL_2(5)\trianglelefteq G_0 \le \Ga L_2(9) < GL_4(3)$ and $G$ has rank $3$ with subdegrees $1$, $40$, $40$.
\end{coroll}

\vspace{4mm}
The layout of the paper is as follows. Theorems \ref{p-conc} and \ref{imprim} are proved in Section \ref{c2sec}; then the proof of Theorem \ref{prim} is given in the next nine sections, culminating in Section \ref{primpf}. The deductions of Theorems \ref{gen},  \ref{tiepnav} and \ref{threehalf} can be found in the final Section \ref{ded}. 

\vspace{4mm}
{\bf Notation:} The following notation will be used throughout the paper. For a vector space~$V$ with subspace ~$U$, an element $g \in GL(V)$ and a subgroup $H \le GL(V)$,
\[
\begin{array}{l}
V^\sharp = V \setminus \{0\} \\
C_V(g) = \{v\in V \mid vg = v\} \\ 
C_V(H) = \{v\in V \mid vh = v \hbox{ for all }h\in H\} \\
V \downarrow H = \hbox{ restriction of }V \hbox{ to }H\\
H_U =\{h\in H\mid U^h=U\}, \textrm{ the setwise stabiliser of $U$ in $H$.}
\end{array}
\]
Moreover, if $H$ stabilises the subspace $U$ then $H^U$ is the subgroup of $GL(U)$ induced by~$H$.
Also $x^G$ denotes the conjugacy class of an element $x$ in a group $G$, 
and $J_k$ denotes a unipotent Jordan block of size $k$.
 
\vspace{4mm}
{\bf Acknowledgements } The first four authors acknowledge the support of an ARC Discovery Grant; Giudici acknowledges the support of an Australian Research Fellowship; Praeger acknowledges the support of a Federation Fellowship;
Tiep acknowledges the support of the NSF (grants DMS-0901241 and DMS-1201374).

\section{Preliminaries}\label{prel}

We begin with a simple observation.

\begin{lemma}\label{normalsubgroups} Let $H\leq GL_d(p)$ be $p$-exceptional on $V=\F_p^d$. 
\begin{enumerate}
 \item[{\rm (i)}] If $K$ is a normal subgroup of $H$ and $p$ divides $|K|$, then $K$ is $p$-exceptional on $V$.
\item[{\rm (ii)}] If $N\leq GL_d(p)$ has order coprime to $p$ and $N$ is normalised by $H$, then $NH$   is $p$-exceptional on $V$. 
\end{enumerate}
\end{lemma}

\pf (i) If $K$ had an orbit, say $\Delta$, in $V$ of length a multiple of $p$, then the $H$-orbit containing $\Delta$ would have length divisible by $|\Delta|$ since $K$ is normal in $H$, contradicting the fact that $H$ is $p$-exceptional.

(ii) Let $L:=NH$, let $v\in V^\sharp$, and consider $v^L$ and  $v^N$, the $L$-orbit and $N$-orbit containing $v$, respectively. Since $N$ is normal in $L$, $v^L$ is the set theoretic union of a subset $\mathcal{B}_0$ of the set $\mathcal{B}$ of $N$-orbits in $V$, and $\mathcal{B}_0$ is an $H$-orbit in its induced action on $\mathcal{B}$. Moreover $v^N\in \mathcal{B}_0$ and $|v^L|=|v^N|.|\mathcal{B}_0|$. As  $|N|$ is coprime to $p$, also $|v^N|$ is coprime to $p$. Since $H$ acts on $\mathcal{B}$, the $H$-orbit $v^H$ consists of a constant number of vectors from each $N$-orbit in $\mathcal{B}_0$. Thus $|\mathcal{B}_0|$ divides  $|v^H|$ and hence $|\mathcal{B}_0|$ is coprime to $p$. 
\hal

\begin{lemma}\label{em:related1}
Let $q=p^a$ with $p$ prime, let $Z = Z(GL_n(q))$ and let $H$ be a subgroup of $\Ga L_n(q)$. Then $H$ is $p$-exceptional if and only if $ZH$ is $p$-exceptional.
\end{lemma}

\pf If $ZH$ is $p$-exceptional then so is $H$, by Lemma \ref{normalsubgroups}(i). The converse follows from Lemma \ref{normalsubgroups}(ii). \hal

\begin{lemma}\label{pexc} 
Let $G\le \Ga L(V) = \Ga L_n(q)$ ($q=p^a$) be $p$-exceptional, and let $G_0 = G \cap GL(V)$. Then one of the following holds:

{\rm (i)} $p$ divides $|G_0|$ and $G_0$ is $p$-exceptional;

{\rm (ii)} $G_0$ is a $p'$-group, and $G$ contains a $p$-exceptional normal subgroup of the
form $G_0\la \s \ra$, where $\s \in \Gamma L(V) \backslash GL(V)$ is a field
automorphism of order $p$.
\end{lemma}

\pf If $p$ divides $|G_0|$ then $G_0$ is $p$-exceptional by Lemma \ref{normalsubgroups}, so (i) holds.
Now assume $p$ does not divide $|G_0|$. As $G/G_0$ is cyclic, we have $G = G_0\la x \ra$ for some $x$ of order divisible by $p$. Taking $\s$ to be a power of $x$ of order $p$, we obtain (ii) by applying Lemma \ref{normalsubgroups} to $G_0\langle\sigma\rangle$. \hal

\vspace{4mm}
The next lemma will be used many times in the proof of Theorem \ref{prim}.

\begin{lemma}\label{bound}
Let $G \le GL(V) = GL_d(q)$ with $q=p^a$ ($p$ prime), and suppose that $G$ is $p$-exceptional and $C_V(O^{p'}(G)) = 0$. Let $t$ be an element of $G$ of order $p$, and let $P \in Syl_p(G)$.

{\rm (i)} Then $d=\dim V \le r_p \log_q |G:N_G(P)|$, where $r_p$ is the minimal number of conjugates of $P$ generating $O^{p'}(G)$.

{\rm (ii)} We have $|V| \le |C_V(t)| \cdot |t^G|$.

{\rm (iii)} Suppose $O^{p'}(G)$ is generated by $\a$ conjugates of $t$. Then $q^{d/\a} \le |t^G|$.
\end{lemma}

\pf As $G$ is $p$-exceptional, every nonzero vector is fixed by some conjugate of $P$, so 
\[
V = \bigcup_{g\in G} C_V(P^g).
\]
Moreover, $\dim C_V(P) \le d(1-\frac{1}{r_p})$, since otherwise the group generated by $r_p$ conjugates of $P$ would have a nonzero centralizer in $V$, contrary to the hypothesis. Hence
\[
q^d = |V| \le |G:N_G(P)|\,q^{d(1-\frac{1}{r_p})}.
\]
This gives (i).

For (ii), observe that every nonzero vector in $V$ is fixed by a conjugate of $t$ (as $G$ is $p$-exceptional), so $V = \bigcup_{g\in G}C_V(t^g)$, which implies (ii). Finally, (iii) follows from (ii) together with the fact that $\dim C_V(t) \le d(1-\frac{1}{\a})$ (which follows from the argument of the first paragraph). \hal

\vspace{2mm}
The next lemma proves the existence of many examples of imprimitive $p$-exceptional linear groups, giving partial converse statements to Theorems \ref{imprim} and \ref{gen}.

\begin{lemma}\label{wreath}
Let $V_1 = \F_p^k$, let $n$ be a positive integer, and let $V = V_1^n = \F_p^{kn}$. Suppose $G_1 \le GL(V_1)$ and $H \le \SSS_n$ are such that one of the following conditions holds:

{\rm (i)} $G_1$ is transitive on $V_1^\sharp$ and $H$ is $p$-concealed,

{\rm (ii)} $G_1$ is $p$-exceptional and $H$ is a $p'$-group.

\no Then the wreath product $G =G_1 \hbox{ wr }H$, acting naturally on $V$, is $p$-exceptional.
\end{lemma}

\pf Suppose (i) holds, and let $0\ne v = (v_1,\ldots ,v_n) \in V_1^n = V$. Let $i_1,\ldots ,i_k$ be the positions $i$ for which $v_i \ne 0$. Then the orbit $v^G$ has size $|V_1^\sharp |^k \cdot \d$, where $\d$ is the size of the orbit of $H$ on $k$-sets containing $\{i_1,\ldots ,i_k\}$. As $H$ is $p$-concealed, $p$ does not divide $\d$, and so $|v^G|$ is coprime to $p$. The argument for (ii) is similar. \hal

\vspace{2mm}
We shall need the following upper bounds on the order of $p'$-subgroups
of $GL_m(q)$ for $m = 2,3$.

\begin{lemma}\label{small}
Let $q=p^f \geq 4$ and let $A$ be a $p'$-subgroup of $GL_m(q)$.

{\rm (i)} If $m = 2$ and $q \neq 5,7,11$ or $19$,  
then $|A| \leq (q^2-1) \cdot (2,q-1)$. 

{\rm (ii)} If $m = 3$ then $|A| \leq (q-1)(q^3-1)$.  
\end{lemma}

\pf
(i) It suffices to show that any $p'$-subgroup of $PGL_2(q)$ ($q \ne 5,7,11,19$) has order at most $(q+1)\cdot
(2,q-1)$. From the list of subgroups of $PGL_2(q)$ in \cite[Chapter XII]{Dickson}, any subgroup of order at least 
$(q+1)\cdot (2,q-1)$ has order dividing one of $2(q+1)$, $q(q-1)$, 24 (if $q=4$ or 9), or 60 (if $q=4,9$ or 29). The assertion follows.

(ii) The bound can be checked directly using \cite{Atlas} for $q \leq 11$, so we will
assume $q \geq 13$. If $A$ is reducible on $\F_q^3$, then $A$ is contained in
a maximal parabolic subgroup $P$ of $GL_3(q)$, and so 
$|A| \leq |P|_{p'} = (q^2-1)(q-1)^2$. If $A$ is irreducible but imprimitive, then
$|A| \leq 6(q-1)^3$. Finally, if $A$ is irreducible and primitive, then
$|A| \leq q^3 \cdot \log_2q^3$ by the main result of \cite{GW}. In all cases $|A| < (q-1)(q^3-1)$ since $q \geq 13$. 
\hal

\vspace{2mm}
Now we consider the case of $p$-exceptional 1-dimensional semilinear groups. Here we identify $V=\F_p^d$ with the field $\F_{p^d}$. Let $\o$ be a primitive element of $\F_{p^d}$ and let $\vp:x\rightarrow x^p$ denote the Frobenius automorphism. The 1-dimensional semilinear groups are subgroups of $\Ga L_1(p^d) = \la \hat\o,\vp\ra$, where $\hat\o$ denotes the multiplication map $x\mapsto x\o$. We determine all such $p$-exceptional groups and show that $p$ divides $d$ and there exists a unique minimal example.

\begin{lemma}\label{lem:1dim}
Suppose that $H\leq  \Ga L_1(p^d)$ and $H$ is $p$-exceptional on $V=\F_{p^d}$. Then $p$ divides $d$ and there is a factorisation $d=p^ks$ for some $k\geq1$ such that $H$ has a  normal subgroup $K=\la \hat\o^{(p^s-1)/j},\vp^s\ra$ of index coprime to $p$, for some $j$ dividing $p^s-1$. Moreover all such subgroups $H$ and $K$ are $p$-exceptional and each contains the $p$-exceptional group $\la \hat\o^{p^{d/p}-1},\vp^{d/p}\ra$. The group $K$ is $\frac{1}{2}$-transitive on $V^\sharp$, having $\frac{p^s-1}{j}$ orbits of length $\frac{p^d-1}{p^s-1}\cdot j$.
\end{lemma}

\pf  Write $H_0=H\cap\la\hat\o\ra = \la\hat\o^c\ra$, say, where $c$ divides $p^d-1$. 
Then $|H_0|=(p^d-1)/c$ is coprime to $p$, and $H/H_0\cong H\la\hat\o\ra/\la\hat\o
\ra$ is isomorphic to a subgroup of $\la\vp\ra$ and hence is cyclic of order dividing $d$.
Since $p$ divides $|H|$ it follows that $p$ divides $d$. Let $d=p^ks$ where $p^k$ is the $p$-part of $|H|$. Then $H$ has a unique normal subgroup $K$ containing $H_0$ such that $|K/H_0|=p^k$. The group $K$ is generated by $\hat\o^c$ and some element $\tau$ of the form $\vp^s\hat\o^b$. We may assume that $|\tau|=p^k$. A routine computation shows that $\tau^{p^k}=\hat\o^{b(p^d-1)/(p^s-1)}$, and hence $p^s-1$ divides $b$, say $b=(p^s-1)b'$. 

By Lemma~\ref{normalsubgroups}, $K$ is $p$-exceptional. This implies in particular that $\tau$ fixes setwise each of the $H_0$-orbits in $V^\sharp$, and these orbits are
the multiplicative cosets $\o^i\la\o^c\ra$ for $0\leq i\leq c-1$. Now $\tau$ maps $\o^i$ to $\o^{ip^s+b}$ and this element must therefore lie in $\o^i\la\o^c\ra$. It follows that
$i(p^s-1)+b=(i+b')(p^s-1)$ is divisible by $c$. Choosing $i$ such that $i+b'\equiv 1\pmod{c}$, we conclude that $c$ divides $p^s-1$, say $c=(p^s-1)/j$. This means that $\hat\o^b\in H_0$, and hence that $K=\la\hat\o^{(p^s-1)/j},\vp^s\ra$.

The computation in the previous paragraph shows that $\vp^s$ fixes each $H_0$-orbit setwise and hence that $K$ has $c$ orbits of length $(p^d-1)/c$ on nonzero elements of $V$. In particular $K$ is $p$-exceptional and hence any subgroup $H$ containing $K$ with index coprime to $p$, and intersecting $\la\hat\o\ra$ in $H_0$ is also $p$-exceptional. Finally each of these subgroups $K$ contains the group $\la \hat\o^{p^{d/p}-1},\vp^{d/p}\ra$, and our arguments (with $k=1$) show that this group is $p$-exceptional. 
\hal

\vspace{4mm}
Next we analyse the possibilities for 2-dimensional semilinear $p$-exceptional groups. We use the following notation: $Z$ denotes the group of scalar matrices in $GL_2(p^f)$; the group of diagonal matrices is denoted by $T$; and $\vp$ denotes the Frobenius map $(a_{ij})\mapsto (a_{ij}^p)$.

\begin{lemma}\label{lem:2dim}
Suppose that $H\leq  \Ga L_2(p^f)$ and $H$ is $p$-exceptional on $V=\F_{p^f}^2$. Then one of the following holds.
\begin{enumerate}
 \item[{\rm (i)}] $H$ contains $SL_2(p^f)$. 
\item[{\rm (ii)}] $p$ divides $2f$, $H\cap GL_2(p^f)$ is contained in $\la T, \begin{pmatrix} 0&1 \\ -1 &0\end{pmatrix} \ra$
if $p$ is odd, or $T$ if $p=2$.
\item[{\rm (iii)}] $p$ divides $2f$, and $H\leq\Ga L_1(p^{2f})$ is as in Lemma~{\rm \ref{lem:1dim}}.
\item[{\rm (iv)}] $p^f=9$ and $SL_2(5)\trianglelefteq H\cap GL_2(9)$.
\end{enumerate}
\end{lemma}

\pf If $H$ contains $SL_2(p^f)$ then $H$ is transitive on $V^\sharp$ so $H$ is $p$-exceptional. Suppose now that this is not the case, and let
$H_0=H\cap GL_2(p^f)$. 
 
Observe that, for a proper subfield $\F_{p^c}$ of $\F_{p^f}$ the group $SL_2(p^c)$ acts regularly on the orbit containing $(1,\o)$ where $\o$ is a primitive element of $\F_{p^f}$. If $SL_2(p^c)$ were normal in $H$ then the $H$-orbit containing $(1,\o)$ would have length a multiple of $p$. Thus $H$ has no normal subgroup conjugate to $SL_2(p^c)$ for any proper divisor $c$ of $f$. Moreover if $f\geq 3c$, then the
stabiliser in $Z\circ SL_2(p^c)$ of the 1-space $\la (1,\o)\ra$ is $Z$. Hence, if $p$ divides $|H_0|$ then $H_0$ is not contained in a conjugate of $Z\circ SL_2(p^c)$ for any divisor $c$ of $f$ such that $f\geq 3c$.

Observe that the $T$-orbits in $V$ have lengths $1, q-1, q-1, (q-1)^2$. Thus if $p$ divides $f$ then any subgroup of 
$T.\la\vp\ra$  containing $T$ and of order divisible by $p$ is $p$-exceptional. If $p$ is odd the same is true for such subgroups of $T.2.\la\vp\ra$. These examples and some of their subgroups are listed in (ii). So suppose that $H_0$ is not conjugate to a subgroup of $T.2$. 

Also if $H_0$ preserves on $V$ the structure of a 1-dimensional space over $\F_{p^{2f}}$, then $H$ is a 1-dimensional semilinear group and we obtain the examples in (iii) by Lemma~\ref{lem:1dim}. 

If $H_0$ has a non-trivial normal $p$-subgroup $K$ then for a vector $v$ not fixed by $K$, the $H$-orbit containing $v$ is a union of some $K$-orbits, each of length a nontrivial power of $p$. Thus no such subgroup exists. 

Consider $\bar H_0\cong H_0Z/Z\leq PGL_2(p^f)$. From our arguments so far, and the classification of subgroups of $PGL_2(p^f)$ \cite[Chapter XII]{Dickson}, we may assume that $\bar H_0\cong \AAA_4, \SSS_4$ (with $p$ odd) or $\AAA_5$ (with $p^f\equiv \pm 1 \hbox{ mod }10$), and that $H_0$ is not realisable  modulo scalars over a proper subfield $\F_{p^c}$ with $f\geq 3c$. In particular then, $p$ is odd and $f$ divides 4. Thus $p$ divides $|H_0|$ and hence $p=3$ (as $p\ne 5$ if $\bar H_0=\AAA_5$). 
If $\bar H_0 = \AAA_4$ or $\SSS_4$ then $H_0 \triangleright SL_2(3)$ which is not the case, so $\bar H_0 = \AAA_5$, and $p^f = 9$ or 81. In the latter case one checks that the orbit of $H_0$ containing the vector $(1,\o)$ has size divisible by 3. Hence $p^f = 9$, which leads to the examples in (iv) since $Z\circ SL_2(5)$ is transitive on the nonzero vectors. \hal

\vspace{2mm}
The next two lemmas concern the usual action of a group $G$ on a quotient group $G/V$ defined by $(Vx)^g = Vx^g$.

\begin{lemma}\label{vnorm}
Let $G$ be a finite group, $p$ a prime, and suppose $G$ has a normal subgroup $V$ which is an elementary abelian $p$-group. If $t\in G$ is a $p'$-element, then $C_G(t)/C_V(t) \cong C_{G/V}(t)$.
\end{lemma}

\pf Assume first that $[V,t] = V$. Let $g\in G$ be a preimage of an element of $ C_{G/V}(t)$, so that $t^g = tv$ for some $v\in V$. By assumption there exists $u\in V$ such that $[t,u]=v$. Then $t^u = tv = t^g$, so $g \in VC_G(t)$. This shows that 
$ C_{G/V}(t) = VC_G(t)/V$, as required.

Now consider the general case. Writing $V$ additively, we have $V = [V,t] \oplus C_V(t)$, by coprime action. Again let $g \in G$ with $t^g = tv$, $v\in V$, and write $v = v_1+v_2$ with $v_1 \in [V,t], v_2 \in C_V(t)$. If $v_2 \ne 0$ then $tv = t(v_1+v_2)$ has order divisible by $p$, a contradiction (as $t$ is a $p'$-element). Hence $v_2=0$, and now we argue as in the first paragraph of the proof. \hal

\begin{cor}\label{vnormgen}
Let $G$ be a finite group, $p$ a prime, and suppose $G$ has a normal $p$-subgroup $V$ such that $V/Z(V)$ is elementary abelian and $Z(V) \le Z(G)$. If $t\in G$ is a $p'$-element, then $C_G(t)/C_V(t) \cong C_{G/V}(t)$.
\end{cor}

\pf Write $\bar G = G/Z(V)$. By Lemma \ref{vnorm}, $\bar V C_{\bar G}(t)/\bar V = C_{\bar G/\bar V}(t)$. If $g \in G$ is a preimage of an element of $C_{\bar G}(t)$ then $t^g = tz$ for some $z \in Z(V)$. Since this has $p'$ order, we must have $z=1$, and the conclusion follows. \hal

\section{Imprimitive groups}\label{c2sec}

In this section we prove Theorems \ref{p-conc} and \ref{imprim}.
First, for Theorem \ref{p-conc},  we determine the primitive {\it $p$-concealed} groups, i.e. 
primitive subgroups $H$ of $\SSS_n$, such that the prime $p$ divides $|H|$ and every orbit
of $H$ on the set of all subsets of $\{1,\ldots ,n\}$ 
has length coprime to $p$.

\vspace{4mm}
\no {\bf Proof of Theorem \ref{p-conc}}

\vspace{2mm}
Let $\Om = \{1,\ldots ,n\}$, and define 
$\Om_k := \{ X \subseteq \Om \mid |X| = k\}$ for $0 \leq k \leq n$. Let $H$ be a primitive subgroup 
of $\SSS_n$ of order divisible by a prime $p$.

Assume first that $H \geq \AAA_n$. Since $p$ divides $|H|$, we have 
$n \geq p$ and $H \ne \AAA_3$ if $(n,p) = (3,2)$, and $H$ has exactly one orbit on $\Om_{p-1}$. Furthermore, 
$p$ is coprime to $|\Om_{p-1}|$ precisely when $p$ does not 
divide any of the 
$p-1$ consecutive integers $n-p+2, n-p+3, \ldots ,n$, that is, $p|(n+1)$.
Now we can write $n = \sum^{s}_{i=1}a_ip^i -1$, where $s \geq 1$, 
$a_s > 0$, and 
$p-1 \geq a_i \geq 0$. Suppose that $n \neq a_sp^s-1$. Choosing 
$k := p^s-1$, we see that 
$$\left\lfloor \frac{n}{p^s} \right\rfloor - \left\lfloor \frac{k}{p^s} \right\rfloor - 
   \left\lfloor \frac{n-k}{p^s} \right\rfloor = a_s - 0 - (a_s-1) = 1,$$
and so $p$ divides $|\Om_k|$. Next suppose that $n = a_sp^s-1$. 
Write any $\ell$ between $0$ and $n$ as $\ell =\sum^{s}_{i=0}b_ip^i$ with 
$0 \leq b_i \leq p-1$.
Then $b_s \leq a_s-1$ and $n-\ell = (a_s-b_s-1)p^{s}+
\sum^{s-1}_{i=0}(p-b_i-1)p^i$. Hence, for $0\le r \le s$,
$$\left\lfloor \frac{n}{p^r} \right\rfloor - \left\lfloor \frac{\ell}{p^r} \right\rfloor - 
   \left\lfloor \frac{n-\ell}{p^r} \right\rfloor = (a_sp^{s-r}-1) - \sum^{s}_{i=r}b_ip^{i-r}
   - (a_s-b_s-1)p^{s-r} - \sum^{s-1}_{i=r}(p-b_i-1)p^{i-r} = 0,$$
and so $p$ does not divide $|\Om_\ell|$. Since
$n \geq 3$, $H$ is transitive on $\Om_\ell$. We have shown that $H$ is 
$p$-concealed if $n = a_sp^s-1$.

From now on we will assume that $H \not\geq \AAA_n$. Clearly, if $H$ 
contains a normal subgroup $K$ which is also 
primitive of order divisible by $p$ 
and has a regular orbit $\Delta$ on $2^{\Om}$, then the 
$H$-orbit containing 
$\Delta$ has length divisible by $p$. Hence, we may assume that 
$H$ has no regular orbit on $2^{\Om}$ and apply \cite[Theorem 2]{Ser} 
to $H$. In all but three of the cases 
listed in \cite[Theorem 2]{Ser} for $H$, either we can find a subgroup 
$K$ with
the prescribed properties, or we can use \textsf{GAP} \cite{gap} or \textsc{Magma} \cite{Mag}
to show directly that $H$ is not a $p$-concealed group. The three exceptional cases give the examples in parts (ii) and (iii) of Theorem \ref{p-conc}. \hal

\vspace{4mm}
The main result of this section is the following:

\begin{thm}\label{c2m}
Assume $G < GL(V)$ is a (not necessarily 
irreducible) $p$-exceptional group which acts primitively 
as a permutation group on the $n$ summands of the direct sum decomposition
$$\Fq^{mn} = V = V_1 \oplus V_2 \oplus \ldots \oplus V_n,$$
where $\dim_{\Fq}V_i = m \geq 1$ and $n \geq 2$. Let $H \leq \SSS_n$ 
be the subgroup induced by this primitive action. Then one of the 
following holds:

{\rm (i)} $H$ is a $p'$-group;

{\rm (ii)} One of {\rm (i)}, {\rm (ii)} and {\rm (iii)} of 
Theorem $\ref{p-conc}$ holds for $H$.
Moreover, $G_{V_1}$ is transitive on $V_1^\sharp$. 
\end{thm}

\pf 
Assume that $p$ divides $|H|$. First we show that $H$ is  
$p$-concealed, and so
(i), (ii) or (iii) of Theorem \ref{p-conc} holds for $H$. 
Indeed, suppose that 
an $H$-orbit $\Delta$ of $H$ on the subsets of 
$\Om := \{V_1, \ldots ,V_n\}$ has
length divisible by $p$. Pick $X = \{V_{j_1}, \ldots ,V_{j_t}\} \in \Delta$,
$0 \neq v_i \in V_{j_i}$ for $1 \leq i \leq t$, and let 
$v := v_{j_1} + \ldots + v_{j_t}$. Then $I := G_v$ preserves $X$ 
and so 
$IK/K \leq H_{X}$ for $K := \cap^{n}_{r=1}G_{V_r}$. 
Since $p$ divides 
$|\Delta|$, it must divide $[G:I] = |v^G|$, contrary to the 
$p$-exceptionality of $G$.

Now we assume that $H$ satisfies one of the conclusions of Theorem \ref{p-conc}, 
but $G_1 := G_{V_1}$ has at least two orbits $u^{G_1}$ 
and $v^{G_1}$ on $V_1^\sharp$. For each $1 \leq i \leq n$, fix 
some $g_i \in G$ such that $V_1g_i = V_i$ (and $g_1 = 1$), and set 
$u_i = ug_i$, $v_i = vg_i$. Choose any non-empty subsets
$X, Y \subset \Om$ with $X \cap Y = \emptyset$. Also consider
$w = \sum^n_{i=1}w_i$, with $w_i = u_i$ if $V_i \in X$, $w_i = v_i$ 
if $V_i \in Y$,
and $w_i = 0$ otherwise. Observe that any $h \in G_w$ 
stabilises both 
$X$ and $Y$. (Indeed, $h$ fixes $X \cup Y$. Assume that 
$V_ih = V_j \in Y$ for some $V_i \in X$. Comparing the 
$V_j$-component of $wh = w$, we see that 
$vg_j = v_j = u_ih = ug_ih$, and so $ug_ihg_j^{-1} = v$. But
$g_ihg_j^{-1}$ stabilises $V_1$, so we conclude that $v \in u^{G_1}$, 
a contradiction.) 
Now the $p$-exceptionality of $G$ implies that $p$ does not divide 
$|H:J|$ for $J := H_{X,Y}$, the subgroup of $H$ consisting of all elements
that stabilise $X$ and stabilise $Y$. 

It remains to exhibit a pair $(X,Y)$ such that $p$ divides $|H:J|$ 
to get the 
desired contradiction. In the case (i) of Theorem \ref{p-conc}, we choose
$X := \{V_1\}$, $Y := \{V_2, \ldots ,V_{p}\}$. Then 
$$|H:J| = |\SSS_n:(\SSS_{p-1} \times \SSS_{n-p})| = p \cdot \frac{n!}{p! \cdot (n-p)!}$$
is divisible by $p$. In the case (ii) of Theorem \ref{p-conc}, we can 
choose $X := \{a\}$ and $Y := \{b,c\}$ for some distinct 
$a,b,c \in \Omega$, 
and check that $|H:J| = 168$ is divisible by $p = 3$. 
In the case (iii) of Theorem \ref{p-conc}, we can 
choose $X := \{a\}$ and $Y := \{b,c\}$ for some distinct 
$a,b,c \in \Omega$, 
and check that $|H:J| = 10$ is divisible by $p = 2$.              
\hal

\vspace{4mm}
\no {\bf Deduction of Theorem \ref{imprim}}

\vspace{2mm}
Theorem \ref{imprim} follows very quickly from the above theorem. Indeed,
suppose $G \le GL(V)$ is irreducible, imprimitive and $p$-exceptional, with
$G = O^{p'}(G)$. Let 
$V = W_1 \oplus ... \oplus W_m$ be an imprimitivity decomposition for 
$G$. 
Coarsen this to a decomposition $V = V_1 \oplus ... \oplus V_n$ such that
$G$ acts as a primitive permutation group on $\{V_1,\ldots ,V_n\}$,
where $V_1$ is, say, $W_1 \oplus ... \oplus W_k$.
As $G = O^{p'}(G)$, conclusion (ii) of Theorem \ref{c2m} holds for the 
action of $G$ on $\{V_1,\ldots ,V_n\}$. In particular, 
$G_{V_1}$ acts transitively on the nonzero vectors of
$V_1$. But it also permutes $W_1, ..., W_k$, so $k = 1$. Theorem \ref{imprim}
follows.

\section{Tensor products I: $\mathcal{C}_4$ case }\label{c4sec}

In this section we handle $p$-exceptional groups preserving tensor product decompositions. These correspond 
to subgroups of groups in class $\mathcal{C}_4$ in Aschbacher's classification of maximal subgroups of classical groups \cite{asch},
hence the title of this (and forthcoming) sections.
If $U$ and $W$ are vector spaces over a field $\F_q$, then a central product $GL(U) \circ GL(W)$ acts
naturally on $V = U \otimes W$. We also denote by ${\Ga L(V)}_{U \otimes W}$ the stabiliser 
of the tensor decomposition, which is a group $(GL(U) \circ GL(W))\la \s \ra$ where $\s$ is a field automorphism
fixing both factors.
As usual, write $Z$ for the group $\F_q^*$ of scalars in $GL(V)$.

\begin{thm}\label{c4lin}
Let $V$ be a vector space over $\Fq$ of characteristic $p$, and write $V = U\otimes W$, a tensor product over $\Fq$ with $\dim U, \dim W \ge 2$.
Let $H \le GL(U) \circ GL(W) < GL(V)$ and suppose $H$ is $p$-exceptional.

Then $p=2$, $\dim U = \dim W = 2$, and $ZH = (GL_1(q^2)\circ GL_1(q^2)).2$, where involutions in $H$ act nontrivially 
on both factors $GL_1(q^2)$. This group is $p$-exceptional, and acts reducibly on~$V$.
\end{thm}

We also need a result for the semilinear case.

\begin{thm}\label{c4semi}
Let $V$ be a vector  space over $\F_{q^p}$ (of characteristic $p$),  and write $V = U \otimes W$ with $2 \le \dim U \le \dim W$. Let $H \le {\Gamma L(V)}_{U \otimes W}$, and assume $H$ is $p$-exceptional and $H \cap GL(V)$ is a $p'$-group.

Then $p = \dim U = 2$. In particular, $H \cap GL(V)$ is not absolutely irreducible on~$V$.
\end{thm}

The proofs are given in the following three subsections.

\subsection{Some theory of tensor decompositions}\label{subtensorsetup}
First we give some general theory for tensor decompositions $V=U\otimes W$ of a vector space $V=\Fq^n$, where $q=p^f$ for a prime $p$, $a=\dim U\geq2, b=\dim W\geq2$, and $n=ab$.

Let $\{u_1,\dots,u_a\}$ be a basis for $U$, $\{w_1,\dots,w_b\}$ be a basis for $W$, and write elements of $GL(U), GL(W)$ as matrices with respect to these bases respectively. Then $V$ has an associated basis $\mathcal{B} :=\{u_i\otimes w_j\,|\,1\leq i\leq a, 1\leq j\leq b\}$, which we refer to as the \emph{standard basis}.
For elements $u=\sum_ia_iu_i\in U$ and $w=\sum_jb_jw_j\in W$ we denote the element $\sum_{i,j}(a_ib_j)(u_i\otimes w_j)$ of $V$ by $u\otimes w$.
A vector $v\in V$ is called \emph{simple}  if it can be expressed as $v=u\otimes w$ for some $u\in U, w\in W$.

The stabiliser $X:={GL(V)}_{U \otimes W}$ in $GL(V)$ of this decomposition is a central product of $X:=GL(U)\circ GL(W)$ and we view elements
of $X$ as   ordered pairs $(A,B)\in GL(U)\times GL(W)$ modulo the normal 
subgroup
$Z_0=\{(\l I,\l^{-1}I)\,|\,\l\in\F_q^*\}$, where $(A,B):u_i\otimes w_j\mapsto u_iA\otimes w_jB$ (extending linearly).
The stabiliser $\hat X:={\Gamma L(V)}_{U \otimes W}$ in $\Gamma L(V)$ is a semidirect product of $X$ and the group $\la \sigma\ra$ of field automorphisms, where $\sigma:\sum_{ij}a_{ij}u_i\otimes w_j\mapsto \sum_{ij}a_{ij}^pu_i\otimes w_j$.

For an arbitrary $v\in V$, the
\emph{weight of $v$} is defined as the minimum number $k$ such that $v$ can be written as a sum of $k$ simple vectors. It is not difficult to prove that elements of $\hat X$ map weight $k$ vectors to weight $k$ vectors; in particular, the notion of a simple vector does not depend on the choice of standard basis. Also, the weight of a vector is well defined since any vector can be written as a sum of $n$ simple vectors (each a scalar multiple of an element of $\mathcal{B}$).

For subspaces $U_0$ of $U$ and $W_0$ of $W$ (not necessarily proper subspaces), $X_{U_0\otimes W_0}$ consists of all elements $(A,B)$ of $X$ such that $U_0$ is $A$-invariant and $W_0$ is $B$-invariant, and $\hat X_{U_0\otimes W_0}$ is generated by $X_{U_0\otimes W_0}$ and a conjugate of $\sigma$.

\begin{lemma}\label{weight}
Let $v\in V$ of weight $k$, and suppose that $v=\sum_{i=1}^kx_i\otimes y_i$, where the $x_i\in U, y_i\in W$. Write $U_0=\la x_1,\dots, x_k \ra$ and $W_0=\la y_1,\dots, y_k\ra$. Then \begin{enumerate} 
\item[{\rm (i)}] $\dim(U_0)=\dim(W_0)=k$, so $k\leq\min\{a,b\}$.
\item[{\rm (ii)}] If also $v=\sum_{i=1}^kx_i'\otimes y_i'$, where the $x_i'\in U, y_i'\in W$, then $\la x_1',\dots, x_k'
\ra=U_0$ and $\la y_1',\dots, y_k'\ra=W_0$.
\item[{\rm (iii)}] Let $A, B$ be the matrices representing the linear transformations $A:x_i\mapsto x_i',\ B:y_i\mapsto y_i'$ (for all $i$) of $U_0, W_0$ with respect to the bases $\{x_1,\dots,x_k\}$ and $\{y_1,\dots,y_k\}$ respectively, where $x_i',y_i'$ are as in {\rm (ii)}. Then $B=A^{-T}\in GL_k(q)$.
\item[{\rm (iv)}] $\hat X_v\leq \hat X_{U_0\otimes W_0}$, the group induced by $X_v$ on $U_0\otimes W_0$ is a diagonal subgroup of the group $GL_k(q)\circ GL_k(q)$ induced by $X_{U_0\otimes W_0}$, consisting of the pairs $(A,A^{-T})$ (modulo $Z_0$) for $A\in GL_k(q)$ (with respect to the bases in (iii)), and $\hat X_v=\la X_v,\sigma'\ra$, where $\sigma'$ is conjugate to $\sigma$ and induces a generator of the group of field automorphisms of $GL(U_0\otimes W_0)$.
\end{enumerate}
\end{lemma}

\pf
Part (i) follows almost from the definition of $k$, since, for example, if $x_k=\sum_{i=1}^{k-1}a_ix_i$, then $v= \sum_{i=1}^{k-1}x_i \otimes(y_i+a_iy_k)$. For the rest of the proof we will assume without loss of generality that $x_i=u_i$ and $y_i=w_i$ for $1\leq i\leq k$.

Suppose also that $v$ has an expression as in part (ii).
For $1\leq i\leq k$, in terms of
the basis $\mathcal{B}$ we have $x_i'=\sum_{j=1}^a a_{ij}u_j$ and $y_i'=\sum_{j=1}^bb_{ij} w_j$. Write $A_0=(a_{ij})\in\M_{k\times a}(q)$ and $B_0=(b_{ij})\in \M_{k\times b}(q)$. Then \[ v=\sum_{i=1}^k\sum_{j=1}^a\sum_{\ell=1}^b a_{ij}b_{i\ell}u_j\otimes w_\ell =\sum_{j,\ell}\left(\sum_{i=1}^k a_{ij}b_{i\ell}\right)u_j\otimes w_\ell = \sum_{j,\ell} (A_0^TB_0)_{j,\ell} u_j\otimes w_\ell.
\]
Then since $v=\sum_{i=1}^ku_i\otimes w_i$, and since $\mathcal{B}$ is a basis for $V$, we deduce that $(A_0^TB_0)_{j\ell}=0$ if at least one of $j>k, \ell>k, j\ne\ell$, and $(A_0^TB_0)_{jj}=1$ for $1\leq j\leq k$. In particular this implies that each $x_i'\in U_0$ and each $y_i'\in W_0$. Moreover by part (i) the $x_i'$ are linearly independent, and also  the $y_i'$ are linearly independent. Thus part (ii) follows. We have also proved part (iii).

Finally if $g=(A_1,B_1)\in X_v$ (modulo $Z_0$), with $A_1=(a_{ij})\in GL(U)$ and $B_1=(b_{ij}) \in GL(W)$, then we have $v=\sum_{i=1}^k(u_iA_1)\otimes(w_iB_1)$. By part (ii), for each $i\leq k$, $u_iA_1\in U_0$ and $w_iB_1\in W_0$, so $g\in X_{U_0\otimes W_0}$. 
Moreover, if $A=A_1|_{U_0},
B=B_1|_{W_0}$, written as matrices
with respect to the bases $\{x_1,\dots,x_k\}$ and $\{y_1,\dots,y_k\}$ respectively, then by part (iii), $B=A^{-T}$ so $g$ is one of the elements described in part (iv). For each $A\in GL_k(q)$, there exists an element $(A_1,B_1)\in X$ with $A=A_1|_{U_0}, A^{-T}=B_1|_{W_0}$, and the fact that this element fixes $v$ follows from the displayed computation above. 
Finally $\hat X_v$ contains a conjugate of $\sigma$ which induces on $U_0\otimes W_0$ the natural field automorphism with respect to the basis formed by the $x_i\otimes y_j$. \hal

\vspace{2mm}
By Lemma~\ref{weight}(ii), the subspaces $U_0$ and $W_0$ are determined uniquely by $v$, and we denote them by $U_0(v)$ and $W_0(v)$ respectively.

\subsection{Proof of Theorem~\ref{c4lin}}

Suppose $H \le \hat X$ preserves a non-trivial tensor decomposition $V=U\otimes W$ of $V=\F_q^n$, and that $\hat X,$ $X, Z_0, q=p^f, a, b$ are as in Subsection~\ref{subtensorsetup}.
Suppose also that $H$ is $p$-exceptional on $V$.
By Lemma~\ref{em:related1} we may assume that $H$ contains $Z:=Z(GL(V))=(Z_U\times Z_W)/Z_0$, where $Z_U=Z(GL(U)), Z_W= Z(GL(W))$.

\medskip
For Theorem~\ref{c4lin} we will have $H\leq X$, but Lemma \ref{firstfacts} is more general and will be used also for the proof of Theorem~\ref{c4semi}. The natural projection maps $\phi_U:\hat X\rightarrow P\Ga L(U)$ and $\phi_W:\hat X\rightarrow P\Ga L(W)$ have kernels $K_U=Z_U\circ GL(W)\cong GL(W)$ and $K_W=GL(U)\circ Z_W\cong GL(U)$ respectively. Also for subspaces $U_0\leq U, W_0\leq W$, we have maps
$\phi_{U_0}: \hat X_{U_0\otimes W}\rightarrow P\Ga L(U_0)$ and $\phi_{W_0}:\hat X_{U\otimes W_0} \rightarrow P\Ga L(W_0)$ with kernels $K_{U_0}, K_{W_0}$ respectively. If an element $x\in\hat X$ or subgroup $L\leq \hat X$ lies in $\hat X_{U_0\otimes W}$ then we write $x^{U_0}=\phi_{U_0}(x)$, $L^{U_0}=\phi_{U_0}(L)$ for the corresponding element or subgroup of $P\Ga L(U_0)$; also by the \emph{fixed point subspace} of $L^U$ we mean the largest $L^U$-invariant subspace $U_0$ of $U$ such that $L^{U_0}=1$. Similarly for subspaces of $W$.
A subgroup $L\leq \hat X_{U_0\otimes W_0}$ is said to \emph{act diagonally on} $U_0\otimes W_0$ if $L\cap K_{U_0}$ and $L\cap K_{W_0}$ both induce only scalar transformations on $U_0\otimes W_0$.

\begin{lemma}\label{firstfacts}
Suppose that $H\leq \hat X$ and $H$ is $p$-exceptional. Let $U_0, W_0$ be $2$-dimensional subspaces of $U, W$ respectively, and let $P$ be a Sylow $p$-subgroup of $H_{U_0\otimes W_0}$. Then $P$ is a Sylow $p$-subgroup of $H$ (and so $P\ne1$), $P$ acts  diagonally on $U_0\otimes W_0$, and moreover $P\cap K_{U_0}=P\cap K_{W_0}=1$, and $P\cap X$ is elementary abelian of order at most $q$.
\end{lemma}

\pf
Choose a weight 2 vector $v\in V$ such that $U_0(v)=U_0$ and $W_0(v)=W_0$. Then by Lemma \ref{weight}, $H_v\leq H_{U_0\otimes W_0}$, and by our assumption $|H:H_v|$ is coprime to $p$. Thus $P$ is a Sylow $p$-subgroup of $H$, and in particular $P\ne1$. Also,  $P$ is conjugate in $H_{U_0\otimes W_0}$ to a Sylow $p$-subgroup $P'$ of $H_v$. Since $P'$ induces a diagonal action on $U_0\otimes W_0$ by Lemma~\ref{weight}(iv), it follows that also $P$ induces a diagonal action on $U_0\otimes W_0$.

Let $Q=P\cap K_{U_0}$ and assume that $Q\ne1$. Since $P$ acts diagonally  on $U_0\otimes W_0$, it follows that $Q^{W_0}=1$ (since $Q$ is a $p$-group), and we deduce that $Q=P\cap K_{W_0}$.
Further, since $Q\ne1$, we may assume without loss of generality that $Q^U\ne1$. We produce an element $g\in Q$ and a 2-dimensional  $\la g^U\ra$-invariant subspace $U_0'$ such that $g^{U_0'}\ne1$ as follows: the fixed point subspace $U_1$ of the nontrivial $p$-group $Q^U$ contains $U_0$ and $U_1\ne U$, and there is a 1-dimensional subspace $U_2/U_1$ of $U/U_1$ left invariant by $Q^U$ in its induced quotient action on $U/U_1$.
Since $Q^{U_1}=1$ and $Q^{U_2}$ is a nontrivial $p$-group, there exists $u\in U_2\setminus U_1$ and $g\in Q$ such that $\la u^g\ra = \la u + x\ra $ for some non-zero $x\in U_1$. Let $U_0'= \la u,x\ra$, and note that $x^g=x$ since $x\in U_1$.
Thus $U_0'$ is invariant under $\la g^U\ra$ and has dimension 2, and $g^{U_0'}\ne1$. Moreover since $g\in Q\subset P$, it follows that $g\in H_{U_0'\otimes W_0}$ and that $g$ does not act diagonally on $U_0'\otimes W_0$ since $g^{W_0}=1$. This contradicts the diagonal action of a Sylow $p$-subgroup of $H_{U_0'\otimes W_0}$ on $U_0'\otimes W_0$.
Thus $P\cap K_{U_0}=P\cap K_{W_0}=1$. In particular $P$ is isomorphic to a subgroup of $P\Ga L(U_0)$ and hence $P\cap X$ is elementary abelian of order at most $q$. \hal

\begin{cor}\label{lem:dim4} Assume that $H\leq X$, and let $1\ne g\in H$ be a $p$-element. Then $g^U$ and $g^W$ are both regular unipotent of order $p$.
\end{cor}

\pf
We may assume that $1\ne g\in P\leq H_{U_0\otimes W_0}$ with $U_0, W_0, P$ as in Lemma~\ref{firstfacts}.
Then $P$ is elementary abelian so $|g|=p$, and $g^{U_0}, g^{W_0}$ are both nontrivial.
Also $g^U, g^W$ have non-zero fixed point subspaces.
Suppose that $U_1$ is a 2-dimensional subspace of the fixed point subspace of $g^U$.
Then $g$ is a nontrivial $p$-element in $H_{U_1\otimes W_0}$, contradicting the diagonal action of $p$-elements proved in Lemma~\ref{firstfacts}.  Thus the fixed point subspace of $g^U$ has dimension 1, so $g^U$ is regular unipotent.  
Similarly $g^W$ is regular unipotent. \hal

\vspace{2mm}
In the next lemma we use the following notation.
For a vector space $V$ of dimension at least $k$, let $P_k(V)$ denote the set of $k$-dimensional subspaces of $V$.
Also let $\Ga L^*_1(q^p) = GL_1(q^p).p$, where the cyclic group of order $p$ on top is generated by a field 
automorphism (so that $\Ga L^*_1(q^p) \le GL_p(q)$).

\begin{lemma}\label{lem:trans}
Assume that $H\leq X$ and let $k\leq \min\{\dim(U),\dim(W)\}$. Then $H$ acts transitively on $P_k(U)\times P_k(W)$ (in its natural product action).
Moreover $p=\dim(U)=\dim(W)\in\{2,3\}$, and $H\leq\Ga L_1^*(q^p)\circ\Ga L_1^*(q^p)$.
\end{lemma}

\pf
By Corollary~\ref{lem:dim4}, each element of order $p$ in $H$ fixes unique $k$-dimensional subspaces of $U$ and of $W$, and this property is true also for each Sylow $p$-subgroup of $H$.
For $i=1,2$, let $U_i, W_i$ be a $k$-dimensional subspace of $U$, $W$ respectively, and let $v_i=\sum_{j=1}^k x_{ij}\otimes y_{ij}\in U_i\otimes W_i$ be a weight $k$-vector so that the $x_{ij}$ span $U_i$ and the $y_{ij}$ span $W_i$. Let $P_i$ be a Sylow $p$-subgroup of $H_{v_i}$. Since $H$ is $p$-exceptional, $P_i$ is a Sylow $p$-subgroup of $H$, and by Lemma~\ref{weight}, $H_{v_i}\leq H_{U_i\otimes W_i}$. Thus there is an element $x\in H$ such that $P_1^x=P_2$, and hence $P_2$ fixes the $k$-subspaces $U_1^x, U_2$ of $U$ and $W_1^x, W_2$ of $W$. By uniqueness, we have $U_1^x=U_2$ and $W_1^x=W_2$. This proves the first assertion.

By \cite[Lemma 4.1]{regsgps}, for the group $\phi_U(H)$ to be transitive on $k$-subspaces, one of the following holds: (i) $\phi_U(H)\geq PSL(U)$, or (ii) $k\in\{1,\dim(U)-1\}$ and $\phi_U(H)\not\geq PSL(U)$, or (iii) $\dim(U)=5, p=2$ and $|\phi_U(H)|=31.5$. Since $p$ divides $|\phi_U(H)|$, case (iii) does not arise, and by Corollary~\ref{lem:dim4}, $\dim(U)=2$ in case (i).
The same comments apply to $\phi_W(H)$. Since we may always take $k=2$ in the previous paragraph, it follows that each of $a:=\dim(U)$ and $b:=\dim(W)$ is at most 3. Then by the classification of transitive linear groups (see \cite[Appendix]{liebaff}), and noting that $p$ divides $|\phi_U(H)|$, either $\phi_U(H)\leq \Ga L_1(q^a)/Z_U$, or $a=2$ and $\phi_U(H)\geq PSL(U)$, or $\AAA_5\leq\phi_U(H)\leq P\Ga L_2(9)$. We have the same three possibilities for $\phi_W(H)$.

Suppose first that $\phi_U(H)\geq PSL(U)$ with $a=2$. Then since $P$ acts diagonally by Lemma~\ref{firstfacts}, it follows that $b=2$ and $H$ has a single composition factor $PSL_2(q)$. However in this case $H$ is not transitive on $P_1(U)\times P_1(W)$. A similar argument rules out the third possibility.
Thus  $\phi_U(H)\leq \Ga L_1(q^a)/Z_U$, and similarly $\phi_W(H)\leq \Ga L_1(q^b)/Z_W$.
Since $p$ divides the order of each of these groups we must have $p=a=b\in\{2,3\}$.
\hal

\vspace{2mm}
Next we show that the case $p=3$ does not yield a $p$-exceptional group.
The proof and the proof of Lemma \ref{lem:podd} use the following simple fact.

\begin{remark}\label{rem:weight}{\rm
Suppose that $a=\dim(U)=\dim(W)$.  Let $\{u_1, u_2,\dots, u_a\}$ be a basis for $U$.
By Lemma~\ref{weight}, each weight $a$ vector $v$ in $U\otimes W$ has a unique representation of the form $v=\sum_{i=1}^a u_i\otimes w_i$ where the $w_i$ form a basis for $W$. Thus the number of weight $a$ vectors is $|GL_a(q)|$.
}
\end{remark}

\begin{lemma}\label{not3}
If $p=3$, and $H\leq\Ga L_1^*(q^p)\circ\Ga L_1^*(q^p)$, then $H$ is not $p$-exception\-al.
\end{lemma}

\pf
Suppose that $H$ is as stated and that $H$ is $p$-exceptional. By Lemma~\ref{em:related1}, we may assume that $Z<H$, and by Lemma~\ref{firstfacts} it follows that a Sylow $3$-subgroup $P$ of $H$ has order 3 and acts diagonally on $U\otimes W$.
Let $r:=q^2+q+1$. By Lemma~\ref{lem:trans}, $r^2$ divides $|H|$. Also, since $\gcd(r,q-1)=1$, it follows that $H=Z\times (Z_r^2.P)$. Now $P$ centralises $Z$ and $N_H(P)=ZP$ so that $H$ has exactly $r^2$ Sylow $3$-subgroups.

The group $P$ in its action on $U$ leaves invariant a unique 2-subspace $U_2$. 
The same is true for the
$P$-action on $W$. Let $\{u_1, u_2, u_3\}$ be a $P$-orbit forming a basis for $U$.
By Lemma~\ref{weight}, each weight 3 vector $v$ in $U\otimes W$ has a unique representation of the form
$v=\sum_{i=1}^3 u_i\otimes w_i$ where the $w_i$ form a basis for $W$, and it is straightforward to show that $v$ is fixed by $P$ if and only if the $w_i$ form a $P$-orbit in $W$; now each such $P$-orbit yields three weight 3 vectors fixed by $P$. Hence the number of weight 3 vectors fixed by $P$ is exactly $q^3-q^2$. By Remark~\ref{rem:weight}, the number of weight 3 vectors is $|GL_3(q)|$, and since each weight 3 vector is fixed by some Sylow $3$-subgroup of $H$, it follows that $H$ has at least $|GL_3(q)|/(q^3-q^2)$ Sylow $3$-subgroups.  Since $r^2 < 
|GL_3(q)|/(q^3-q^2)$, this is a contradiction.
\hal

\vspace{2mm}
Finally we show that the case $p=2$ does lead to (reducible) $p$-exceptional examples.
This result completes the proof of Theorem~\ref{c4lin}.

\begin{lemma}\label{yes2}
If $p=2$, and $H\leq\Ga L_1^*(q^2)\circ\Ga L_1^*(q^2)$, then $ZH$ is $2$-exceptional if and only if $ZH=(GL_1(q^2)\circ GL_1(q^2)).2$, an index $2$ subgroup of $\Ga L_1^*(q^2)\circ\Ga L_1^*(q^2)$, with the Sylow $2$-subgroups acting diagonally. This group $ZH$ is reducible and $2$-exceptional, with two orbits of length $q^2-1$ and $q-1$ orbits of length $(q^2-1)(q+1)$ on non-zero vectors.
\end{lemma}

\pf Suppose that $H$ is as stated and that $H$ is $2$-exceptional. By Lemma~\ref{em:related1}, we may assume also that $Z<H$. Arguing as in the proof of Lemma~\ref{not3}, a Sylow $2$-subgroup $P$ has order 2 and acts diagonally on $U\otimes W$, $N_H(P)=ZP$, and $H=Z\times (Z_{q+1}^2.P) =(GL_1(q^2)\circ GL_1(q^2)).2$.

Conversely suppose that $H=(GL_1(q^2)\circ GL_1(q^2)).2$ with each Sylow $2$-subgroup $P$ acting diagonally.
Identify $U$ and $W$ with $\F_{q^2}$ and let $H_1$ be the index 2 subgroup of $H$ so that $H_1$ acts by field multiplication on both factors. Then $H=H_1\langle \tau\rangle$ where $\tau$ acts on both factors as the field automorphism of order 2. Let $\omega\in\F_{q^2}$ have order $q+1$ and note that $\omega^2=1+\lambda\omega$ where $\lambda=\omega+\omega^q\in\F_q$.

Let $v=1\otimes 1+\omega\otimes\omega$ and $w=1\otimes\omega+\omega\otimes 1$. Let $X=GL_1(q^2)\circ I$ and note that $v^X=\{\xi\otimes 1+\xi\omega\otimes\omega\mid\xi\in\F_{q^2}^*\}$ is a set of size $q^2-1$ and forms the set of nonzero vectors of  an $\F_q$-subspace $U_1$ of $V$. Similarly, $w^X=\{\xi\otimes\omega+\xi\omega\otimes 1\mid\xi\in\F_{q^2}^*\}$ is also the set of $q^2-1$ nonzero vectors of an $\F_q$-subspace $U_2$. Moreover, $U_1\cap U_2=\{0\}$. Note that an element of $X$ induces multiplication by the same element of $\F_{q^2}^*$ on each of $U_1$ and $U_2$.

Let $g\in H$ be the element that multiplies by $\omega$ in both the first  and second factors. Then
\[
\begin{array}{ll}
 w^g &= \omega\otimes\omega^2+\omega^2\otimes\omega \\
     &  =\omega\otimes(1+\lambda\omega)+(1+\lambda\omega)\otimes\omega \\
     & =\omega\otimes 1+ (\lambda\omega+1+\lambda\omega)\otimes\omega=w
\end{array}
\]
Thus  $|w^{H_1}|\leq q^2-1$ but since $w^X\subseteq w^{H_1}$ it follows that $w^{H_1}$ has size $q^2-1$. Also 
\[
\begin{array}{ll}
(\xi\otimes \omega+\xi\omega\otimes 1)^\tau &= \xi^q\otimes\omega^q+\xi^q\omega^q\otimes 1\\
                                             &= \xi^q\otimes (\lambda+\omega)+\xi^q(\lambda+\omega)\otimes 1\\
                                             &=\xi^q\otimes\omega+ (\xi^q\lambda+\xi^q\lambda+\xi^q\omega)\otimes 1\\
                                             & =\xi^q\otimes\omega+\xi^q\omega\otimes 1\in U_2
\end{array}
\]
Hence $U_2$ is $H$-invariant.

Similar calculations, taking $g$ to be the element that multiplies the first factor by $\omega$ and the second by $\omega^q$, show that $U_1$ is also $H$-invariant. 

Let $C$ be the subgroup of $GL_1(q^2)\times GL_1(q^2)$ acting on $U_1\oplus U_2$ given by 
$$C := \{ (\alpha,\beta) \mid \alpha,\beta \in \F_{q^2}^*,~
  (\alpha\beta^{-1})^{q+1} = 1 \}.$$
Note that $|C|=(q^2-1)(q+1)=|H_1|$. We have already seen that elements of $X$ are elements of $C$. Now consider elements $y=(1,\xi)$ of $Y=I\circ GL_1(q^2)\le H_1$ and let $a,b\in\F_q$ such that $\xi=a+b\omega$. Then
\[
\begin{array}{ll}
 (1\otimes 1+\omega\otimes\omega)^y &=1\otimes(a+b\omega)+\omega\otimes (a+b\omega)\omega \\
                                   &=a\otimes 1+b\otimes\omega+a\omega\otimes\omega+b\omega\otimes  (1+\lambda\omega)\\
                                   &= (a+b\omega)\otimes 1+(b+a\omega+b\lambda\omega)\otimes \omega\\
                                   &= \xi\otimes 1+ \xi\omega\otimes\omega                                   
\end{array}
\]
and so $y$ induces multiplication by $\xi$ on $U_1$. Similarly, $(1\otimes\omega+\omega\otimes 1)^y=\xi^q\otimes\omega+\xi^q\omega\otimes 1$ and so $y$ induces multiplication by $\xi^q$ on $U_2$. Thus the elements of $Y$ are also elements of $C$. Since each element of $H_1$ is the product of an element of $X$ and an element of $Y$, comparing orders yields $C=H_1$. We have already seen that $H_1$ has two orbits of length $q^2-1$  ($U_1^\sharp$ and $U_2^\sharp$), and using the fact that $H_1=C$ we see that it also has
$q-1$ orbits of length $(q^2-1)(q+1)$, namely  
$$\Delta_{\lambda} := \{ (u,v) \mid u,v \in \F_{q^2}^*,~
  (uv^{-1})^{q+1} = \lambda \}$$
for $\lambda \in \F_q^*$, on $V^\sharp$. 
Also consider $s_1,s_2 \in GL(V)$, where $s_1$ sends $(u,v) \in V$ to $(u^q,v)$  
and $s_2$ sends $(u,v)$ to $(u,v^q)$. Then 
$\tilde{C} := \langle C,s_1,s_2 \rangle \cong C:2^2$ has the same orbits as 
$C$ does on $V$. Moreover, $H=C:\langle s_1s_2\rangle$.
 \hal

\subsection{Proof of Theorem \ref{c4semi}}

Suppose now that $H\leq \Ga L_{n}(q^p)$ acting on $V=\F_{q^p}^n$, is 
$p$-exceptional, 
that $p=|H:H\cap GL(V)|$, and that $|H\cap GL(V)|$ is coprime to $p$.  Suppose 
moreover that $H\leq X:={\Ga L(V)}_{U\otimes W}$, where 
$a=\dim(U)\geq2, b=\dim(W)\geq2$ with $a\leq b$, and set $r=q^p$. By 
Lemma~\ref{em:related1}, we may assume that $H$ contains $Z=Z(GL(V))$. 
Theorem \ref{c4semi} follows from the following lemma.

\begin{lemma}\label{lem:podd}
$p=a=2$.
\end{lemma}

\pf
Choose 2-dimensional subspaces $U_0, W_0$ of $U, W$ respectively, set 
$V_0:= U_0\otimes W_0$,
and consider the subgroup $L$ of ${\Ga L(V_0)}_{U_0\otimes 
W_0}$ induced by $H_{U_0\otimes W_0}$.
Let $\Delta$ denote the set of weight 2 vectors of $U_0\otimes W_0$ 
(considered as vectors of $V$),
and let $v\in\Delta$. By Lemma~\ref{weight}, $H_v\leq H_{U_0\otimes 
W_0}$, and since $H$ is $p$-exceptional
$|H:H_v|$ is coprime to $p$. It follows from Lemma~\ref{firstfacts} that 
$p$ divides $|L|$ and
a Sylow $p$-subgroup $P$ of $L$ acts diagonally on $U_0\otimes W_0$.
We may assume that $P$ acts as a group of field automorphisms of order $p$.
Then the set of fixed points of $P$ in $U_0\otimes W_0$ forms an 
$\F_q$-space
$U'_0\otimes W'_0=\F_q^2\otimes \F_q^2$.
In particular we may choose a basis $u_1,u_2$ for $U_0$ from $U'_0$.

By Remark~\ref{rem:weight}, $|\Delta|=|GL_2(r)|$. Also each 
$v'\in\Delta$ has a unique expression
as $v'=u_1\otimes w_1+ u_2\otimes w_2$ where $w_1, w_2$ span $W_0$, and 
it is straightforward to
prove that $P$ fixes $v'$ if and only if $P$ fixes $w_1$ and $w_2$, that 
is to say, if and only if $v'\in U'_0
\otimes W'_0$. Thus $P$ fixes exactly $|GL_2(q)|$ vectors in $\Delta$. 
Since each $v'\in\Delta$
is fixed by at least one Sylow $p$-subgroup of $L$ (by 
$p$-exceptionality), it follows that
the number $|L:N_L(P)|$ of Sylow $p$-subgroups of $L$
is at least
\[
y:= \frac{|GL_2(r)|}{|GL_2(q)|}=\frac{(r^2-1)(r^2-r)}{(q^2-1)(q^2-q)}.
\]
Let $Z_0=Z(GL(V_0))$, and note that $Z_0$ is contained in $L_0:= L\cap 
GL(V_0)$
since $H$ contains $Z$. Now $N_L(P)\cap Z_0\cong Z_{q-1}$, and $N_L(P)$ 
contains
$L\cap (GL(U_0')\circ GL(W_0'))$.

Recall the definitions of the maps $\phi_{U_0},\phi_{W_0}$ at the 
beginning of Section 4.2.
 From the classification of subgroups of $PGL_2(r)$ 
(see \cite[Chapter XII]{Dickson}), it 
follows that the
$p'$-group $L_0$ is such that each of $\phi_{U_0}(L_0), \phi_{W_0}(L_0)$ is
either a subgroup of $D_{2(r\pm1)}$ or equals one of $\AAA_4, \SSS_4, \AAA_5$. In 
the latter three cases,
$p$ would be odd, and such subgroups would lie in a subfield subgroup 
$PGL_2(q)$, and be centralised by
$P$. 

Suppose for instance that $\phi_{U_0}(L_0)$ is $\AAA_4$, $\SSS_4$, or $\AAA_5$. As $L_0$ is 
$p'$-group, we get $p \geq 5$, $r=q^p \geq 32$, so $2(r+1)$ is the largest 
possible order for $\phi_{W_0}(L_0)$. Also, $L_0 \geq Z_0$. Hence
for $K_0 = \Ker(\phi_{U_0}) \cap L_0$, we have that $|K_0/Z_0| \leq 2(r+1)$.
On the other hand, we noted that $L_0/K_0$ is centralized by $P$, and
$L=PL_0$. It follows that $L=NK_0 = NZ_0K_0$ for $N=N_L(P)$. Now
$$\begin{array}{ll}|L:N| &= |L:NZ_0||NZ_0:N| = |NZ_0K_0:NZ_0||NZ_0:N|\\
      & = |K_0/(NZ_0 \cap K_0)||Z_0/(N \cap Z_0)|\\
    &   \leq |K_0/Z_0||Z_0/(N \cap Z_0)|\\
     &  \leq 2(r+1)(r-1)/(q-1)
\end{array}$$
which is strictly less than $y$,
giving a contradiction.
Thus $L_0/Z_0\leq D_{2(r+\ve)}\times D_{2(r+\ve')}$, for some $\ve, 
\ve'\in\{1,-1\}$.

Suppose first that $p$ is odd. Then $N_{L_0}(P)/N_{Z_0}(P)\leq 
D_{2(q+\ve)}\times D_{2(q+\ve')}$, and we find
\[
|L:N_L(P)|\leq \frac{r-1}{q-1}.\frac{r+\ve}{q+\ve} 
.\frac{r+\ve'}{q+\ve'} \leq \left(\frac{r-1}{q-1}\right)^3
\]
which is less than $y$, contradiction. Hence $p=2$.

Suppose now that $a=\min\{a,b\}\geq 3$. We may repeat the above analysis 
with 3-dimensional subspaces
$U_0,W_0$ and $\Delta$ the set of weight 3 vectors in $V_0=U_0\otimes 
W_0$: the
cardinality  $|\Delta|$ is $|GL_3(r)|$, $P$ fixes $|GL_3(q)|$ vectors in 
$\Delta$, and  the number
of Sylow $2$-subgroups of $L$ is at least
\[
y': = \frac{|GL_3(r)|}{|GL_3(q)|}=\frac{(r^3-1)(r^3-r)(r^3-r^2)}{(q^3-1)(q^3-q)(q^3-q^2)}. 
\]
Again $L_0:=L\cap GL(V_0)$ is such that each of $\phi_{U_0}(L_0)$ and 
$\phi_{W_0}(L_0)$ is an odd
order subgroup of $GL_3(r)$, and hence is completely reducible. Thus 
each of these subgroups is a
subgroup of one of $Z_{r^2+r+1}.3$, $(Z_{r^2-1}\times Z_{r-1})/Z_0$ or 
$Z_{r-1}^3/Z_0$, and it follows that
\[
|L:N_L(P)|\leq   
\frac{r-1}{q-1}.\left(\max\left\{\frac{r^2+r+1}{q^2+q+1}, (r+1) 
\frac{r-1}{q-1}, \left(\frac{r-1}{q-1}\right)^2\right\}\right)^2
\]
which equals 
$(q^2+1)^2(q+1)^3$ (recall that $r=q^2$ here). However this quantity is less than 
$y'$ and we have a contradiction. Thus $a=2$.
\hal

\section{Tensor products II: $\mathcal{C}_7$ case}\label{c7sec}

In this section we classify $p$-exceptional groups which preserve tensor-induced decompositions. By this we mean the following. Let $V_1$ be a vector space over $\F_q$, and let $V = V_1^{\otimes t} = V_1\otimes V_2 \otimes \cdots \otimes V_t$, a tensor product of $t$ spaces isomorphic to $V_1$. The group $(GL(V_1) \circ \cdots \circ GL(V_t)).\SSS_t$ acts on $V$, where all centres are identified in the central product and the group $\SSS_t$ permutes the tensor factors. If $G$ is a subgroup of this group we say that $G$ preserves the tensor-induced decomposition $V = V_1^{\otimes t}$.

\begin{thm}\label{c7m}
Assume $G < GL(V)$ is a (not necessarily 
irreducible) $p$-exceptional group which preserves a tensor-induced decomposition
$$\F_q^n = V = (V_1)^{\otimes t} = V_1 \otimes V_2 \otimes \cdots \otimes V_t,$$
where $\dim_{\F_q}V_i = m \geq 2$ and $t \geq 2$. Then $p = 2$, and one of the 
following holds:

{\rm (i)} $t = 4$ and $m = q = 2$;

{\rm (ii)} $t = 3$, and $m = 2,3$. Moreover, if $m = 3$ then $q = 2$, and 
if $m = 2$ then $q \leq 4$;

{\rm (iii)} $t = 2$, and $m = 2,3$. Moreover, if $m = 3$ then $q \leq 8$.
\end{thm}

We shall also need the following result for the semilinear case.

\begin{thm}\label{c10m}
Let $G \leq \Gamma L(V)$, and assume that $G_0 \trianglelefteq  G = \langle G_0,\sigma \rangle$, where
\begin{enumerate}
\item[\rm (i)] $G_0$ is an absolutely irreducible $p'$-subgroup of $GL(V)$ which 
preserves a tensor-induced decomposition
$$V = (V_1)^{\otimes t} = V_1 \otimes V_2 \otimes \cdots \otimes V_t,$$
with $V = \F_q^n$, $\dim_{\F_q}V_i = m \geq 2$, $t \geq 2$, and
\item[\rm (ii)] $q = p^{pf}$, and $\sigma$ induces the field automorphism 
$x \mapsto x^{p^f}$ of $V$ modulo $GL(V)$.
\end{enumerate}
\no Then $G$ is not $p$-exceptional.
\end{thm}

The following result classifies the $p$-exceptional examples occurring in the cases left over by Theorem \ref{c7m}.

\begin{prop}\label{tensoregs}
Let $G < GL(V)$ be $p$-exceptional as in the hypothesis of Theorem $\ref{c7m}$, and suppose $t,m,q$ are as in one of conclusions {\rm (i)--(iii)} of the theorem. Suppose also that $G$ is irreducible on $V$. Then one of the following holds:

{\rm (a)} $m=3,t=2,q=2$: there are two irreducible $2$-exceptional groups $G$, of the form
$7^2.\SSS_3$ (orbit lengths $1, 21, 49^7, 147$) and $(7.3)^2.2$ (orbit lengths $1, 21, 49, 147^3$); both are imprimitive.

{\rm (b)} $m=2,t=2$: any irreducible $2$-exceptional group $G$ in this case is conjugate to a subgroup of $GL_2(q^2)$, hence is given by Lemma $\ref{lem:2dim}$.
\end{prop}

The proofs of these results are presented in the following three subsections.

\subsection{Proof of Theorem \ref{c7m}}

Throughout this section we assume that $G \leq GL(V)$ is a (not necessarily 
irreducible) $p$-exceptional group which preserves a tensor-induced decomposition
$$V = (V_1)^{\otimes t} = V_1 \otimes V_2 \otimes \cdots \otimes V_t,$$
where $V = \F_q^n$, $\dim_{\F_q}V_i = m \geq 2$, $t \geq 2$, and 
$(m,t,p) \neq (2,2,2)$. Let $B := G \cap (GL(V_1) \circ \ldots \circ GL(V_t))$ 
be the base group and let $H = G/B \leq \SSS_t$ be the permutation group 
induced by the action of $G$ on the $t$ tensor factors $V_i$.

\subsubsection{First reduction} 

We begin with some elementary observations. Recall that a {\it rational} element of a finite group is an element which is conjugate to all of its powers which have the same order.

\begin{lemma}\label{trivial}    
Under the above hypothesis, the following statements hold.

{\rm (i)} $B$ is a $p'$-group.

{\rm (ii)} $H$ is a transitive subgroup of $\SSS_t$ of order divisible by $p$.
In particular $t \geq p$.

{\rm (iii)} Let $1 \neq h \in G$ be any $p$-element and let $Q \leq G$ be 
any $p$-subgroup containing $h$. Then 
$$|G:N_G(Q)| > |V/C_V(h)|.$$

{\rm (iv)} $G \setminus B$ contains an element $g$ of order $p$, and for such an element,
$$\frac{|G|}{p \cdot |C_B(g)|} > |V/C_V(g)|.$$  
If in addition the element $gB$ is rational in $H = G/B$, then 
$$\frac{|G|}{p(p-1) \cdot |C_B(g)|} > |V/C_V(g|.$$ 
\end{lemma}

\pf
If $p$ divides $|B|$ then $B$ is $p$-exceptional by Lemma \ref{normalsubgroups}, which contradicts Theorem \ref{c4lin} (since we are assuming that $(m,t,p)\ne (2,2,2)$). Part (i) follows.
Likewise, if $H$ is intransitive, then $G$ preserves 
a nontrivial tensor decomposition of $V$ and we get a contradiction by the 
same result; hence (ii) holds. Next, the $p$-exceptionality of $G$ implies that any nonzero 
element $v \in V$ is fixed by a Sylow $p$-subgroup of $G$ and so by  
a conjugate of $Q$ as well. Hence, 
$$|V|-1 = |V^\sharp| \leq |G:N_G(Q)| \cdot |C_V(Q)| \leq 
  |G:N_G(Q)| \cdot |C_V(h)|,$$
and (iii) follows. Since $p$ divides $|G|$ and $B$ is a $p'$-group, we can find
$g \in G \setminus B$ of order $p$. Now we choose $h := g$ and 
$Q := \langle g \rangle$ in (iii). Observe that $C_G(g)$ contains $g$ and the 
$p'$-subgroup $C_B(g)$, whence the first inequality in (iv) follows.
Finally, since $B$ is a $p'$-group, the rationality of $gB$ in $G/B$ implies
$g$ is rational in $G$ (see e.g. \cite[Lemma 4.11]{TZ2}), in which case 
we have $|N_G(Q)| = (p-1)|C_G(g)|$. Hence the second 
inequality in (iv) follows.  
\hal

\vspace{2mm}
We fix the element $g$ in Lemma \ref{trivial}(iv) and bound 
$\kappa := (\dim C_V(g))/(\dim V)$. Observe that $|V/C_V(g)| = |V|^{1-\kappa}$.
Replacing $G$ by some conjugate subgroup, we may assume that $g$
permutes $V_1, \ldots ,V_p$ cyclically:
$$g~:~V_1 \mapsto V_2 \mapsto V_3 \mapsto \ldots \mapsto V_p \mapsto V_1.$$
Choose a basis $(e^1_j \mid 1 \leq j \leq m)$ of $V_1$ and let 
$e^i_j := (e^1_j) g^{i-1}$ for $1 \leq i \leq p$. Since $|g| = p$, we see
that $(e^p_j)g = e^1_j$ and $(e^i_j \mid 1 \leq j \leq m)$ is a basis of $V_i$
for $1 \leq i \leq p$. Clearly, 
$$(e^1_{j_1} \otimes e^2_{j_2} \otimes \cdots \otimes e^p_{j_p} \mid 
  1 \leq j_1, \ldots ,j_p \leq m)$$ 
is a basis of $U := V_1 \otimes \cdots \otimes V_p$. Let $J_k$ denote the 
Jordan block of size $k$ with eigenvalue $1$. Then $g\downarrow U$ permutes the above basis
vectors of $U$ in $(m^p-m)/p$ cycles of length $p$, hence has
Jordan canonical form $(J_1^a ,J_p^b)$, where
$a := m$, $b:= (m^p-m)/p$. Also let $W := V_{p+1} \otimes \cdots \otimes V_t$ 
so that $V = U \otimes W$.

\begin{lemma}\label{for-g} We have 
$$\kappa = \frac{\dim C_V(g)}{\dim V} \leq 
  \frac{1}{p} + \frac{1-\frac{1}{p}}{m^{p-1}}.$$  
\end{lemma}   

\pf
Consider any indecomposable direct summand $W'$ of the 
$\langle g \rangle$-module $W$. Suppose $g$ acts on $W'$ via Jordan block $J_k$.
Then $g$ acts on $U \otimes W'$ with Jordan canonical form 
$(J_1^a,J_p^b) \otimes J_k = (J_k^a, J_p^{bk})$. It follows using the values of $a,b$ above, that 
$$\frac{\dim C_{U \otimes W'}(g)}{\dim (U \otimes W')} = 
  \frac{a+bk}{k(a+bp)} \leq \frac{a+b}{a+bp} = \frac{\dim C_U(g)}{\dim U} = 
  \frac{1}{p} + \frac{1-\frac{1}{p}}{m^{p-1}}.$$
Applying this observation to every indecomposable direct summand $W'$ of the 
$\langle g \rangle$-module $W$, we get 
$\kappa \leq (a+b)/(a+bp)$, yielding the desired inequality.   
\hal
 
Next we estimate $|B:C_B(g)|$.

\begin{lemma}\label{b-size}
Let $X$ be a $p'$-subgroup of $PGL_m(q)$ of largest possible order.

{\rm (i)} Let $h \in G$ be an arbitrary element. Then 
$$|B:C_B(h)| \leq |X|^t \leq (|PGL_m(q)|_{p'})^t.$$

{\rm (ii)} If, in addition, $g$ acts trivially on $W$, then 
$$|B:C_B(g)| \leq |X|^{p-1} \leq (|PGL_m(q)|_{p'})^{p-1}.$$ 
\end{lemma}

\pf
Recall that for a tensor product space $\F_q^k \otimes \F_q^l$, if $A, C \in GL_k(q)$ and $D,E \in GL_l(q)$ are such that 
$A \otimes D = C \otimes E$, then $C = \a A$ and $E = \a^{-1}D$ for some
$\a \in \F_q^*$. It follows that the map
$$f~:~B \to PGL(V_1) \times \cdots \times PGL(V_t),$$
defined by $f(x) = (\bar A_1, \ldots, \bar A_t)$ if 
$x = A_1 \otimes A_2 \otimes \cdots \otimes A_t \in B$, $A_i \in GL(V_i)$, and
$\bar A_i$ denotes the coset containing $A_i$ in $PGL(V_i)$, is a well-defined
homomorphism. Observe that each fibre of $f$ is contained in 
exactly one $C_B(h)$-coset in $B$. 
Indeed, if $f(x) = f(x')$, then $x' = \beta x$ for some $\beta \in \F_q^*$,
and so $x'x^{-1} \in C_B(h)$. Furthermore, since each $x \in B$ is a $p'$-element
by Lemma \ref{trivial}(i),
the elements $\bar A_i$ are $p'$-elements in $PGL(V_i)$. Composing $f$ with
the projection $PGL(V_1) \times \cdots \times PGL(V_t) \to PGL(V_i)$, we 
therefore get a homomorphism $f_i~:~B \to PGL(V_i)$ with $f_i(B)$ being a 
$p'$-group. It follows that $|f_i(B)| \leq |X| \leq |PGL_m(q)|_{p'}$. Now   
$|f(B)| \leq \prod^{t}_{i=1}|f_i(B)|$, whence (i) follows.

For (ii), notice that
\begin{equation}\label{yidef}
Y_i := \{ C \in GL(V_i) \mid 
  C = A_i \mbox{ for some }h = A_1 \otimes \cdots \otimes A_t \in B \}
\end{equation}
is a $p'$-subgroup of $GL(V_i)$.
Given the action of $g$ on $V_1, \ldots, V_p$, we can identify 
$Y_i$, $1 \leq i \leq p$, with $Y_1$.
Consider the homomorphism 
$$\fs~:~Y := Y_1 \times \cdots \times Y_t
  \to Y_1 \circ \cdots \circ Y_t$$
given by $\fs(A_1, \ldots ,A_t) = A_1 \otimes \cdots \otimes A_t$, 
and note that $B \le f^*(Y)$. Let 
$$\begin{array}{rl}
  K & := \{(\underbrace{A, \ldots ,A}_{p},D_{p+1}, \ldots ,D_t) \mid 
  A \in Y_1,D_i \in Y_i \},\\
  Z_0 & := \left\{(a_1I_m, \ldots ,a_tI_m) \mid a_i \in \F_q^* \right\}.
  \end{array}$$ 
Since the element $g$ of order $p$ acts trivially on $W$ and permutes 
$V_1, \ldots, V_p$ cyclically, $f^*(KZ_0)$ centralizes $g$. Thus 
$f^*(KZ_0) \le C_{f^*(Y)}(g)$ and
\[
|B:C_B(g)| = |g^B| \le |g^{f^*(Y)}| \le \frac{|f^*(Y)|}{|f^*(KZ_0)|} = \frac{|Y|}{|KZ_0|} =
  \left(\frac{|Y_1|}{q-1}\right)^{p-1}.
\]
It remains to observe that $|Y_1| \leq (q-1)|X|$.         
\hal

\begin{lemma}\label{c71}
Under the above assumptions, one of the following holds.

{\rm (i)} $p = 3$ and $(m,t) = (2,4)$, $(2,3)$. 

{\rm (ii)} $p = 2$. Furthermore, either $t = 2$ or 
$(m,t) = (4,3)$, $(3,3)$, $(2,6)$, $(2,5)$, 
$(2,4)$, $(2,3)$. 
\end{lemma}

\pf
 By Lemma \ref{b-size}, for $g$ the element defined before Lemma \ref{for-g} we have
$$\begin{array}{ll}|G:C_B(g)| & \leq |H| \cdot |B:C_B(g)] 
   \leq (t!) \cdot (|PGL_m(q)|_{p'})^t \\ \\
   & < \left( \frac{t+1}{2} \cdot q^{m(m+1)/2-1} \right)^t 
     < |V|^{f(m,t,q)},\end{array}$$
where 
$$f(m,t,q) = t \cdot \frac{m(m+1)/2 + \log_q\frac{t+1}{2}-1}{m^t}.$$ 
In particular, if $t \geq 5$, then $f(m,t,q) < 0.6$. 
By Lemma \ref{trivial}(iv), 
\[
|G:C_B(g)| > p|V/C_V(g)| > q^{\dim V - \dim C_V(g)} = |V|^{1-\kappa },
\]
and so $f(m,t,q)+\kappa > 1$.

If $p \geq 5$, then $t\ge 5$ by Lemma \ref{trivial}(ii) and so 
$f(m,t,q) < 0.6$. Then by Lemma \ref{for-g}, $\kappa \leq 1/4$, a contradiction. 

Now assume that $p = 3$, so $t\ge 3$ by Lemma \ref{trivial}(ii). If $t \geq 5$, then $f(m,t,q) < 0.47$.
If $t = 4$ and $m \geq 3$, then $f(m,t,q) < 0.3$.  
If $t = 3$ and $m \geq 4$, then $f(m,t,q) < 0.46$. In all these cases
$\kappa \leq 1/2$, and we arrive at a contradiction as above. 

Consider the case $m=t=3$ (still with $p=3$). If $q \geq 9$, then 
$f(m,t,q)+\kappa < 0.5907 + 11/27 < 1$,
again a contradiction. Assume $q = 3$. Then $|C_V(g)| \leq 3^{11}$ by Lemma 
\ref{for-g}. On the other hand, by Lemma \ref{b-size} we have 
$$|G:C_B(g)| \leq 2 \cdot (|PGL_3(3)|_{3'})^3 < 3^{16},$$
contradicting Lemma \ref{trivial}(iv). The remaining cases are listed in (i).

Now we consider the case $p=2$. If $t \geq 7$, then 
$f(m,t,q) + \kappa < 0.22 + 0.75 < 1.$
If $m \geq 3$ and $t \geq 4$, then 
$f(m,t,q) + \kappa < 0.32 + 2/3 < 1.$
And if $m \geq 5$ and $t \geq 3$, then 
$f(m,t,q) + \kappa < 0.36 + 0.6 < 1.$ The remaining cases are listed in (ii). \hal

\subsubsection{The case $p = 3$}

\begin{prop}\label{c73}
The case $p = 3$ is impossible.
\end{prop}

\pf
Observe that the subgroup $X$ in Lemma \ref{b-size} has order
$2(q+1)$ by Lemma \ref{small}(i), and $\kappa \leq \frac{1}{2}$ by Lemma \ref{for-g}.
So $|V/C_V(g)| \ge |V|^{1-\kappa} \ge q^{m^t/2}$.
By Lemma \ref{c71} we need to consider two cases.

(1) Suppose $(m,t) = (2,4)$. In this case, Lemma \ref{trivial}(iv)
and Lemma \ref{b-size} imply
$$q^8 \leq |V/C_V(g)| < 8 \cdot (2(q+1))^4,$$
whence $q = 3$. Assume in addition that $g$ acts nontrivially on 
$V_4 = \F_3^2$, i.e. $V_4\downarrow g = J_2$. Then 
$$V\downarrow g = (J_1^2,J_3^2) \otimes J_2 = (J_2^2,J_3^4)$$
and so Lemma \ref{trivial}(iv) implies
$$q^{10} \leq |V/C_V(g)| < 8 \cdot (2(q+1))^4,$$
again a contradiction. 

Thus $g$ acts trivially on $V_4 = W$. In this case, by Lemma \ref{b-size}
we have $q^8 < |G:C_G(g)| \leq 8 \cdot (2(q+1))^2$, again a contradiction. 

(2) Assume now that $(m,t) = (2,3)$. In particular, $H = \AAA_3$ or $\SSS_3$,
and in the latter case $g$ is rational. Also, $g$ acts trivially on $W$. Hence
by Lemma \ref{trivial}(iv) we have 
$$q^4 \leq |V/C_V(g)| < |B:C_B(g)| \leq (2(q+1))^2,$$
a contradiction as $q \geq 3$.  
\hal

\subsubsection{The case $p = 2$}\label{tp2}

\begin{lemma}\label{456}
Assume $p = 2$ and $t \geq 4$. Then $t = 4$ and $m = q = 2$, i.e. 
$V = \F_2^{16}$.
\end{lemma}

\pf
 Recall that $H$ is an even-order transitive subgroup of $\SSS_t$. We claim
that $H$ contains an involution $h = dk$ with $d \in \SSS_4$ a double 
transposition and $k \in \SSS_{t-4}$ disjoint from $d$. (If not, then we may 
assume $H \ni x = (1,t)$. Since $H$ is transitive, for any $2 \leq i \leq t-1$ 
we can find $u \in H$ with $1^u = i$, and so $H \ni x^u = (i,t^u)$. If 
$t^u \neq 1,t$ then $H \ni x \cdot x^u = (1,t)(i,t^u)$, contrary to our
assumption. If $t^u = t$, then $x^u = (i,t) \in H$. If $t^u = 1$, then 
$H \ni x^{ux} = (i,t)$. We have shown that $(i,t) \in H$ for all $i$ with  
$1 \leq i \leq t-1$, and so $H = \SSS_t \ni (12)(34)$, again a contradiction.)

Without loss of generality we may now assume that $G$ contains an involution $h$ 
which permutes $V_1$ with $V_2$, $V_3$ with $V_4$, and acts on 
$\{V_5, \ldots ,V_t\}$. Arguing as in the 
discussion about $g$ preceding Lemma \ref{for-g}, we see that 
$$(V_1 \otimes V_2)\downarrow h = (V_3 \otimes V_4)\downarrow h = 
   (J_1^a,J_2^b).$$
Then setting  $M := V_1 \otimes \cdots \otimes V_4$ and arguing as in the proof of Lemma \ref{for-g} we obtain 
$$\gamma := \frac{\dim C_V(h)}{\dim V} \leq 
  \frac{\dim C_M(h)}{\dim M} = \frac{a^2+2ab+2b^2}{m^4} = 
  \frac{1}{2} + \frac{1}{2m^2}.$$
 In particular, $\gamma \leq 5/8$.
Now we will apply Lemmas \ref{trivial} and \ref{b-size}(i) to $h$ instead of 
$g$, and treat the cases described in Lemma \ref{c71} separately.

Assume first that $(m,t) = (2,6)$. Then $|C_V(h)| \leq |V|^{\gamma} \leq q^{40}$, 
so $|V/C_V(h)| \ge q^{24}$. On the 
other hand, by Lemmas \ref{b-size} and \ref{small}(i), we have
$$\frac{1}{2}|G:C_B(h)| \le |G:C_G(h)| \leq 360 \cdot (q^2-1)^6 < q^{24},$$
contradicting Lemma \ref{trivial}(iv).

Next assume that $(m,t) = (2,5)$. Then $|C_V(h)| \leq |V|^{\gamma} \leq q^{20}$, so $|V/C_V(h)| \ge q^{12}$. 
Hence by Lemmas \ref{trivial}(iv) and \ref{b-size} we must have
$$q^{12} < \frac{1}{2}|G:C_B(h)| \le |G:C_G(h)| \leq 60 \cdot (q^2-1)^5,$$
and so $q = 2$ or $4$. If $q = 4$, then by Lemmas \ref{small} and \ref{b-size} we have
$|B:C_B(h)| \leq (q+1)^5$, whence 
$$|G:C_G(h)| \leq 60 \cdot (q+1)^5  < q^{12},$$
a contradiction. If $q = 2$, then for $Q \in Syl_2(G)$ we have  
$$|G:N_G(Q)| \leq 15 \cdot 3^5 < 2^{12},$$
again a contradiction by Lemma \ref{trivial}(iii).       

Finally, we assume that $(m,t) = (2,4)$ and $q \geq 4$. Then 
$|C_V(h)| \leq |V|^{\gamma} \leq q^{10}$. Also, any $2'$-subgroup of 
$PGL_2(q)$ has order at most $q+1$ by Lemma \ref{small}(i), and $|H| \leq |\SSS_4|$.  
In particular, the involution $h$ is central in some 
$Q \in Syl_2(G)$. Furthermore, $C_B(h)$ has odd order, so
$|C_G(h)| \geq |Q| \cdot |C_B(h)|$. It follows that 
$$|G:C_G(h)| \leq |H:Q| \cdot |B:C_B(h)| \leq 3 \cdot (q+1)^4  < q^{6},$$
a contradiction for $q \geq 4$.    
\hal

\begin{lemma}\label{34}
Suppose that $p = 2$ and $t = 3$. Then either $m = 3$ and $q = 2$, or 
$m = 2$ and $q \leq 4$.
\end{lemma}

\pf
By Lemma \ref{c71} we need to distinguish two cases.

(1) Assume that $(m,t) = (4,3)$. Then $\kappa \leq 5/8$ by Lemma \ref{for-g}
and so $|C_V(g)| \leq q^{40}$. Observe that any $2'$-subgroup $X$ of
$PGL_4(q) \cong SL_4(q)$ has order $\leq (q^3-1)(q^2-1)(q-1)$. (Indeed, this follows from 
Lemma \ref{small} if $X$ acts reducibly on $V_1 = \F_q^4$. Suppose that this action is
irreducible. Then $\Hom_X(V_1) \cong \F_{q^a}$ for some $a|4$, and 
$V_1$ is a $(4/a)$-dimensional absolutely irreducible $\F_{q^a}X$-module. 
Since $X$ is soluble, any irreducible Brauer character of $X$ lifts to 
a complex character by the Fong-Swan Theorem \cite[72.1]{dorn}; in particular,
$4/a$ divides $|X|$ and so $a = 4$. This in turn implies that 
$X \leq GL_1(q^4)$, and so we are done.) 
Hence, applying Lemma \ref{b-size} to the element $g$ defined before Lemma \ref{for-g}, we have
$$q^{24} < |G:C_G(g)| \leq 3 \cdot ((q^3-1)(q^2-1)(q-1))^3,$$
a contradiction.

(2) Consider the case $(m,t) = (3,3)$ and $q \geq 4$. Note that any
$2'$-subgroup of $PGL_3(q)$ has order $\leq q^3-1$ by Lemma \ref{small}(ii).
Assume in addition that the involution
$g$ acts nontrivially on $V_3 = \F_q^3$, i.e. $V_3\downarrow g = (J_1,J_2)$. Then 
$$V \downarrow g = (J_1^3,J_2^3) \otimes (J_1,J_2) = (J_1^3,J_2^{12})$$
and so $|C_V(g)| = q^{15}$. Lemmas \ref{trivial}(iv) and \ref{b-size} now imply that
$$q^{12}  < |G:C_G(g)| \leq 3 \cdot (q^3-1)^3,$$
a contradiction. 

Thus $g$ acts trivially on $V_3 = W$. In this case, by Lemma \ref{b-size}(ii)
we have 
$$q^9 \leq |V/C_V(g)| < |G:C_G(g)| \leq 3 \cdot (q^3-1),$$ 
again a contradiction. 

(3) Now assume $(m,t) = (2,3)$ and $q \geq 8$. Note that any
$2'$-subgroup of $PGL_2(q)$ has order $\leq q+1$ by Lemma \ref{small}(i).
Assume in addition that the involution
$g$ acts nontrivially on $V_3 = \F_q^2$, i.e. $V_3\downarrow g = J_2$. Then 
$$V \downarrow g = (J_1^2,J_2) \otimes J_2 = J_2^4$$
and so $|C_V(g)| = q^4$. Lemmas \ref{trivial}(iv) and \ref{b-size} now imply that
$$q^{4}  < |G:C_G(g)| \leq 3 \cdot (q+1)^3,$$
a contradiction. 
Thus $g$ acts trivially on $V_3 = W$. In this case, by Lemma \ref{b-size}(ii)
we have 
$$q^2 \leq |V/C_V(g)| < |G:C_G(g)| \leq 3 \cdot (q+1),$$ 
again a contradiction. 
\hal

\vspace{4mm}
The rest of this subsection is devoted to the case $t = 2$, so 
$V = V_1 \otimes V_2$ and $p = 2$ by Lemma \ref{trivial}(ii),
and $V_1$, $V_2$ are interchanged by $g$. 
As before, we fix the basis $(e_j := e^1_j \mid 1 \leq j \leq m)$ of 
$V_1$ and $(f_j := e^2_j = e_jg \mid 1 \leq j \leq m)$ of $V_2$. 
Then we can identify both $V_1$ and $V_2$ with $\F_q^m$.
Consider the subgroups $Y_1 \cong Y_2$ of $GL_m(q)$ 
defined in (\ref{yidef}). 

The key observation in the case $t=2$ is the following:

\begin{lemma}\label{t21}
Suppose $t = 2$. Then the subgroup $Y_1 < GL_m(q)$ is transitive on 
$k$-dimensional subspaces of $\F_q^m$ for any $k \leq m-1$.
\end{lemma}

\pf
 Recall that if $0 \neq v \in V$ has weight $k$: 
$v = \sum^{k}_{i=1}x_i \otimes y_i$, then 
$$[v]_1 := \langle x_1, \ldots ,x_k \rangle_{\F_q},~~
  [v]_2 := \langle y_1, \ldots ,y_k \rangle_{\F_q}$$
are $k$-dimensional subspaces of $\F_q^m$ uniquely determined by $v$. 
In particular, if $vg = v$, then 
$$\sum^{k}_{i=1}x_i \otimes y_i = v = vg = \sum^{k}_{i=1}y_i \otimes x_i,$$
and so $([v]_1)g = [v]_2$. 

 Next we show that if $0 \neq v \in V$ has weight $k$ and is fixed by some 
involution $g'$ of $G$, then there is some element $a \otimes b \in B$ such that 
$[v]_2b = ([v]_1a)g$. Indeed, since $|G|=2|B|$ and $|B|$ is odd, we must have 
$g' = xgx^{-1}$ for some $x = a \otimes b \in B$. Writing 
$v = \sum^{k}_{i=1}x_i \otimes y_i$, we see that 
$w := vx = \sum^{k}_{i=1}x_ia \otimes y_ib$ is fixed by $g$. Hence, according to
the first paragraph we then have 
$$([v]_1a)g = ([w]_1)g = [w]_2 = ([v]_2)b,$$
as stated.

 Now we set $L := \langle e_1, \ldots ,e_k \rangle_{\F_q}$, 
$M := \langle f_1, \ldots ,f_k \rangle_{\F_q}$, and consider any $k$-dimensional 
subspace $N := \langle u_1, \ldots ,u_k \rangle_{\F_q}$ of $V_1$. Then 
$v = \sum^{k}_{i=1}u_i \otimes f_i \in V$ has weight $k$, with
$[v]_1 = N$, $[v]_2 = M$. The $p$-exceptionality of $G$ implies that $v$ is fixed
by some involution $g' \in G$. By the previous paragraph, there is some 
$x = a \otimes b \in B$ such that 
$$(Na)g = ([v]_1a)g = ([v]_2)b = Mb = (Lg)b,$$
i.e. $Na = Lgbg^{-1}$. Writing $f_jb = \sum^m_{i=1}b_{ij}f_i$ for some 
$b_{ij} \in \F_q$, we have 
$(e_j)gbg^{-1} = \sum^m_{i=1}b_{ij}e_i$. Thus, under our identification of $V_1$ and 
$V_2$ with $\F_q^m$, we have $Na = Lgbg^{-1} = Lb$, and so $Lba^{-1} = N$. 
It remains to observe that $B$ contains $x^gx^{-1} = (ba^{-1}) \otimes (ab^{-1})$, 
whence $ba^{-1} \in Y_1$.     
\hal

\begin{prop}\label{t22}
Assume $t = 2$ ( so $p=2$). Then either $m = 2$, or $m = 3$ and 
$q \leq 8$.
\end{prop}

\pf 
By Lemma \ref{t21}, $Y_1$ is transitive on $k$-spaces for all $k$, and has odd order. Hence \cite[Lemma 4.1]{regsgps}
shows that if $m>2$ then either $m=3$ or $(m,q) = (5,2)$. In the latter case
$Y_1 = \Gamma L_1(32)$, $V = \F_2^{25}$, and $|C_V(g)| \le 2^{15}$ by Lemma \ref{for-g}. Applying Lemma 
\ref{trivial}(iv), $|B:C_B(g)| > 2^{10}$. On the other hand,
the proof of Lemma \ref{b-size}(ii) shows that 
$|B:C_B(g)| \leq |Y_1| = 155$, a contradiction.    

Next suppose that $m = 3$. By Lemma \ref{for-g} and 
Lemma \ref{trivial}(iv) we have 
$|B:C_B(g)| > |V/C_V(g)| \ge q^3$. Now the proof of Lemma \ref{b-size}(ii) 
shows that $q^3 < |B:C_B(g)| \leq |Y_1|/(q-1) \leq 3f(q^2+q+1)$, which can happen 
only when $q \leq 8$.  
\hal

\vspace{4mm}
This completes the proof of Theorem \ref{c7m}.

\subsection{Proof of Theorem \ref{c10m}}

Throughout this section we assume that $G \leq \Gamma L(V)$ is a  
$p$-exceptional group such that $G_0 \lhd G = \langle G_0,\sigma \rangle$, where

(i) $G_0$ is an absolutely irreducible $p'$-subgroup of $GL(V)$ which 
preserves a tensor-induced decomposition
$$V = (V_1)^{\otimes t} = V_1 \otimes V_2 \otimes \cdots \otimes V_t,$$
with $V = \F_q^n$, $\dim_{\F_q}V_i = m \geq 2$, $t \geq 2$, and

(ii) $q = q_0^p$, $q_0 = p^f$, and $\sigma$ induces the field automorphism 
$x \mapsto x^{q_0}$ of $V$ modulo $GL(V)$. 

\no Note that 
$|C_V(\sigma)| \leq |V|^{1/p}$. Indeed, since 
${\mathcal {G}} := GL(V \otimes_{\F_q}\overline{\F}_q)$ is connected
(where $\overline{\F}_q$ is the algebraic closure of $\F_q$), 
$\sigma$ is ${\mathcal {G}}$-conjugate to the standard Frobenius
morphism $\sigma_0~:~x \mapsto x^{q_0}$. Hence 
$$|C_V(\sigma)| \leq |C_{(V \otimes_{\F_q}\overline{\F}_q)}(\sigma)| = 
|C_{(V \otimes_{\F_q}\overline{\F}_q)}(\sigma_0)| = |V|^{1/p}.$$     

\subsubsection{First reductions}

The following lemma simplifies further computations.

\begin{lemma}\label{red1}
Under the above assumptions, 

{\rm (i)} $m^t$ divides $|G_0|$; in particular, $p$ does not divide $m$, and

{\rm (ii)} $|G_0| > |V|^{1-1/p}$. 
\end{lemma}

\pf
Part (i) follows from the assumptions that $G_0$ is a $p'$-group and absolutely irreducible.
For (ii), one can argue as in the proof of Lemma \ref{trivial}(iii), taking 
$Q = \langle \sigma \rangle$. 
\hal

\vspace{2mm}
Next we rule out most of the cases using Lemma \ref{red1}:

\begin{prop}
Under the above assumptions, $t = 2$.
\end{prop}

\pf
Assume to the contrary that $t \geq 3$.
The proof of Lemma \ref{b-size} implies that the base subgroup
$G_0 \cap (GL(V_1) \circ \ldots \circ GL(V_t))$ of $G_0$ has order at most
$(q-1) \cdot |X|^t$, where $X$ is a $p'$-subgroup of largest possible 
order of $PGL_m(q)$. Hence
$$|G_0| \leq t! \cdot (q-1) \cdot (|PGL_m(q)|_{p'})^t < 
  \left(\frac{t+1}{2}\right)^t \cdot q^{t(\frac{m(m+1)}{2}-1)+1},$$ 
and Lemma \ref{red1}(ii) yields that 
$$f(m,t,p) := m^{t}\left(1-\frac{1}{p}\right) - t \cdot \frac{m^2+m-2}{2} - 1 
  -t \cdot \log_p \frac{t+1}{2} < 0.$$
Note that the function $f(m,t,p)$ is non-decreasing for each of its variables. 
Now direct computations show that the latter inequality is impossible unless
one of the following holds.

(a) $t = 5$, $m = p = 2$. This is ruled out by Lemma \ref{red1}(i).

(b) $t = 4$, $m = 2$, $p \leq 3$. By Lemma \ref{red1}(i) we must have $p = 3$ and
so $q = q_o^p \geq 27$. Since for $m = 2$ we have $|X| \leq 2(q+1)$ 
by Lemma \ref{small}(i), 
$$|G_0| \leq 24 \cdot (2(q+1))^4 \cdot (q-1) < q^{32/3} = |V|^{2/3},$$
a contradiction by  Lemma \ref{red1}(ii).

(c) $(t,m,p) = (3,3,2)$. In this case $q \geq 4$, and $|X| \leq q^3-1$ as $m=3$
by Lemma \ref{small}(ii). It follows that   
$$|G_0| \leq 6 \cdot (q^3-1))^3 \cdot (q-1) < q^{27/2} = |V|^{1/2},$$
again contradicting Lemma \ref{red1}(ii).

(d) $(t,m) = (3,2)$. By Lemma \ref{red1}(i) we must have $p \geq 3$ and so
$q \geq 27$. Also, by Lemma \ref{small}(i) we have $|X| \leq 2(q+1)$ as $m = 2$. Hence 
$$|G_0| \leq 6 \cdot (2(q+1))^3 \cdot (q-1) < q^{16/3} \leq |V|^{1-1/p},$$
again a contradiction.
\hal

\subsubsection{The case $t = 2$}\label{t2cas}

Throughout this subsection we assume $t = 2$. 

\begin{lemma}\label{baseobs}
We have $G_0 \le GL(V_1)\circ GL(V_2)$ and $p=2$.
\end{lemma}

\pf Let $B = G_0 \cap (GL(V_1)\circ GL(V_2))$, and suppose $G_0 \ne B$. Then $G_0 = B\la s \ra = B.2$, where $s$ interchanges $V_1$ and $V_2$.
Since $G_0$ is a $p'$-group this implies that $p>2$. As $\s$ has order $p$ it therefore fixes $V_1$ and $V_2$, and so 
$B\la \s \ra$ is a normal subgroup of index 2 in $G$, hence is $p$-exceptional. This contradicts Theorem \ref{c4semi}.

Hence $G_0 = B$. If $p>2$ then again $\s$ fixes $V_1$, $V_2$ and we contradict Theorem \ref{c4semi}. Hence $p=2$, completing the proof. \hal

\begin{prop}\label{c104}
The case $t=p=2$ cannot occur.
\end{prop}

\pf
Assume to the contrary that $t = p = 2$. 
By Lemma \ref{baseobs} we have $G_0 = B$ and $G = B\la \s \ra$ with $\s$ of order 2, and
$G_0$ is an absolutely irreducible 
$2'$-group on $V$. By Lemma \ref{red1}, $m$ is odd; in particular, $m \geq 3$. 
If $\s$ fixes both $V_1$ and $V_2$, then Theorem \ref{c4semi} gives a contradiction. 
So $\s$ interchanges $V_1$ and 
$V_2$ and also it is semilinear: $(\lambda v)^\s = \lambda^rv$ with 
$q = 2^{2f} = r^2$. 

We will now follow the arguments in Subsection \ref{tp2} for the corresponding case 
in Theorem \ref{c4lin}, and indicate  necessary modifications because of the 
semilinearity of $\s$. We fix the basis $(e_j := e^1_j \mid 1 \leq j \leq m)$ of 
$V_1$ and $(f_j := e^2_j = \s(e_j) \mid 1 \leq j \leq m)$ of $V_2$. 
Then we can identify both $V_1$ and $V_2$ with $\F_q^m$.
Consider the subgroups $Y_1 \cong Y_2$ of $GL_m(q)$ 
defined in (\ref{yidef}). Note that if 
$x = X \otimes Y \in B$, then $x^{\sigma} = Y^{(r)} \otimes X^{(r)}$, 
where $X^{(r)} = (x_{ij}^r)$ if $X = (x_{ij})$. In particular, if
$X \in Y_1$ then $X^{(r)} \in Y_2$ and vice versa, whence $Y_1 \cong Y_2$. 
Now the proof of Lemma \ref{t21} can be carried over verbatim, except that
we have to replace $ba^{-1}$ by $b^{(r)}a^{-1}$. 
Thus $Y_1$ (and $Y_2$) is transitive
on $k$-spaces of $\FQ^m$ for all $k$, and $|Y_1|$ is odd (and contains all scalar
matrices). 
Now \cite[Lemma 4.1]{regsgps} implies that $m=3$ or $(m,q) = (5,2)$.
The latter is impossible as $q=r^2$.

Hence $m = 3$. Now we consider the homomorphism
$\fs~:~Y_1 \times Y_2 \to Y_1 \otimes Y_2$ defined by 
$\fs(X,Y) = X \otimes Y$, and the subgroup $K := \{(X,X^{(r)}) \mid X \in Y_1\}$ 
of $Y_1 \times Y_2$. Note that $\s$ centralizes $\fs(K)$. 
As in the proof of Lemma \ref{b-size}(ii), this implies that $|B:C_B(\s)|$ 
has order at most $|(Y_1 \times Y_2):K| = |Y_1|$. On the other hand,
$|C_V(\s)| = |V|^{1/2} = q^{9/2}$. Hence the $2$-exceptionality of $G$ implies
$|G:C_G(\s)| \geq |V/C_V(\s)| = q^{9/2}$ by Lemma \ref{bound}(ii). Since $|C_G(\s)| = 2|C_B(\s)|$, we get
$$q^{9/2} \leq |Y_1| \leq 6f(q^3-1),$$
again a contradiction.  
\hal

\vspace{6mm}
This completes the proof of Theorem \ref{c10m}.

\subsection{Proof of Proposition \ref{tensoregs}}

Now we prove Proposition \ref{tensoregs} by treating the cases left in conclusions (i)-(iii) of Theorem \ref{c7m}.
Let $G$ be as in the hypothesis of the proposition. If $(m,t) \ne (2,2)$ then there are just a small number of cases for $(m,t,q)$ to consider, all with $q \le 8$ and $m^t \le 16$, and a \textsc{Magma} computation shows that the only irreducible $p$-exceptional examples are those in part (a) of Proposition \ref{tensoregs}. (For $(m,t)=(3,2)$ and $q=4, 8$ we first use Lemmas \ref{trivial}(i) and \ref{t21} to reduce to $G$ being a subgroup of $(\Gamma L_1(q^m) \circ \Gamma L_1(q^m)).2$ before doing the \textsc{Magma} computations.)

The remaining case $(m,t)=(2,2)$ is handled by the following result.

\begin{prop}\label{c72}
Suppose $G < GL(V)$ is an irreducible $p$-exceptional subgroup
satisfying the assumptions of Theorem $\ref{c7m}$ with $m = t = p = 2$. Then 
$G$ is conjugate to a subgroup of $GL_2(q^2)$ (and so 
is known by Lemma $\ref{lem:2dim}$). 
\end{prop}

\pf
Without loss of generality, we may assume that $G$ contains $Z = Z(GL(V))$. 
Consider the base subgroup $B = G \cap (GL(V_1) \circ GL(V_2))$ of index 
at most 2 in $G$. Then all $B$-orbits on $V$ have odd lengths. Hence, by 
Theorem \ref{c4lin}, either $|B|$ is odd, or $B$ is (conjugate to) the 
group $H$ appearing in the conclusion of that theorem.

Consider the latter case. We may identify $B$ with a subgroup of index $2$ in 
the group $\tilde{C}$ defined in the proof of Lemma \ref{yes2}. 
Adopt the notation of that proof. 
Since $G$ normalizes $C = O_{2'}(B)$, 
$G$ permutes the two $C$-orbits $U_1^\sharp$ and 
$U_2^\sharp$ of length $q^2-1$. Since $G$ is irreducible, $G$ cannot
fix either of them. Thus $G$ interchanges them, and so $G$ has an orbit 
of length $2(q^2-1)$, contradicting the $2$-exceptionality of $G$.

Hence $|B|$ is odd. Thus $|G|$ is not divisible by $4$ and so $G$ is soluble. Also 
$G$ is irreducible on $V$. Hence by the Fong-Swan theorem, the dimension $d$ of 
$V$ over $\End_G(V) \supseteq \FQ$ cannot be $4$ (but divides $4$), and so 
it is either $1$ or $2$. If $d = 1$, then $|G|$ divides $|GL_1(q^4)|$ and so
it is odd, a contradiction. Thus $d = 2$ and $G \leq GL_2(q^2)$ as in the conclusion.
\hal

\vspace{4mm}
This completes the proof of Proposition \ref{tensoregs}.

\section{Subfields}\label{c5sec}

It turns out that a $p$-exceptional group $H\leq \Ga L_n(q)$ cannot be realisable modulo scalars over a proper subfield $\F_{q_0}$ of $\F_q$.
Such groups are conjugate to subgroups of $Z\circ GL_n(q_0).\la\phi\ra$, for some proper subfield $\F_{q_0}$ of $\F_q$, where $\phi$ generates the group of field automorphisms of $GL_n(q)$ and $Z = Z(GL_n(q))$. We prove

\begin{thm}\label{sub}
Let $V = \Fq^n$, and $q = q_0^s$ with $s>1$. Suppose $H \le (Z \circ GL_n(q_0))\la \phi \ra < \Gamma L(V)$, 
where $\phi$ generates the group of field automorphisms of $GL_n(q)$. Then $H$ is not $p$-exceptional.
\end{thm}

\pf
Suppose that $H$ is $p$-exceptional on $V=\F_q^n$. By Lemma~\ref{normalsubgroups}, $ZH$ is also $p$-exceptional, so we may assume that $Z\subseteq H$. Then $H\cap GL_n(q)=Z\circ H_0$, where $H_0=H\cap GL_n(q_0)$, and note that $H_0$ is normal in $H$.
Let $\{b_1\dots,b_s\}$ be a basis for $\F_q$ as an $\F_{q_0}$-vector space. Then each $v\in V$ may be written as a sum $\sum_{i=1}^s b_i v_i$ with $v_i\in V_0=\F_{q_0}^n$.  If the non-zero $v_i$ are linearly independent in $V_0$ then the stabiliser $(H_0)_v$ must fix each of the $v_i$. 

Suppose first that $p$ divides $|H_0|$, or equivalently that $p$ divides $|H\cap GL_n(q)|$. Then $H_0$ is $p$-exceptional on $V$, and hence also on $V_0$, by Lemma~\ref{normalsubgroups}.
Let $x\in H_0$ have order $p$. Then there exists a 2-dimensional $x$-invariant $\F_{q_0}$-subspace $\la v_1,v_2\ra_{\F_{q_0}}$ of $V_0$ which is not fixed pointwise by $x$. 
Let $v=b_1v_1+b_2v_2$.   
Then as we observed in the previous paragraph, $(H_0)_v$ fixes both $v_1$ and $v_2$ and hence leaves the $\F_q$-space $U:= \la v_1,v_2\ra_{\F_{q}}$ invariant, fixing it pointwise. Moreover $(H_0)_v=(H_0)_{v_1,v_2}$ is the kernel of the action of $(H_0)_U$ on $U$.
It follows that $x\in (H_0)_U\setminus (H_0)_v$ and hence that $p$ divides $|(H_0)_U:(H_0)_v|$, which divides $|H_0:(H_0)_v|$, a contradiction to $p$-exceptionality.

Thus $|H_0|$ is coprime to $p$, and $H$ has a normal subgroup which contains $H\cap GL_n(q)$ as a subgroup of index $p$. By Lemma~\ref{normalsubgroups}, this normal subgroup is $p$-exceptional and so,  without loss of generality, we may assume that $H=(H\cap GL_n(q)).\la x\ra=(Z\circ H_0).\la x\ra$, where $x$ is a field automorphism of order $p$, and so $x$ induces a (possibly trivial) field automorphism on $GL_n(q_0)$. Thus each element of $H$ has the form $x^ich$, for some $i$, with $c\in Z$ and $h=(a_{ij})\in GL_n(q_0)$. Recall that $q=q_0^s$, and note that $x$ fixes pointwise (at least) an $\F_p$-subspace $\F_p^n$ of $V_0$. Let $q_1:=q^{1/p}$, so $\F_{q_1}$ is the fixed  field of $x$.


Choose $v_1, v_2$ linearly independent vectors in $V_0$ fixed by $x$, let $U:=\la v_1,v_2\ra_{\F_q}$, and define
\[
X(U):=\{c_1v_1+c_2v_2\,|\,c_1, c_2\in \F_q^*, c_1c_2^{-1}\not\in\F_{q_0}\}. 
\]
First we show that (i) for $v\in X(U)$, the stabiliser $H_v\subseteq H_U$, and (ii) $X(U)$ is $H_U$-invariant. Let $v=c_1v_1+c_2v_2\in X(U)$ and extend $\{v_1,v_2\}$ to a basis $\{v_1,v_2,\ldots,v_n\}$ for $V$ over $\F_q$.

(i) Let $x^ich\in H_v$ with $c,h$ as in the previous paragraph, with $h$ represented with respect to the basis $\{v_1,\ldots,v_n\}$.
Since $x$ and $Z$ leave $U$ invariant it is sufficient to show that $h$ does also. Now 
\[
v=v^{x^ich}=c\sum_{j=1}^n (c_1^{x^i}a_{1j}+c_2^{x^i}a_{2j})v_j 
\]
and so, for each $j\geq3$, $c_1^{x^i}a_{1j}+c_2^{x^i}a_{2j}=0$ which implies that $(c_1c_2^{-1})^{x^i}a_{1j}=-a_{2j}$. Since $a_{1j}, a_{2j}\in\F_{q_0}$  while  $(c_1c_2^{-1})^{x^i}\not\in\F_{q_0}$, it follows that $a_{1j}=a_{2j}=0$. Thus $h$ leaves $U$ invariant, proving part (i). Also we have, for $j=1,2$, that  
$c_j=c(c_1^{x^i}a_{1j}+c_2^{x^i}a_{2j})$, 
whence 
\begin{equation}\label{eq:c}
c^{-1} = c_1^{-1}c_1^{x^i}a_{11}+c_1^{-1}c_2^{x^i}a_{21}=
c_2^{-1}c_1^{x^i}a_{12}+c_2^{-1}c_2^{x^i}a_{22}  
\end{equation}
and since $h$ is nonsingular, $a_{11}a_{22} - a_{12}a_{21}\ne0$.

(ii) Since $c_1^{x^i}(c_2^{x^i})^{-1} = (c_1c_2^{-1})^{x^i}\not\in\F_{q_0}$, and since $(cc_1)(cc_2)^{-1} = c_1c_2^{-1}\not\in\F_{q_0}$, it follows that $x^i$ and each element of $Z$ leaves $X(U)$ invariant. It remains to consider $v^h$ where $h=(a_{ij})\in H_U$. Now
\[
 v^h = (c_1a_{11}+c_2a_{21})v_1 + (c_1a_{12}+c_2a_{22})v_2.
\]
If the coefficient of $v_2$ were 0 we would have $a_{22}=-a_{12}(c_1c_2^{-1})\in\F_{q_0}$ and hence $a_{22}=a_{12}=0$. However $a_{1j}=a_{2j}=0$ for each $j\geq3$, and this would imply that $h$ is singular, a contradiction. A similar argument shows that  the coefficient of $v_1$ is also nonzero. Suppose for a contradiction that $v^h\not\in X(U)$. Then
\[
d:=(c_1a_{11}+c_2a_{21})(c_1a_{12}+c_2a_{22})^{-1}\in\F_{q_0}^* 
\]
and we have $c_1(a_{11}-da_{12}) = c_2(da_{22}-a_{21})$. If $da_{22}-a_{21}\ne0$ then $c_2c_1^{-1}=(da_{22}-a_{21})^{-1}(a_{11}-da_{12})\in\F_{q_0}$ which is a contradiction. Thus $a_{21}=da_{22}$, and hence also $a_{11}=da_{12}$. This again implies that $h$ is singular, and finally we conclude that $v^h\in X(U)$, proving (ii).

We have $\la x\ra\leq H_U$, and since the $p$-part of $|H|$ is $p$, it follows that $\la x\ra$ is a Sylow $p$-subgroup of $H_U$. 
By $p$-exceptionality, $|H_v|$ is divisible by $p$ and hence, by Sylow's Theorem, it follows that $x$ must fix at least one vector of $X(U)$. Without loss of generality we may assume that $x$ fixes $v$.
Now observe that $v=v^x$ if and only if $c_j=c_j^x$ for $j=1,2$,  and hence $(c_1c_2^{-1})^x=c_1c_2^{-1}\in\F_{q_1}\setminus\F_{q_0}$. Thus $\F_{q_1}\not\subseteq\F_{q_0}$. 

Now we consider a special element of $X(U)$, namely $v'=v_1+bv_2$, where $b$ is a primitive element of $\F_{q_1}$ (which we have just shown does not lie in $\F_{q_0}$).
Suppose that $x^ich$ fixes $v'$, and note that $b^{x^i}=b$ since $b$ lies in the fixed field $\F_{q_1}$ of $x$. Then equation (\ref{eq:c}) in the computation in (i) shows that 
\[
c^{-1} =a_{11}+ba_{21}=b^{-1}a_{12}+a_{22}
\]
which implies that $b^2a_{21} +b(a_{11}-a_{22})-a_{12}=0$. Since $b\not\in\F_{q_0}$, it follows that $b$ has minimal polynomial over $\F_{q_0}$ of degree 2. The extension field $\F_{q_0}(b)$ contains the maximal subfield $\F_{q_1}$ of $\F_q$ as a proper subfield, and hence  $\F_{q_0}(b)=\F_q$ and $s=2$. Also, since $\F_{q_1}\not\subseteq\F_{q_0}$, $p$ is odd.

We may therefore write $q=r^{2p}, q_0=r^p, q_1=r^2$, and then $\F_{q_1}\cap\F_{q_0}=\F_r$.
Now $|X(U)|=(q-1)(q-q_0)=r^p(r^{2p}-1)(r^p-1)$. We showed above that the vectors of $X(U)$  fixed by $x$ are $c_1v_1+c_2v_2$ with each $c_i\in\F_{q_1}$ but $c_1c_2^{-1}\not\in\F_{q_0}$. Thus $x$ fixes precisely $(r^2-1)(r^2-r)$ vectors in $X(U)$. Since each Sylow $p$-subgroup of $H_U$ fixes the same number of vectors in $X(U)$, the number of Sylow $p$-subgroups of $H_U$ is at least $\frac{r^p(r^{2p}-1)(r^p-1)}{r(r^{2}-1)(r-1)}>r^{p-1}r^{2p-2}r^{p-1}=r^{4(p-1)}$.

We now consider the induced group $H_U^U\leq\Ga L(U)$. Since $|H\cap GL(V)|$ is coprime to $p$, and since $x$ acts nontrivially on $U$, it follows that $H_U^U$ is of the form
$(Z^U\circ L).\la x^U\ra$ with $L$ a $p'$-subgroup of $GL_2(q_0)$, where $Z^U$ is the group induced by $Z$. Note that the normaliser $N$ of $\la x^U\ra$ in $H_U^U$ intersects $Z^U$ in a subgroup of order $r^2-1$. Thus $|H_U^U:NZ^U|= |H_U^U:N|/(\frac{r^{2p}-1}{r^2-1})$, which is at least 
$\frac{r^p(r^p-1)}{r(r-1)}> r^{2p-2}$. Since $NZ^U$ contains $Z^U.\la x^U\ra$, we have $|L/(L\cap Z^U)|\geq |H_U^U:NZ^U|>r^{2p-2}\geq 3^4$. It then follows from the classification of subgroups of $PGL_2(q_0)$, that the $p'$-group
$L/(L\cap Z^U)$ is contained in a dihedral group $D_{2(q_0-1)}$ or $D_{2(q_0+1)}$. Thus 
the number of Sylow $p$-subgroups of $H_U^U$ is at most the number of Sylow $p$-subgroups of $H_U^U$ in the case where $L/(L\cap Z^U)= D_{2(q_0\pm1)}$. In this case, the normaliser of $\la x^U\ra$ in $H_U^U$ is $Z_{r^2-1}\cdot D_{2(r\pm1)}\cdot\la x^U\ra$, and so the number of Sylow $p$-subgroups is 
$\frac{(r^{2p}-1)(r^p\pm1)}{(r^2-1)(r\pm1)}<4r^{3(p-1)}$, which is less than the lower bound $r^{4(p-1)}$. This contradiction completes the proof. \hal

%
%

\section{Extraspecial type normalizers: $\mathcal{C}_6$ subgroups}\label{c6sec}

Let $r$ be a prime, $m$ a positive integer, and let $R$ be an $r$-group of symplectic type such that $|R/Z(R)| = r^{2m}$, $R$ is of exponent $r\cdot (2,r)$, and $R$ is as in Table 1. Let $V = V_d(q)$ be a faithful, absolutely irreducible $\F_qR$-module, where $r$ does not divide $q$. Then $d = \dim V = r^m$, and $N_{GL(V)}(R)$ is as in the table. Assume further that $R$ is not realised over a proper subfield of $\F_q$. Then $q$ is a minimal power of the characteristic $p$, subject to the conditions in the last column of the table.

Here we prove

\begin{thm}\label{c6case}
Let $r$ be a prime, and assume that $R$ and $V=V_d(q)$ are as above. 
Suppose $R \trianglelefteq G \le N_{\Gamma L(V)}(R)$ and $G$
is $p$-exceptional, and is not transitive on $V^\sharp$. Then 
$G$ is imprimitive on $V$, and one of the following holds.
\begin{enumerate}
\item[\rm (i)] $r=2,q=3,d=4$ and $G = 2^{1+4}\AAA_4$ or $2^{1+4}\SSS_4$, with orbits
on vectors of sizes $1,16,64$.
\item[\rm (ii)] $r=3,q=4,d=3$: here there are five $2$-exceptional groups of the form $3^{1+2}.X$; they are
$3^{1+2}.2$, $3^{1+2}.6$, $3^{1+2}.\SSS_3$ (two such groups), and $3^{1+2}.D_{12}$; the first has orbit sizes $1,9^4,27$, the rest have $1, 9, 27^2$.
\item[\rm (iii)] $r=2,q=3,d=8$: here there are five $3$-exceptional groups; they are 
 $2^{1+6}_+.X$ with $X = L_3(2)$, $2^3.L_3(2)$ or $2^3.7.3$
(all with orbit sizes $1$, $16$, $112$, $128^2$, $224$, $448$, $896$, $1024$, $1792^2$),
and $2^{1+6}_-.Y$ with $Y=2^4.\AAA_5$ or $2^4.\SSS_5$
(all with orbit sizes  $1$, $160$, $1280$, $5120$).
\end{enumerate}
\end{thm}

\begin{table}
\caption{}
\label{c6poss}
\[
\begin{array}{|l|l|l|l|}
\hline
R & d & N_{GL(V)}(R)/RZ & q \\
\hline
r^{1+2m},\,r \hbox{ odd} & r^m & Sp_{2m}(r) & q \equiv 1 \hbox{ mod }r \\
4 \circ 2^{1+2m} & 2^m & Sp_{2m}(2) &  q \equiv 1 \hbox{ mod }4 \\
2^{1+2m}_{\pm} &2^m & O^\pm_{2m}(2) &  q \equiv 1 \hbox{ mod }2 \\
\hline
\end{array}
\]
\end{table}

\subsection{Reductions}

Let $G$ be a $p$-exceptional group as in the hypothesis of Theorem \ref{c6case}, and assume that $G \le GL(V)$. We shall handle the case where $G \le \Ga L(V)$ and $G \not \le GL(V)$ at the end of the proof in Section \ref{semilinc6}.

We begin with a technical lemma concerning the Jordan block structure of certain elements of $N_{GL(V)}(R)$. Write $J_k$ for a unipotent Jordan block of size $k$.

\begin{lemma}\label{jordan}
Let $R$ be as in Table $\ref{c6poss}$.
Assume that the characteristic $p$ is a primitive prime divisor of $r^{2m}-1$ 
or a primitive prime divisor of $r^m-1$ with $m$ odd, and also that $p=7$ 
when $(r,m) = (2,3)$. Let $t$ be an element of order $p$ in $N_{GL(V)}(R)$.
Then $t$ acts on $V$ as $(J_p^k,J_\ell)$, where $r^m = kp+\ell$ and 
$0\le \ell < p$.
\end{lemma}

\pf 
First we consider the case where $m$ is odd and 
$p$ is a primitive prime divisor of $r^m-1$ (so $\ell=1$). If $R = 2^{1+2m}_{\pm}$, 
then this in particular implies that $R = 2^{1+2m}_{+}$. In all cases,
embedding $tRZ$ in a subgroup $GL_m(r)$ of $N_{GL(V)}(R)/RZ$, 
one can check that there exists a $t$-stable elementary 
abelian $r$-subgroup $A < R$ of order $r^m$ such that $V\downarrow A$ affords the 
regular representation of $A$ and moreover $t$ acts fixed-point-freely on the 
nontrivial irreducible characters of $A$. Thus $t$ also permutes 
fixed-point-freely the nontrivial $A$-eigenspaces in $V$. Since $|t| = p$ and 
$\dim C_V(A) = 1$, it follows that $t$ acts on $V$ as $(J_p^k,J_1)$ as stated. 

Assume now that $p$ is a primitive prime divisor of $r^{2m}-1$ (so $\ell=p-1$).
Note that $Z(R)Z$ acts trivially on $V \otimes V^*$, $B := R/Z(R)$ is 
elementary abelian of order $r^{2m}$, $(V \otimes V^*) \downarrow B$ affords the 
regular representation of $B$, and $t$ acts fixed-point-freely on the 
nontrivial irreducible characters of $B$. Thus $t$ also permutes 
fixed-point-freely the nontrivial $B$-eigenspaces in $V \otimes V^*$. As 
before, it follows that $t$ acts on $V \otimes V^*$ as $(J_p^n,J_1)$ with 
$n := (r^{2m}-1)/p$. 

Recall \cite[Theorem VIII.2.7]{Fe} that 
the Jordan canonical form of 
$J_a \otimes J_b$ equals 
$$(J_{a+b-1},J_{a+b-3},\ldots ,J_{b-a+1})$$
if $1 \leq a \leq b < a + b \leq p$, and
$$(J_p^{(a+b-p)},  J_{2p-a-b-1},J_{2p-a-b-3}, \ldots ,J_{b-a+1})$$
if $1 \leq a \leq b < p < a+b$. It follows, that  
$J_a \otimes J_b$ contains no block of size between
$2$ and $p-1$ only when $a = b \in \{1,p-1\}$. Applying this observation 
to $(V\otimes V^*) \downarrow t$, 
we see that $V \downarrow t = (J_p^c,J_e^d)$ for some $c,d$ and some 
$e \in \{1,p-1\}$. In fact, since $(V \otimes V^*) \downarrow t$ contains 
$J_1$ only once, $d \leq 1$. But $pc + de = r^m = kp+(p-1)$, whence 
$(d,e) = (1,p-1)$, as stated. 
\hal 

\vspace{2mm}

The next lemma reduces the number of possibilities for $R$ to a finite number.

\begin{lemma}\label{finite}
The possibilities for $R$ are as follows:
\[
\begin{array}{l}
4\circ 2^{1+2m},\;m\le 10 \\
2^{1+2m}_{\pm},\;m\le 11 \\
3^{1+2m},\;m\le 5 \\
5^{1+2m},\; m\le 2 \\
7^{1+2m},\; m\le 2 \\
11^{1+2} 
\end{array}
\]
\end{lemma}

\pf First consider $R = 4\circ 2^{1+2m}$. Assume $m\ge 11$. Here $N_{GL(V)}(R)/RZ = Sp_{2m}(2)$ and $G \le N_{GL(V)}(R)$. 
Let $t \in G$ be an element of order $p$. By \cite[4.3]{gusa}, there are $m+3$ conjugates of $t$ which generate a subgroup of 
$N_{GL(V)}(R)$ covering $Sp_{2m}(2)$. In fact these conjugates generate the whole of $N_{GL(V)}(R)$ since otherwise they would generate a covering group of $Sp_{2m}(2)$; but there is no such group in dimension $2^m$ by \cite{LS}.
It follows that $\dim C_V(t) \le 2^m(1-\frac{1}{m+3})$. Hence, as $G$ is $p$-exceptional, Lemma \ref{bound}(ii) implies that
\begin{equation}\label{largebd}
q^{2^m} \le  q^{2^m(1-\frac{1}{m+3})}|t^G|.
\end{equation}
It follows that $ q^{2^m/(m+3)} \le |t^G| \le 2^{2m}|Sp_{2m}(2)|$. Since $q\ge 5$ in this case (as it is 1 mod 4), this is a contradiction for $m \ge 11$.

An entirely similar argument shows that $m\le 11$ in the case where $R = 2^{1+2m}_{\pm}$.

Now consider $R = r^{1+2m}$ with $r$ odd. As above, take an element $t$ of 
order $p$ in $G$. By \cite{gusa}, there are $m+3$ conjugates of $t$ which 
generate a subgroup of $ N_{GL(V)}(R)$ which covers $Sp_{2m}(r)$. Such a 
subgroup fixes no nonzero vectors in $V$ (note that the restriction of $V$ to a 
subgroup $Sp_{2m}(r)$ is the sum of two irreducible Weil modules -- 
see \cite{ward}), and so again 
$\dim C_V(t) \le r^m(1-\frac{1}{m+3})$, giving
\begin{equation}\label{use}
 q^{r^m/(m+3)} \le |t^G| \le r^{2m}|Sp_{2m}(r)|.
\end{equation}
Also $r$ divides $q-1$. The bound (\ref{use}) implies that $R$ is as in the conclusion, except that $R = 13^{1+2}$ is also possible; but this can be ruled out by noting that only 3 conjugates of $t$ are required (rather than 4), by \cite[3.1]{gusa}. \hal

\begin{lemma} We have $r < 5$.
\end{lemma}

\pf Suppose $r\ge 5$. Then $r=5,7$ or 11 by Lemma \ref{finite}.

First consider $r=5$. Here $m\le 2$. Suppose $m=1$, so $G/Z \le 5^2.Sp_2(5)$, $d = \dim V = 5$  and $q \equiv 1 \hbox{ mod }5$. As $p$ divides $|G|$ we have $p=2$ or 3 and $q = 16$ or 81. For $p=2$, Lemma \ref{jordan} shows that an involution $t \in G$ acts on $V$ as $(J_2^2,J_1)$, so that $\dim C_V(t) = 3$. Hence Lemma \ref{bound}(ii) gives $16^2 \le |t^G|$. However, 
$|t^G|  \le 5^2$, so this is a contradiction. And for $p=3$, an element $t \in G$ of order 3 acts on $V$ as $(J_3,J_2)$ and we argue similarly.

Now suppose $m=2$ (still with $r=5$). Here  $G/Z \le 5^4.Sp_4(5)$, $d = \dim V = 25$  and $q \equiv 1 \hbox{ mod }5$.
As $p$ divides $|G|$ we have $p = 2,3$ or 13 and $q=16$, 81 or $13^4$. Let $t\in G$ be an element of order $p$. For $p=2$, involutions in $G$ lie in a subgroup $5^{1+2}Sp_2(5) \circ 5^{1+2}Sp_2(5)$ acting on $V$ as a tensor product of two 5-spaces, and hence act on $V$ as 
either $(J_2^2,J_1) \otimes (J_2^2,J_1)$ or $(J_2^2,J_1) \otimes I_5$; hence $\dim C_V(t)\le 15$. For $p=3$, elements of order 3 in $G$ also lie in this subgroup, so act as $(J_3,J_2) \otimes (J_3,J_2)$ or $(J_3,J_2) \otimes I_5$; hence $\dim C_V(t)\le 10$. And for $p=13$, elements of order 13 in $G$ act as $(J_{13},J_2)$, hence $\dim C_V(t) =2$. We conclude that 
$\dim C_V(t) \le 15$ in any characteristic. Now Lemma \ref{bound}(ii) gives a contradiction.

The cases $r=7$ and $r=11$ are ruled out in entirely similar fashion. \hal

\subsection{The case $r=2$}

In this section we handle the case where $r=2$. By Lemma \ref{finite}, $R$ is either $4\circ 2^{1+2m}$ with $m \le 10$, or $2^{1+2m}_{\pm}$ with $m\le 11$.

\begin{lemma}\label{downto5}
We have $m\le 5$.
\end{lemma}
 
\pf Consider the case $m=6$. Here $d = \dim V = 64$, $R = 2^{1+12}_{\pm}$ or $4\circ 2^{1+12}$ and $G/RZ \le O^\pm _{12}(2)$ or $Sp_{12}(2)$, respectively. 
Moreover $p = 3,5,7,11, 13, 17$ or 31.

Suppose first that $p=3$ and $R = 2^{1+12}_{\pm}$. Then $q=3$. We shall use Lemma \ref{bound}(ii). Let $t \in G$ have order 3. Modulo $R$, $t$ is conjugate to an element $t_k$ of order 3 in a subgroup $O_2^-(2)^k$ of $O_{12}^\pm (2)$, where $1\le k \le 6$, projecting nontrivially to each factor. We can work out the Jordan form of $t_k$ on $V$ as an element of $2^{1+2}O_2^-(2) \otimes \cdots \otimes 2^{1+2}O_2^-(2) \otimes I_{2^{6-k}}$: so on $V$, $t_k$ acts as $J_2 \otimes \cdots \otimes J_2 \otimes I_{2^{6-k}}$ where there are $k$ factors $J_2$. The Jordan forms of these tensor products are easily worked out, and we find that the number of Jordan blocks of $t_k$ is as follows:
\[
\begin{array}{l|llllll}
k & 1&2&3&4&5&6 \\
\hline
\dim C_V(t_k) & 32 & 32 & 24 & 24 & 22 & 22 
\end{array}
\]
By Corollary \ref{vnormgen}, modulo $R$ the centralizer of $t_k$ in $R.O_{12}^\pm (2)$ is
$C_{O_{12}^-(2)}(t_k)$, which contains $O_{12-2k}(2)\times 3^k.\SSS_k$. Hence
\[
|t_k^G| \le \frac{2^{12}|O_{12}^\pm (2)|}{2^{12-2k}|O_{12-2k}(2)|\cdot 3^k\cdot k!}.
\]
For each $k$ between 1 and 6, we check that $|V:C_V(t_k)| = 3^{64-\dim C_V(t_k)}$ is greater than $|t_k^G|$. So this contradicts Lemma \ref{bound}(ii).

When $p=3$ and $R=4\circ 2^{1+12}$ we have $q=9$ (as $q \equiv 1 \hbox{ mod }4$ in this case). Here an element $t \in G$ of order 3 lies in a subgroup $Sp_2(2)^k$ of $Sp_{12}(2)$ (modulo $R$) for some $k$ with $1 \le k \le 6$, and we argue as above.

Now suppose $p=5$. Here $q=5$ for both cases $R =  4\circ 2^{1+12}$ and $R=2^{1+10}_{\pm}$. In the former case, an element $t \in G$ of order 5 is conjugate to an element $t_k$ in a subgroup $Sp_4(2)^k$ of  $Sp_{12}(2)$, where $1\le k \le 3$, projecting nontrivially to each factor. Then $t_k$ acts on $V$ as an element of a tensor product of $k$ factors $4\circ 2^{1+4}.Sp_4(2)$ acting in dimension 4. By Lemma \ref{jordan}, an element of order 5 in such a factor acts as $J_4$, so $t_k$ acts on $V$ as 
$J_4 \otimes \cdots \otimes J_4 \otimes I_{64/4^k}$. Working out these tensor products of Jordan blocks, it follows that 
$\dim C_V(t_k) = 16, 16, 13$ according as $k=1,2,3$ respectively. Moreover as above 
\[
|t_k^G| \le  \frac{2^{12}|Sp_{12}(2)|}{2^{12-4k}|Sp_{12-4k}(2)|\cdot 5^k\cdot k!}.
\]
It now follows that $|V:C_V(t_k)| = 5^{64-\dim C_V(t_k)}$ is greater than $|t_k^G|$ for each $k$, contradicting Lemma \ref{bound}(ii) again. The $R=2^{1+10}_{\pm}$ case is entirely similar.

For $p=7$, we see as above that an element $t\in G$ of order 7  is conjugate to an element $t_k$ ($k = 1$ or 2) 
acting on $V$ as $(J_7,J_1) \otimes I_8$ or $(J_7,J_1) \otimes (J_7,J_1)$
(using Lemma \ref{jordan} to see that the action of an element of order 7 in $4 \circ 2^{1+6}.Sp_6(2) < SL_8(7^a)$ is $(J_7,J_1)$). Hence $\dim C_V(t_k) = 16$ or 10 and we contradict Lemma \ref{bound} as before. 

For larger values of $p$ there is only one class of elements of order $p$ in $N_{GL(V)}(R)$, and its action on $V$ can be computed as above using Lemma \ref{jordan}. We find that for $t \in G$ of order $p$, $\dim C_V(t)$ is $6, 5, 4$ or 4 according as $p = 11,13,17$ or 31. In each case $|V:C_V(t)|$ is much bigger than $|N_{GL(V)}(R)|$, contradicting Lemma \ref{bound} once more.

This completes the argument for $m=6$. For $m\ge 7$ the same method applies. \hal

\begin{lemma}\label{m5}
$m$ is not $5$.
\end{lemma}

\pf Suppose $m=5$. The method is very similar to the previous proof, but the bounds are tighter and more work is needed for the case $q=3$. 

Consider first the case where $p=3$. Let $R = 2^{1+10}_{\pm}$. Then $q=3$.  If $t \in G$ is an element of order 3, then $t$ is conjugate to an element $t_k$ lying in a subgroup $O_2^-(2)^k$ of $O_{10}^\pm (2)$ and acting on $V$ as $J_2\otimes \cdots \otimes J_2 \otimes I_{32/2^k}$. We compute that $\dim C_V(t_k) = 16, 16, 12, 12, 11$ for $k=1,2,3,4,5$ respectively. If $k=1$ then Lemma \ref{bound}(ii) gives $3^{16} < |t_1^G|$; however $|t_1^G| \le |O_{10}(2) : 2^8.(3\times O_8(2)|$ which is less than $3^{16}$. 
Hence (again by Lemma \ref{bound}(ii)), $G/R$ is a subgroup of $O_{10}^\pm (2)$ containing no conjugates of $t_1$ and 
at least $3^{32-\dim  C_V(t_k)}/2^{2k}$ conjugates of $t_k$ for some $k \ge 2$. However a \textsc{Magma} computation shows that there is no such subgroup. 

In the case where $p=3$ and $R =  4\circ 2^{1+12}$ we have $q=9$. Here an element $t \in G$ of order 3 is conjugate to some $t_k$ as in the previous paragraph but no \textsc{Magma} computation is required as $9^{32-\dim  C_V(t_k)} > |t_k^G|$ for all $k$, which contradicts Lemma \ref{bound}.

Now consider $p=5$. Here $q=5$ and an element $t \in G$ of order 5 is conjugate to an element $t_k$ ($k=1$ or 2) lying in a subgroup $O_4^-(2)^k$ of $O_{10}(2)$ or $Sp_{10}(2)$ and acting on $V$ as $J_4 \otimes I_8$ (for $k=1$) or $J_4\otimes J_4 \otimes I_2$ (for $k=2$). So $\dim C_V(t_k) = 8$ for $k=1,2$. However we check as above that $5^{24} > |t_k^G|$ for $k=1,2$, which contradicts Lemma \ref{bound}(ii).

The other possible values of $p$ are $7,11,17$ and 31. For these values there is only one class of elements of order $p$ in $N_{GL(V)}(R)$, and its action on $V$ can be computed as above using Lemma \ref{jordan}. In each case $|V:C_V(t)|$ is much bigger than $|N_{GL(V)}(R)|$, contradicting Lemma \ref{bound} once more. \hal

\begin{lemma}\label{m4}
$m$ is not $4$.
\end{lemma}

\pf Suppose $m=4$. Then $p$ divides $|Sp_8(2)|$, so is $3,5,7$ or 17. The cases $p=7$ or 17 are easily handled as in the last paragraph of the previous proof.

Consider now the case $p=5$. Here $q=5$ and $G \le (4\circ 2^{1+8}).Sp_8(2) < GL_{16}(5)$. An element $t \in G$ of order 5 is conjugate to an element $t_k$ ($k=1$ or 2) lying in $Sp_4(2)^k$ and acting on $V$ as $J_4\otimes I_4$ or $J_4\otimes J_4$; so  $\dim C_V(t_k) = 4$ for both $k=1,2$. Hence $5^{12}<|t^G|$ by Lemma \ref{bound}(ii). It follows that $G/R$ contains no conjugates of $t_1$ and contains at least $5^{12}/2^8$ conjugates of $t_2$. Using \cite{Atlas}, one checks that the only possible maximal subgroup of $Sp_8(2)$ containing $G/R$ is $O_8^+(2)$, and then that the only subgroup of this containing enough conjugates of $t_2$ are $\Om_8^+(2)$ and $O_8^+(2)$. But these contain conjugates of $t_1$, a contradiction.

Now suppose $p=3$. Consider first the case where $R = 4\circ 2^{1+8}$. Here $q=9$. An element $t \in G$ of order 3 is conjugate to some $t_k$ ($1 \le k \le 4$ lying in $Sp_2(2)^k$ and acting on $V$ as $J_2\otimes \cdots \otimes J_2 \otimes I_{16/2^k}$, and we find that $\dim C_V(t_k) = 8,8,6,6$ according as $k = 1,2,3,4$ respectively. We compute that 
$9^{\dim V - \dim C_V(t_k)} > |t_k^G|$ for $k=1,3,4$. Hence $G/R$ is a subgroup of $Sp_8(2)$ containing no conjugates of $t_1,t_3,t_4$, and at least $9^8/2^4$ conjugates of $t_2$. One checks using \cite{Atlas} that there are no such subgroups.

Finally suppose $p=3$ and $R = 2^{1+8}_\e $, so $G \le 2^{1+8}_\e . O_8^\e (2) < GL_{16}(3)$. Here the usual bounding methods do not work and we use \textsc{Magma} computation along the following lines. 
For $\e = -$, we compute that the group $R$ has a regular orbit on vectors that has 48960 images under 
$R.O_8^-(2)$. Moreover, $O_8^-(2)$ has 4683 subgroups of order 
divisible by 3, and all of them have an orbit of length divisible by 3 in 
the action on 48960 points.
Similarly, for $\e=+$ the group $R$ has an orbit of length 128 that has 1575 images under 
$R.O_8^+(2)$, and $O^+_8(2)$ has 5988 subgroups with order 
divisible by 3;  all of them have an orbit of length divisible by 3 on 
the set of size 1575. \hal

\begin{lemma}\label{m3}
If $m=3$, then $q=3$ and $G$ is as in conclusion {\rm (iii)} of Theorem $\ref{c6case}$.
\end{lemma}

\pf Suppose $m=3$. Then $p$ divides $|Sp_6(2)|$, so is $3,5$ or 7. We make heavy use of computation in this proof.

Let $p=3$. First consider $R = 2^{1+6}_+$, so $G \le 2^{1+6}_+.O_6^+(2) < GL_8(3)$. 
We compute that the group $R$ has an orbit of length 16 on vectors,  that has 30 images under 
$X:=R.O_6^+(2)$. 
Let $\D$ be this set of 30 images, and let $G = R.Y \le X$ with $Y\le O_6^+(2)$.
Then $|Y|$ is divisible by 3. If $Y$ has an orbit $\D_0$ on $\D$ of length divisible by 3, then $G$ has 
an orbit on vectors of length $16|\D_0|$, contrary to 3-exceptionality.
We compute that the group $O_6^+(2)$ has 176 subgroups of order 
divisible by 3, and all but 12 have an orbit of length divisible by 3 in 
the action on $\D$. Hence $Y$ is one of the remaining 12 groups.  
When pulled back to subgroups of $X$ containing $R$, 
all but three of these have an orbit on nonzero vectors of length 
divisible by 3. The three remaining groups are 3-exceptional -- they are
$R.L_3(2)$, $R.2^3.L_3(2)$ and $R.2^3.7.3$. All three are imprimitive on $V$, preserving a decomposition into eight 1-dimensional spaces, and are as in Theorem \ref{c6case}(iii).

Now consider $R = 2^{1+6}_-$, so $G \le 2^{1+6}_-.O_6^-(2) < GL_8(3)$. 
Here $R$ has an orbit on vectors of length 32 that has 45 images under 
$X:=R.O^-(6,2)$.  Now $O^-(6,2)$ has 238 subgroups of order 
divisible by 3 and all but 7 have an orbit of length divisible by 3 in 
the action on 45 points.  Pulling these remaining 7 subgroups back to subgroups of $X$ 
containing $R$ we see that all but two have an orbit on nonzero vectors 
of length divisible by 3. The two groups are $R. (2^4.\AAA_5)$ and 
$R. (2^4.\SSS_5)$, as in Theorem \ref{c6case}(iii). Both are imprimitive on $V$ with a
decomposition into two 4-dimensional spaces.

Finally, if $R= 4\circ 2^{1+6}$, then $q=9$ and
$R$ has an orbit of length 32 that has 270 images under 
$X:=R.Sp(6,2)$. Then $Sp_6(2)$ has 2777 subgroups of order 
divisible by 3, and all but 13 of these have an orbit of length divisible 
by 3 in the action on 270 points. Pulling back these 13 subgroups to 
subgroups of $X$ containing $R$, we find that they all have an orbit of 
length divisible by 3 on vectors.

Now let $p=5$, so $G \le X:= 4\circ 2^{1+6}.Sp_6(2) < GL_8(5)$. 
Here $R$ has an orbit on vectors of length 32 that has 135 images under 
$X:=R. Sp_6(2)$. Now $Sp_6(2)$ has 82 subgroups with order 
divisible by 5 and all have an orbit of length divisible by 5 in the 
action on 135 points.

Lastly let $p=7$. If $q=49$ the usual method using Lemma \ref{bound} yields a contradiction, 
so suppose $q=7$ and $R = 2^{1+6}_+$ (note that 7 does not divide $|O_6^-(2)|$ so only the $+$ type is possible).
An element $t \in G$ of order 7 acts on $V$ as $(J_7,J_1)$, so Lemma \ref{bound}(ii) implies that $G/R$ is a subgroup of $O_6^+(2)$ containing at least $7^6/2^6$ conjugates of $t$. This implies that $G/R$ contains $\Om_6^+(2)$. 
Now one computes that $R.\Om_6^+(2)$ has an orbit on vectors of length divisible by 7.
 \hal

\begin{lemma}\label{m2}
If $m=2$ then $q=3$ and $G$ is a subgroup $2^{1+4}_-.\AAA_4$ or $2^{1+4}_-.\SSS_4$ of $GL_4(3)$. These subgroups are imprimitive and $3$-exceptional, with orbits on vectors of sizes $1,16,64$.
\end{lemma}

\pf Here $G$ is a subgroup of one of $2^{1+4}_\e.O_4^\e(2) < GL_4(3)$, or $4\circ 2^{1+4}.Sp_4(2) < GL_4(9)$ or $GL_4(5)$. A routine computation of all subgroups of $R.O_4^\e (2)$ and $R.Sp_4(2)$ containing $R$ gives the conclusion. \hal

\vspace{4mm}
This completes the proof of Theorem \ref{c6case} in the case where $G \le GL(V)$ and $r=2$.

\subsection{The case $r=3$}

Now suppose $r=3$. By Lemma \ref{finite}, we have $R = 3^{1+2m}$ with $m \le 5$. Also 
$G \le R.Sp_{2m}(3) < GL(V) = GL_{3^m}(q)$ with $q \equiv 1 \hbox{ mod }3$.

\begin{lemma}\label{mle2}
We have $m\le 2$.
\end{lemma}

\pf First consider $m=5$. The bound (\ref{use}) forces $q =4$ or 7. If $q=7$ and $t \in G$ is an element of order 7, then modulo $R$, $t$ lies in a subgroup $Sp_6(3) \times Sp_4(3)$ of $Sp_{10}(3)$. Hence by Lemma \ref{jordan}, $t$ acts on $V$ as $(J_7^{27},J_6^9)$ and $\dim C_V(t) = 36$. Now Lemma \ref{bound} gives a contradiction. Similarly, when $q=4$ an involution $t\in G$ lies in a subgroup $Sp_2(3)^k$ and acts as a tensor product of $k$ factors $(J_2,J_1)$ with an identity  matrix, whence we see that $\dim C_V(t) \le 162$, and again Lemma \ref{bound} is violated.

Now let $m=4$. Then $p$ divides $|Sp_8(3)|$. If $t\in G$ has order $p$, we see using Lemma \ref{jordan} in the usual way that $\dim C_V(t) \le 54$. Now Lemma \ref{bound}(ii) forces $p=2$ and $q=4$. The involution $t$ is conjugate to an element $t_k$ lying in $Sp_2(3)^k$ ($1\le k\le 4$) and acting on $V$ as a tensor product of $k$ factors $(J_2,J_1)$ with an identity matrix. Hence $\dim C_V(t_k) = 54,45,42,41$ according as $k=1,2,3,4$ respectively. We find that for each $k$ we have 
$|V:C_V(t)| > |t_k^G|$, contrary to Lemma \ref{bound}.

The case $m=3$ is dealt with in exactly the same fashion, and we leave this to the reader. \hal

\begin{lemma}\label{m12}
We have $m=1$, $q=4$ and $G = 3^{1+2}.2$ or $3^{1+2}.6$ in $GL_3(4)$, as in Theorem $\ref{c6case}${\rm (iii)}. Both are imprimitive on $V$, with orbit sizes $1,9^4,27$ or $1,9,27^2$, respectively.
\end{lemma}

\pf Suppose $m=2$. Then $p$ divides $|Sp_4(3)|$, so $p=2$ or $5$. The usual arguments rule out $p=5$ (for which $q=25$), so $p=2$, $q=4$ and $G \le R.Sp_4(3) < GL_9(4)$. A computation of all the subgroups of $Sp_4(3)$ of even order shows that there are no 3-exceptional groups in this case.

Hence $m=1$. Then $p$ divides $|Sp_2(3)|$, so  $p=2$, $q=4$ and a computation gives the examples in the lemma. \hal

\vspace{4mm}
This completes the proof of Theorem \ref{c6case} in the case where $G \le GL(V)$.

\subsection{The semilinear case}\label{semilinc6}

We now complete the proof of Theorem \ref{c6case} by handling the case where the $p$-exceptional group $G$ in the hypothesis  lies in $\Ga L(V)$ and not in $GL(V)$, where $V = V_d(q)$. Let $G_0 = G\cap GL(V)$. If $p$ divides $|G_0|$ then $G_0$ is $p$-exceptional, hence is given by the linear case of the theorem which we have already proved. The only possibility is that $d=3,q=4$ and $G_0$ is one of the two groups in the conclusion of Lemma \ref{m12}. Computation in $\Ga L_3(4)$ reveals one further 2-exceptional group $G$ in this case, the group $3^{1+2}.D_{12}$ in Theorem \ref{c6case}(iii).

Assume now that $p$ does not divide $|G_0|$. By Lemma \ref{pexc}, $G$ has a $p$-exceptional normal subgroup $G_0\la \s \ra$, where $\s \in \Ga L(V)\backslash GL(V)$ is a field automorphism of order $p$. Hence $q=p^{kp}$ for some integer $k$. If $r=2$ then
$q=p$ or $p^2$ (as $G$ is not realisable over a proper subfield of $\F_q$), and also $q$ is odd; so
it is impossible to have $q = p^{kp}$. Therefore $r\ge 3$.
Also the field automorphism $\s$ acts on $V$, fixing pointwise a subset $V_d(q^{1/p})$, and so $|V:C_V(\s)| = q^{(1-\frac{1}{p})r^m}$. Hence Lemma \ref{bound}(ii) gives
\begin{equation}\label{crude}
q^{(1-\frac{1}{p})r^m} \le |\s^G| \le r^{2m}|Sp_{2m}(r)|.
\end{equation}
If $p\ne 2$ then $q \ge 27$ and one checks that (\ref{crude}) cannot hold. Hence $p=2$ and $q=4^k$.

If $k\ge 2$ then (\ref{crude}) implies that $m=1$, $q=16$ and $r=5$. For $k=1$, we have $r=3$ and (\ref{crude}) gives $m\le 3$. Moreover if $m=3$ then $G_0/R$ is an odd order subgroup of $Sp_6(3)$. Computation shows that the largest such subgroup has order $3^7\cdot 13$.  Hence $|G_0/R| \le 3^7\cdot 13$, and so (\ref{crude}) gives $2^{27} \le 3^{14}\cdot 13$, which is false. Hence $m\le 2$. The possibilities remaining are as follows, where we write $K$ for the odd order group $G_0/R\le Sp_{2m}(r)$:
\[
\begin{array}{ll}
(1) & G = 5^{1+2}.K.2 < \Ga L_5(16) \\
 (2) & G = 3^{1+2}.K.2 < \Ga L_3(4) \\
(3) &  G = 3^{1+4}.K.2 < \Ga L_9(4)
\end{array}
\]
In cases (1) and (3) we check computationally that there are no 2-exceptional examples, and in case (2) we get the two  examples $3^{1+2}.\SSS_3$ in the conclusion of Theorem \ref{c6case}(iii). 

\vspace{4mm}
This completes the proof of Theorem \ref{c6case}.

\section{$\mathcal{C}_9$ case I: preliminaries and ${\rm Lie}(p)$}\label{c91}

Define $G \le \Ga L(V)=\Ga L_d(q)$ ($q=p^a$) to be in class $\mathcal{C}_9$ if $G/Z$ is almost
simple, with socle absolutely irreducible and not realisable over a proper
subfield of $\F_q$.  Let $L = {\rm soc}(G/Z)$, a simple group, and let $\hat L$ be a perfect preimage of $L$ in $G$. 

The case where the simple group $L$ is of Lie type in characteristic $p$ turns out to be very easy. Define ${\rm Lie}(p)$ to be the set of simple groups of Lie type in characteristic $p$, excluding $Sp_4(2)'$, $G_2(2)'$, $^2\!F_4(2)'$ in characteristic 2, and $^2\!G_2(3)'$ in characteristic 3. The first two of these and the last are dealt with in Sections \ref{c92} and \ref{c94} (in their guises as $\AAA_6$, $U_3(3)$ and $L_2(8)$) , and $^2\!F_4(2)'$ in the remark below after Corollary \ref{liepcor}. 

\begin{lemma}\label{liep}
Suppose that $G \le \Ga L(V)=\Ga L_d(q)$ ($q=p^a$) is in class $\mathcal{C}_9$ and $L = {\rm soc}(G/Z) \in {\rm Lie}(p)$. If $G$ is $p$-exceptional, then $G$ (and $\hat L$) is transitive on $P_1(V)$.
\end{lemma}

\pf From the structure of the exceptional Schur multipliers of simple groups in 
${\rm Lie}(p)$ (see \cite[Table 5.1.D]{KL} for example), 
we see that the perfect group 
$\hat L$ must be an image of the simply connected group of type $L$ (since the 
extra parts of such multipliers are always $p$-groups which hence act trivially on 
irreducible modules in characteristic $p$). 
Therefore it follows from \cite[Theorem 4.3(c)]{curtis} 
that a Sylow $p$-subgroup $P$ of $\hat L$ fixes a unique $1$-space in $V$. If $G$ 
is $p$-exceptional then so is $\hat L$, and so every nonzero vector in $V$ is 
fixed by an $\hat L$-conjugate of $P$. For 1-spaces $\la v \ra$, $\la w \ra$ fixed by 
$P^g,P^h$, we then have $\la v\ra g^{-1}h = \la w \ra$, and the conclusion 
follows. \hal

\begin{cor}\label{liepcor}
If $G$ is as in Lemma $\ref{liep}$ and is $p$-exceptional, then $G$ is transitive on $V^\sharp$.
\end{cor}

\pf By Lemma \ref{liep}, $\hat L$ is transitive on $P_1(V)$. Hence $\hat L$ is given by Hering's Theorem (see \cite[Appendix 1]{liebaff}): so $\hat L$ is $SL_d(q)$, $Sp_d(q)$ or $G_2(q)$ ($q$ even, $d=6$), and in each case $\hat L$ is transitive on $V^\sharp$. \hal

\vspace{2mm}
\no {\bf Remark } The case where $L = {\rm soc}(G/Z) = \,^2\!F_4(2)'$ and $p=2$ is quickly ruled out, as follows. Suppose $G$ is 2-exceptional. The 2-modular characters of $L$ are given in \cite{JLPW}, and the bound in Lemma \ref{bound}(i) implies that $d=26,q=2$. As 13 does not divide $|V^\sharp| = 2^{26}-1$, an element $t \in L$ of order 13 fixes a nonzero vector $v$. By 2-exceptionality, $v$ is also fixed by a Sylow 2-subgroup $P$ of $L$. But $t$ and $P$ generate $L$ (see \cite[p.74]{Atlas}), so this is impossible.

\section{$\mathcal{C}_9$ case II: Alternating groups}\label{c92}

In this section we deal with the case where the simple group $L = {\rm soc}(G/Z)$ is an alternating group.
We prove

\begin{thm}\label{altcase}
Let $G \le \Ga L(V) = \Ga L_d(q)$ ($q=p^a$) be such that $G/Z$ is almost simple with socle $L \cong \AAA_c$, an alternating group with $c \ge 5$, and suppose the perfect preimage $\hat L$ of $L$ in $G$ is absolutely irreducible on $V$ and realisable over no proper subfield of $\F_q$. Assume $G$ is $p$-exceptional and not transitive on $V^\sharp$. Then one of the following holds:

{\rm (i)} $G = \AAA_c$ or $\SSS_c$ with $c = 2^r-2$ or $2^r-1$, with $V$ the deleted permutation module for $G$ over $\F_2$, of dimension $d = c-2$ or $c-1$ respectively;

{\rm (ii)} $\hat L  = SL_2(5) \trianglelefteq G < \Ga L_2(9)$, with orbit sizes $1,40,40$ on vectors.

\no Conversely, the groups $G$ in {\rm (i)} are $p$-exceptional.
\end{thm}

We shall need the following.

\begin{prop}\label{acgen}
If $r$ is a prime with $r \le c$ (and $c\ge 5$), then $\AAA_c$ is generated by two of its Sylow $r$-subgroups.
\end{prop}

\pf If $r = 2$ this follows from \cite{gural}, so assume that $r>2$. 
Write $c = kr+s$ for integers $k,s$ with $0 \le s \le r-1$. 

First consider the case where $k=1$. If $s=0$, observe that $(1\,2\,\cdots r)$ and $(1\,2\,\cdots r)\,(1\,2\,3)$ are both $r$-cycles, and the group they generate contains a 3-cycle, hence is $\AAA_r$. If $s=1$ then $c=r+1$ and the $r$-cycles 
$(1\,2\,\cdots r)$, $(2\,3\,\cdots r+1)$ generate $\AAA_c$ as their commutator contains a 3-cycle. A similar argument applies for $s=2$, taking $(1\,2\,\cdots r)$ and $(3\,4\,\cdots r+2)$: the group these generate is a primitive group containing both an $r$-cycle and a 5-cycle, hence is $\AAA_c$. Finally, if $s>2$ then $(1\,2\,\cdots r)$ and $(s+1\,s+2\,\cdots c)$ generate a primitive group containing an $r$-cycle fixing more than 2 points, hence generate $\AAA_c$ by Jordan's Theorem \cite[13.9]{W}.

Now suppose $k\ge 2$. Take $R$ to be a Sylow $r$-subgroup of $\AAA_c$ containing the $r$-cycles $(1\cdots r)$, 
$(2r-1 \cdots 3r-2)$, $\ldots $, $(2ir-2i+1\cdots (2i+1)r-2i)$ up to the maximal $i$ such that $(2i+1)r-2i \le c$, 
and take $S$ to be a Sylow $r$-subgroup containing the $r$-cycles $(r \cdots 2r-1)$, $(3r-2 \cdots 4r-3)$, $\ldots$,
$((2i-1)r-2i+2 \cdots 2ir-2i+1)$. Then the group generated by $R$ and $S$ contains a subgroup $\AAA_{c-r+1}$, and now we can include a further $r$-cycle in $S$ to generate $\AAA_c$. \hal

\vspace{4mm}
Now we embark upon the proof of Theorem \ref{altcase}. Let $G \le \Ga L(V) = \Ga L_d(q)$ ($q = p^a$) be as in the hypothesis, with $L = {\rm soc}(G/Z) \cong \AAA_c$. 
 Note that $|{\rm Out}(L)| = 2$ or 4, so $p$ divides $G\cap GL(V)$ which is therefore $p$-exceptional by Lemma \ref{normalsubgroups}. Hence we may replace $G$ by $G\cap GL(V)$ and assume that $G \le GL(V)$.
By Lemma \ref{bound}, together with Lemma \ref{acgen}, we have
\begin{equation}\label{gory}
d = \dim V \le 2\log_q(c!) - 2\log_q((\frac{c!}{2})_p),
\end{equation}
where $n_p$ denotes the $p$-part of an integer $n$.

Since $\AAA_5\cong L_2(5) \cong L_2(4)$, $\AAA_6 \cong L_2(9)$ and $\AAA_8 \cong L_4(2)$, Lemma \ref{liep} shows that we may exclude these groups from consideration in these characteristics.

\begin{lemma}\label{dichot}
One of the following holds:

{\rm (i)} $\hat L = L$ and $V$ is the deleted permutation module over $\F_p$;

{\rm (ii)} $c \le 13$ and $d \le 32$.
\end{lemma}

\pf Suppose $V$ is not the deleted permutation module. Assume that $d > 250$. Then the bound (\ref{gory}) forces $c>40$.
By \cite[Theorem 7]{james} (for $\hat L = \AAA_c$) and \cite{wag} (for $\hat L = 2.\AAA_c$), we have $d \ge \frac{1}{4}c(c-5)$. But this is greater than the upper bound for $d$ given by (\ref{gory}) when $c>40$. Hence $d \le 250$.
Now all the possibilities for the $\F_q \hat L$-module $V$ are given by 
\cite{HM}.
We check that the only possibilities for which the bound (\ref{gory}) is satisfied have $c \le 13$ and $d \le 32$. \hal

\begin{lemma}\label{delete}
Suppose $V$ is the deleted permutation module for $\hat L = \AAA_c$ over $\F_p$. Then $p=2$, 
$c = 2^r-1$ or $2^r-2$, and $V$ has dimension $d = c-1$ or $c-2$ respectively. Moreover, in these representations $\AAA_c$ and $\SSS_c$ are $2$-exceptional.
\end{lemma}

\pf Here $V = S/(S\cap T)$ where $S = \{(a_1,\ldots ,a_c) \in \F_p^c \mid \sum a_i = 0\}$, $T = \{(a,\ldots ,a) \mid a \in \F_p\}$. 
If $p\ge 5$ define 
\[
v = (1,2,-3,1,-1,2,-2,\ldots , \frac{p-3}{2}, -\frac{p-3}{2},0,\ldots ,0) \in S
\]
 and let $x = v+(S\cap T) \in V$. Then $(\hat L)_x = \AAA_{c-p}\times H$ where $H \le \AAA_p$ and $|H|$ is coprime to $p$; hence $p$ divides $|x^{\hat L}|$, contradicting the fact that $G$ is $p$-exceptional. If $p=3$ and $c\ne 6$ the same argument applies, taking $v = (1,1,-1,-1,0,\ldots ,0)$; and if $p=3$ and $c=6$ then taking the same $v$, we have $|\hat L_x|_3 = 3$, so $|x^{\hat L}|$ is divisible by $|\hat L|_3/3 = 3$, again a contradiction.

Hence $p=2$. Here $d = \dim V = c-(c,2)$. For $c$ odd, the orbit sizes of $\hat L = \AAA_c$ on nonzero vectors are $\binom{c}{i}$ for $1\le i \le c-1$; and for $c$ even the orbit sizes are $\binom{c}{2i}$ ($i \ne c/2$) and also $\frac{1}{2}\binom{c}{c/2}$ when $4|c$. It follows that we get $2$-exceptional examples when $c = 2^r-1$ or $2^r-2$. And when $c$ is not of one of these forms, let $2^e$ be the smallest power of 2 that is missing in the binary expansion of $c$; then $\binom{c}{2^e}$ is even, so $\hat L$ (and hence $G$) is not 2-exceptional, a contradiction. \hal

\vspace{4mm}
Assume from now on that $V$ is not the deleted permutation module.

\begin{lemma}\label{cless8}
If $V$ is not the deleted permutation module then $c \le 7$.
\end{lemma}

\pf From \cite{HM} and (\ref{gory}), we see that the possibilities for $c,d,p$ and $q$ are as in the following table:
\[
\begin{array}{|l|l|l|l|l|}
\hline
\hbox{Case}&c&d&p&q \\
\hline
1&13&32&3&3 \\
2&12&16&2,3&4,3 \\
3&11& 16 & 2,3,5,7 & p\hbox{ or }p^2 \\
4&10&8&5&5 \\
5&  & 16 & 2,3 & 2,3 \\
6& 9 & 8 & 2,3,5,7 & p \\
7&  & 21 & 3 & 3 \\
8& 8 & 8 & 3,5,7 & p \\
9&   & 13 & 3 & 3 \\
\hline
\end{array}
\]
In cases 1, 2, 3, 4 and 7 we adopt a unified approach: for the respective primes $r = 13,11,11,7$ and 7 we observe that $r$ does not divide $q^d-1$, and hence there is a nonzero vector $v \in V$ fixed by an element $t \in \hat L$ of order $r$. As $G$ is $p$-exceptional, $v$ is also fixed by a Sylow $p$-subgroup $P$ of $\hat L$. However it is easy to see that $\la t,P \ra = \hat L$ in all these cases, so this is a contradiction.

Now consider case 5. Here $c=10$, $p=2$ or 3 and $d = \dim V = 16$. The Brauer character $\chi$ of $V$ is given in \cite{JLPW}. For $p=3$ we have $\chi(t) = 1$ where $t$ is a $5B$-element of $\hat L = 2.\AAA_{10}$; this implies that $t$ fixes a nonzero vector $v$ -- however $\la t,P\ra = \hat L$ for any Sylow $3$-subgroup $P$, so this is not possible. For $p=2$ we apply the same argument with $t$ a $9A$-element of $\hat L = \AAA_{10}$, noting that $t$ and any Sylow 2-subgroup generate $\hat L$.

Next consider case 6. For $p=2$ we have $\hat L = \AAA_9 < O_8^+(2) = O(V)$, and $\hat L$ has two orbits on nonzero vectors (see \cite{Atlas}), one of which has even size, contrary to 2-exceptionality. For $p=3$ or 5 we argue as above that an element $t$ in class $7A$ or $9A$ fixes a nonzero vector, and generates $\hat L$ with any Sylow $p$-subgroup, a contradiction. The case $p=7$ requires a little more argument. Let $S$ be a Sylow 3-subgroup of $\hat L = 2.\AAA_9$. Since $|V| = 7^8 \equiv 4 \hbox{ mod }9$, there is an orbit of $S$ on nonzero vectors of size 1 or 3, and hence there is a nonzero vector $v$ which is fixed by a subgroup $S_0$ of $S$ of index 3. However, $S_0$ and any Sylow 7-subgroup generate $\hat L$, so this is impossible.

Now consider case 8. For $p=5$ we argue with the Brauer character that an element in class $4A$ in $\hat L = 2.\AAA_8$ fixes a nonzero vector; but this element generates $\hat L$ with any Sylow 5-subgroup. For $p=3$, the Brauer character $\chi$ of $V$ shows that an element $t \in \hat L$ of order 7 fixes a nonzero vector $v$. Then $v$ is stabilised by $\la t,P \ra$ for some Sylow 3-subgroup $P$ of $\hat L$, and this contains $2.\AAA_7$. However $\chi \downarrow 2.\AAA_7$ is a sum of two irreducibles of degree 4, so $2.\AAA_7$ does not fix a nonzero vector, a contradiction. For $p=7$ we argue similarly using a $4B$-element $t$.

Finally, consider case 9. Here a $7A$ element $t \in \hat L = \AAA_8$ fixes a nonzero vector $v$, and so $v$ is fixed by 
 $\la t,P \ra$ for some Sylow 3-subgroup $P$ of $\hat L$. This contains $\AAA_7$. However we see from \cite{JLPW} that $\AAA_7$ acts irreducibly on $V$, so this is impossible. \hal

\begin{lemma}\label{c7cas}
If $V$ is not the deleted permutation module then $c$ is not $7$.
\end{lemma}

\pf Suppose $c=7$. Here \cite{HM} and (\ref{gory}) show that $d,q$ are as in the following table of possibilities:
\[
\begin{array}{|l|l|l|}
\hline
\hbox{Case} & d &  q \\
\hline
1&15  &  3 \\
2&13  & 3,5 \\
3&9  & 7 \\
4&8  & 5 \\
5&6  & 4,9,7 \\
6&4  & 2,9,25,7 \\
7&3  & 25 \\
\hline
\end{array}
\]
In cases 1--4, an element $t$ of order 5 or 7 fixes a nonzero vector, and generates $\hat L$ with any Sylow $p$-subgroup, giving a contradiction as usual. 

Consider case 5. For $p=7$, we argue as before with an element of order 5. For $p=3$, the Brauer character shows that an element $t$ of order 5 in $\hat L = 2.\AAA_7$ fixes a nonzero vector $v$, and so $v$ is fixed by$\la t,P \ra$ for some Sylow 3-subgroup $P$ of $\hat L$. This contains $2.\AAA_6$, which is irreducible on $V$, a contradiction. Finally for $p=2$, we have $\hat L = 3.\AAA_7 < 3.M_{22} < SU_6(2) < SL(V)$ (see \cite[p.39]{Atlas}). From \cite{Atlas} we see that $M_{22}$ is transitive on the set of 672 non-isotropic 1-spaces in $V$, and $\AAA_7$ has 3 orbits on these, so that one orbit must have even size, contrary to $p$-exceptionality.

Now consider case 6. For $p=3$ or 5, we argue as usual with an element of order 7 which fixes a vector. For $p=2$ we have $\hat L = \AAA_7 < SL_4(2) = SL(V)$ and $\hat L$ is transitive on $V^\sharp$, contrary to assumption. And for $p=7$ we have $\hat L = 2.\AAA_7 < SL_4(7)$; by \cite[Appendix 2]{liebaff}, $\hat L$ has two orbits on nonzero vectors, one of which has size divisible by 7.

Finally, in case 7 the Brauer character shows that there is an element $t$ of order 3 fixing a nonzero vector $v$, so $v$ is fixed by $\la t,P \ra$ for some Sylow 5-subgroup $P$, and this contains $\AAA_5$; but a subgroup $\AAA_5$ must be irreducible on $V$, a contradiction. \hal 

\begin{lemma}\label{c6cas}
If $V$ is not the deleted permutation module then we have $c=5$, $d=2$, $q=9$ and $\hat L = SL_2(5) < SL_2(9) = SL(V)$, with orbit sizes $40,40$ on nonzero vectors.
\end{lemma}

\pf We know that $c=5$ or 6. Recall that $p=3$ when $c=5$ , and $p\ne 3$ when $c=6$. Hence \cite{HM} and (\ref{gory}) show that $d,q$ are as in the following table of possibilities:
\[
\begin{array}{|l|l|l|l|}
\hline
\hbox{Case} & c& d  & q \\
\hline
1&6 &4  & 2,5 \\
2&6 & 3  & 4,25 \\
3&5 & 6 &3 \\
4&5 & 3 &9 \\
5&5 & 2 &9 \\
\hline
\end{array}
\]
Consider case 1. If $p=2$ then $\hat L < SL_4(2)$ is transitive on nonzero vectors. And if
$p=5$, the Brauer character shows that there is an element $t$ of order 3 fixing a nonzero vector $v$, so the stabiliser of $v$ contains $\la t,P\ra$ for some Sylow 5-subgroup $P$ of $\hat L = 2.\AAA_6$. So this stabiliser contains $2.\AAA_5$; but a subgroup $2.\AAA_5$ is irreducible on $V$, a contradiction. The same argument deals with case 2 for $p=5$; and in case 2 for $p=2$ we have $\hat L = 3.\AAA_6 < SL_3(4) = SL(V)$, and \cite{liebaff} shows that $\hat L$ has 2 orbits on nonzero vectors, one of which has even size.

In cases 3 and 4 we observe that an element of order 5 in $\hat L$ fixes a vector and generates $\hat L$ with any Sylow 3-subgroup. Finally, case 5 gives the example in the conclusion. \hal

\vspace{4mm}
This completes the proof of Theorem \ref{altcase}.

\section{$\mathcal{C}_9$ case III: Sporadic groups}\label{c93}

In this section we deal with the case where the simple group $L = {\rm soc}(G/Z)$ is a sporadic group.
We prove

\begin{thm}\label{sporcase}
Let $G \le \Ga L(V) = \Ga L_d(q)$ ($q=p^a$) be such that $G/Z$ is almost simple with socle $L$ a sporadic group, and suppose the perfect preimage $\hat L$ of $L$ in $G$ is absolutely irreducible on $V$ and realisable over no proper subfield of $\F_q$. 
Assume $G$ is $p$-exceptional  on $V$. Then one of the following holds:

{\rm (i)} $M_{11}\trianglelefteq G < GL_5(3)$, orbit sizes $1,22,220$;

{\rm (ii)} $M_{23}=G < GL_{11}(2)$, orbit sizes $1,23,253,1771$.
\end{thm}

We shall need the following result, analogous to Proposition \ref{acgen}.

\begin{prop}\label{2genspor}
Let $T$ be a sporadic simple group and $p$ a prime dividing $T$. Then there exists Sylow $p$-subgroups $S_1$, $S_2$ such that $T=\langle S_1,S_2 \rangle$.
\end{prop}

\pf 
For all sporadic simple groups except for $Th$, $J_4$, $Ly$, $B$ and $M$, permutation representations exist in \textsc{Magma} \cite{Mag} (for the larger ones using the generators given in the online Atlas \cite{atlas}). This allows a
Sylow $p$-subgroup to be constructed and a conjugate that generates $T$ can then be found. 

For the remaining five sporadics we use three strategies. Let $H$ be a maximal subgroup of $T$ such that $p$ divides $|T:H|$ and suppose that $H$ has Sylow $p$-subgroups $S_1$ and $S_2$ such that $\langle S_1,S_2\rangle= H$ or a normal subgroup of index coprime to $p$. Then for each $i\in\{1,2\}$ we can find a Sylow $p$-subgroup $S_i'$ of $T$ properly containing $S_i$ and since $H$ is maximal in $T$ it follows that $\langle  S_1',S_2'\rangle=T$. Here we list the sporadic groups, primes $p$ and maximal subgroups $H$ for which we used this method:
\[
\begin{array}{|l|l|l|}
\hline
T & p & H\\
\hline
Th  & 3, 5, 7 &   2^5.L_5(2) \\
      & 2       & \SSS_5 \\
J_4  & 2, 3, 11 & L_2(32):5 \\
Ly  & 2, 3  &   5^3.L_3(5) \\
      & 5     & 2\AAA_{11} \\
B &  2, 3, 5 & Th \\
   & 7   & Fi_{23} \\
M &2, 3, 5, 7, 11,13 & 2B\\
\hline
\end{array}
\]

Next, if the only maximal subgroup of $T$ containing a Sylow $p$-subgroup $S$ normalises $S$ then $T$ will be generated by any pair of Sylow $p$-subgroups. This is true for $(T,p)=(J_4,43)$,  $(Ly,67)$, and  $(B,47)$. Finally, if the order of a Sylow $p$-subgroup of $T$ is $p$ and there are two elements of order $p$ whose product has order a prime $r$ such that there are no maximal subgroups of $T$ with order divisible by $pr$ then $T$ is generated by these two elements of order $p$. The existence of such primes can be checked by either doing random searches using the matrix representations available in the online Atlas \cite{atlas} or by doing character table calculations in \textsf{GAP} \cite{gap}. This method was used for the groups and primes listed below. 
\[
\begin{array}{|l|l|l|}
\hline
T & p & r\\
\hline
Th & 13, 19& 31\\
     &31 &19\\
J_4 &5, 23, 29, 31, 37 & 43\\
      & 7  & 37\\
Ly &7, 31, 37 & 67\\
     & 11 & 37\\
B &11, 13, 17, 19, 31 & 47\\
    & 23 & 31\\
M & 17, 19, 23, 29, 31, 41, 47, 59  &71\\ 
    & 71  & 59\\
\hline
\end{array}
\]\hal

\vspace{4mm}
Now we embark upon the proof of Theorem \ref{sporcase}. Let $G \le \Ga L(V) = \Ga L_d(q)$ ($q = p^a$) be as in the hypothesis, with $L = {\rm soc}(G/Z)$ a sporadic group. Note that $|{\rm Out}(L)| \le 2$, so as in the previous section we may assume that $G \le GL(V)$.  By Lemma \ref{bound}, together with Proposition \ref{2genspor}, we have
\begin{equation}\label{sporgory}
d = \dim V \le 2\log_q|G:N_G(P)|,
\end{equation}
where $P \in Syl_p(G)$. 

\begin{lemma} $L$ is not $M$, $BM$, $Fi_{24}$, $Fi_{23}$, $Th$, $Ly$, $HN$ or $O'N$.
\end{lemma}

\pf Suppose $L$ is one of these groups other than $M$. Then (\ref{sporgory}) implies that $d < 250$, whence $L$ and $V$ are in the list in \cite{HM}. We check that for all possibilities, (\ref{sporgory}) fails. And if $L = M$ then (\ref{sporgory}) gives $d <360$, whereas any nontrivial representation of $M$ has dimension greater than this (see \cite[5.3.8]{KL}). \hal

\begin{lemma} $L$ is not $J_4$, $Fi_{22}$, $Co_1$, $Co_2$, $Co_3$, $Suz$, $Ru$ or $He$.
\end{lemma}

\pf Suppose $L$ is one of these groups. Then $d < 250$ by (\ref{sporgory}), so \cite{HM} together with (\ref{sporgory}) imply that $L,d,q$ are as in the following table:

\[
\begin{array}{|l|l|l|}
\hline
L & d & q \\
\hline
J_4 & 112 & 2 \\
Fi_{22} & 27,\,78 & 4, 2 \\
Co_1 & 24 & p\,(\hbox{any}) \\
Co_2 & 22 &  2 \\
         & 23 &  p\,(p>2) \\
Co_3 & 22 & 2,3 \\
         & 23 & p\,(p>3) \\
Suz & 12 & p \hbox{ or }p^2 \\
Ru & 28 & 2 \\
He & 51 & 2 \\
\hline
\end{array}
\]

We deal with each of these in turn.

Let $L = J_4$. From \cite[p.190]{Atlas} we see that there is a maximal subgroup $H = 2^{10}.SL_5(2)$ fixing a nonzero vector $v \in V$. Hence $v$ is stabilised by the subgroup generated by $H$ and a Sylow 2-subgroup of $L$, which is $L$, a contradiction.

Now let $L = Fi_{22}$. Here $(d,q) = (27,4)$ or $(78,2)$. Since $q^d-1$ is not divisible by 13, there is an element $t$ of order 13 fixing a vector, and $L$ is generated by $t$ together with any Sylow 2-subgroup.

Next consider $L = Co_1$. A subgroup $H = Co_2$ of $L$ fixes a nonzero vector $v$, and together with any Sylow $p$-subgroup generates $L$, provided $p\ne 11$ or 23. For $p=11,23$ consider $V \downarrow M_{24}$. Using \cite{JLPW} we see that this is $V_1 \oplus V_{23}$, where $V_{23}$ is the deleted permutation module for $M_{24}$. Hence for an element $t \in L$ of order $p$ we have $\dim C_V(t) = 4$ or 2 according as $p=11$ or 23. Now Lemma \ref{bound}(ii) gives a contradiction.

Now let $L = Co_2$. If $p=2$ we see from \cite[p.154]{Atlas} that there is a vector stabilised by $U_6(2).2$, and this generates $L$ with any Sylow 2-subgroup. If $p=3$ or 5 then there is a vector stabilised by an element of order 23, and this generates $L$ with any Sylow 2-subgroup. Finally, for $p\ge 7$, consideration of  $V \downarrow M_{23}$ shows that $\dim C_V(t)\le 5$ for an element $t$ of order $p$, and now Lemma \ref{bound} gives a contradiction.

Now consider $L=Co_3$. For $p=2$, there is a cyclic subgroup $H = C_{11}$ fixing a vector; and for $p=3,5$  subgroups $H = McL, M_{23}$, respectively, fix a vector. Now observe that in each case $H$ generates $L$ with any Sylow $p$-subgroup. And for $p \ge 7$ we argue just as for $Co_2$.

Next let $L = Suz$. For $p\ne 3,11$, there is an element of order 11 fixing a vector, and this generates $L$ with any Sylow $p$-subgroup. For $p=3$ there is a subgroup $U_5(2)$ fixing a vector and generating $L$ with any Sylow 3-subgroup. Finally, if $p=11$ then $q=121$ and (\ref{sporgory}) fails.

Now let $L=Ru$. Here a 13-element fixes a vector and generates $L$ with any Sylow 2-subgroup.

Finally, for $L = He$, a 17-element fixes a vector and generates $L$ with any Sylow 2-subgroup. \hal 

\begin{lemma} $L$ is not $McL$, $HS$, $J_1$, $J_2$ or $J_3$.
\end{lemma}

\pf  Suppose $L$ is one of these groups. Then \cite{HM} together with (\ref{sporgory}) imply that $L,d,q$ are as in the following table:

\[
\begin{array}{|l|l|l|}
\hline
L & d & q \\
\hline
McL & 21 & 3,5 \\
      & 22 & 2,7 \\
HS & 20 & 2 \\
     & 21 & 5 \\
     & 22 & 3 \\
J_3 & 9 & 4 \\
      & 18 & 4 \\
J_2 & 6 & 4,9,5,49 \\
      & 14 & 4,3,5,7,49 \\
      & 36 & 2 \\
J_1 & 7 & 11 \\
      & 20 & 2 \\
\hline
\end{array}
\]

Let $L = McL$. For each $q$ we produce a subgroup $H$ of $L$ stabilizing a nonzero vector: for $q=2,3$, take $H = 5^{1+2}$ (as 5 is coprime to $|V^\sharp|$); for $q=5$ take $H = C_{11}$; and for $q=7$ take $H = M_{11}$, noting that from \cite{JLPW} the Brauer character of $L$ on $V$ restricts to $M_{11}$ as $1+\chi_{10}+\chi_{11}$. From \cite[p.100]{Atlas}, we check that each of these subgroups $H$ generates $L$ together with any Sylow $p$-subgroup, which is a contradiction.

Now consider $L = HS$. For $q=2$, we see from \cite{JLPW} that the value of the Brauer character of $L$ on elements in classes $5B,5C$ is 0, and hence such elements fix a nonzero vector in $V$. By \cite[p.80]{Atlas}, $L$ has only two classes of maximal subgroups of odd index, and between them they meet only one class of elements of order 5. This is a contradiction. For $q=3$ observe that the restriction of $V$ to a subgroup $M_{22}$ is $V_1+V_{21}$ where $V_{21}$ is the deleted permutation
module. Now $M_{22}$ contains a Sylow 3-subgroup of $L$, and if $t \in M_{22}$ is an element of order 3 acting on 22 points with cycle-type $(3^6,1^4)$, then $\dim C_V(t) = 10$. This leads to a contradiction via Lemma \ref{bound}(ii). Finally for $q=5$, $V$ restricts to $M_{22}$ as the deleted permutation module, and hence $M_{21} = L_3(4)$ fixes a vector; but $L_3(4)$ generates $L$ with any Sylow 5-subgroup.

Next let $L = J_3$. Here an element of order 17 fixes a vector, and generates $L$ with any Sylow 2-subgroup.

Now let $L = J_2$. For $d= 14$ or 36 we produce a subgroup $H$ of $L$ fixing a vector and generating $L$ with any Sylow $p$-subgroup, as follows: if $p=2$ or 5, $H = C_7$ (if $p=2$ use the Brauer character value); if $p=3$, $H = 5^2$; and if $p=7$, $H = C_5$ generated by a $5A$-element (which fixes a vector by consideration of its Brauer character value). For $d=6$ we have $L < PSp_6(q)$ for $q = 4,5,9$ or 49. In the first three cases the orbit sizes of $\hat L$ on vectors can be found in \cite{liebaff} ($q=4,5$) and \cite[Table 5]{affdtg} ($q=9$), and in each case there is a size divisible by $p$, a contradiction. And for $q=49$ a $5A$-element fixes a vector and generates $L$ with any Sylow 7-subgroup.

Finally let $L = J_1$. Here an element of order 19 fixes a vector in both modules, and generates $L$ with any Sylow $p$-subgroup. \hal

\begin{lemma} If $L$ is a Mathieu group then one of the following holds:

{\rm (i)} $M_{11}\trianglelefteq G < GL_5(3)$, orbit sizes $1,22,220$;

{\rm (ii)} $M_{23}=G < GL_{11}(2)$, orbit sizes $1,23,253,1771$.
\end{lemma}

\pf Here (\ref{sporgory}) and \cite{HM} imply that the possibilities for $L, d, q$ are as follows:

\[
\begin{array}{|l|l|l|}
\hline
L & d & q \\
\hline
M_{24} & 11,44 & 2 \\
             & 22 & 3 \\
             & 23 & 5 \\
M_{23} & 11,44 & 2 \\
            & 22 & 3 \\
M_{22}  & 10 & 2,9,25,7,11 \\
             & 6,15,34 & 4,4,2 \\             
              & 21 & 3 \\
M_{12} &10 & 2 \\
            & 6,10,15 & 3 \\
            & 11 & 5 \\
M_{11} & 10 & 2 \\
              & 5,10 & 3 \\
\hline 
\end{array}
\]

Let $L = M_{24}$. For $q=2$, an element of order 11 fixes a vector and generates $L$ with any Sylow 2-subgroup. And for $q=3$ or 5, $V$ is the deleted permutation module and the orbits of the vectors $(1,1,-1,-1,0^{20})$ ($q=3$) and 
$(1,1,2,2,-1,0^{19})$ ($q=5$) have size divisible by $p$.

Now consider $L=M_{23}$. For $q=2$ and $d=11$ there are two possible modules. In both cases the orbit sizes are given in \cite[p.170]{affdtg}, and one of them gives the 2-exceptional example in conclusion (ii) of the lemma. For $d=44$ we argue as follows. The group $L$ has one class of involutions, and can be generated by three of them; hence $\dim C_V(t) \le \frac{2}{3} \dim V < 30$. Now Lemma \ref{bound}(ii) gives a contradiction. Finally, for $q=3$, $V$ is the deleted permutation module and the orbit of the vector $(1,1,-1,-1,0^{19})$ has size divisible by 3.

Next let $L = M_{22}$. Consider $p=2$. For $d=10,15$ or 34, an element of order 7 fixes a vector and generates $L$ with any Sylow 2-subgroup. And for $d=6$ we have $\hat L = 3.M_{22} < SU_6(2)$ and from \cite[p.39]{Atlas} we see that $\hat L$ is transitive on the 672 non-isotropic vectors in $V$, contrary to 2-exceptionality. For $p=3$ and $d=21$, $V$ is the deleted permutation module, dealt with in the usual way. It remains to consider $d=10$ in characteristics $3,5,7$ and 11. For $p=11$, an element of order 7 fixes a vector and generates $L$ with any Sylow $p$-subgroup. And for $p=3,5,7$ we have $\hat L = 
2.M_{22} < SL_{10}(q)$ with $q = 9,25,7$ respectively, and a \textsc{Magma} computation reveals the existence of orbits of size divisible by $p$. 

Now let $L=M_{12}$. For $q=2$ or 5, $V$ is the deleted permutation module, dealt with in the usual way. Next let $q=3$. 
For $d=10$ or 15, an element of order 5 fixes a vector and generates $L$ with any Sylow 3-subgroup. And for $d=6$ we have $\hat L = 2.M_{12} < SL_6(3)$ and \cite[p.170]{affdtg} shows that there is an orbit size divisible by 3.   

Finally, let $L = M_{11}$. For $q=2$,  $V$ is the deleted permutation module. For $q=3$ and $d=5$ the orbit sizes are given by \cite{liebaff}, giving the example in part (i) of the lemma. And for $q=3$, $d=10$ there are three possible modules $V$, 
and a \textsc{Magma} computation shows that for each of these there is an orbit size divisible by 3. \hal

\vspace{4mm}
This completes the proof of Theorem \ref{sporcase}.

\section{$\mathcal{C}_9$ case IV: ${\rm Lie}(p')$}\label{c94}

In this section we deal with the case where the simple group $L = {\rm soc}(G/Z)$ is a simple group of Lie type in $p'$-characteristic. We prove

\begin{thm}\label{liecase}
Let $G \le \Ga L(V) = \Ga L_d(q)$ ($q=p^a$) be such that $G/Z$ is almost simple with socle $L$ a simple group of Lie type in $p'$-characteristic, and $L$ is not isomorphic to an alternating group. Suppose the perfect preimage $\hat L$ of $L$ in $G$ is absolutely irreducible on $V$ and realisable over no proper subfield of $\F_q$. 
Assume $G$ is $p$-exceptional  on $V$ and is not transitive on $V^\sharp$. Then $L_2(11) \trianglelefteq G < GL_5(3)$, with orbit sizes $1,22,110,110$.
\end{thm}

As in previous sections, for the proof it is useful to know that the group $L$ is generated by any two of its Sylow $p$-subgroups. For $p=2$ this is true by \cite{gural}. It is presumably true for all primes, but we do not prove this here, and just cover the following groups for which we need the result.

\begin{lemma}\label{punygen}
Let $L$ be one of the following simple groups:
\[
\begin{array}{l}
L_n(r):\;n=2,r\le 37; \hbox{ or }n=3,r\le 5; \hbox{ or }n\le 5,r=3; \hbox{ or }n\le 8,r=2 \\
PSp_{2n}(r):\; n=2,r\le 11; \hbox{ or }n=3,r\le 7; \hbox{ or }n=4,r=5; \hbox{ or }n\le 6, r\le 3 \\
U_n(r):\; n=3,r\le 5; \hbox{ or }n=4,r\le 4; \hbox{ or }n=5,r\le 3; \hbox{ or }n\le 11,r=2 \\
\hbox{and }\Om_7(3),\Om^\pm_8(2),\Om_{10}^\pm(2),F_4(2),^2\!F_4(2)',^3\!D_4(2),G_2(3),G_2(4),^2\!B_2(8).
\end{array}
\]
Let $X$ be a group such that $L \le X \le {\rm Aut}(L)$. If $p$ is a prime dividing $|X|$, then $O^{p'}(X)$ is generated by two of its Sylow $p$-subgroups.
\end{lemma}

\pf This was verified computationally using \textsc{Magma}. \hal

\vspace{4mm}
Now we embark upon the proof of Theorem \ref{altcase}. Let $G \le \Ga L(V) = \Ga L_d(q)$ ($q = p^a$) be as in the hypothesis, with $L = {\rm soc}(G/Z)$ a simple group of Lie type in $p'$-characteristic not isomorphic to an alternating group. We assume now that $G\le GL(V)$; we shall handle the case where $G\le \Ga L(V)$ at the end of the proof.

\begin{lemma}\label{list}
The simple group $L$ is one of the following:
\[
\begin{array}{l}
L_2(r)\,(r\le 37),\,L_3(4),\,L_3(5),\,L_4(3), \\
PSp_4(r)\,(r\le 9,r\ne 8),\,PSp_6(r)\,(r=2,3,5),\,PSp_8(3),\,PSp_{10}(3), \\
U_3(r)\,(r=3,4,5),\,U_4(2),\,U_5(2),\,U_6(2), \\
\Om_7(3),\,\Om_8^\pm(2), \\
F_4(2),\,^2\!F_4(2)',\,^3\!D_4(2),\,G_2(3),\,G_2(4),\,^2\!B_2(8).
\end{array}
\]
\end{lemma}

\pf By Lemma \ref{bound}(i) we have $d = \dim V \le r_p\log_q|G:N_G(P)|$, where $r_p$ is the minimal number of Sylow $p$-subgroups required to generate $O^{p'}(G)$. We have $r_2 = 2$ by \cite{gural}, while upper bounds for $r_p$ with $p$ odd are provided by \cite[Sections 3,4,5]{gusa}. On the other hand, lower bounds for $d$ are given in \cite{LS, SZ}.
These, together with the above bound, imply that $L$ is one of the groups listed in Lemma \ref{punygen}. By that lemma, we have $r_p=2$ for these groups, and so we now have the bound
\begin{equation}\label{liegory}
d = \dim V \le 2\log_q|G:N_G(P)|.
\end{equation}
Applying this with the above-mentioned lower bounds for $d$ eliminates many of the groups in the list, and leaves just the groups in the conclusion. \hal

\begin{lemma} \label{l2}
If $L = L_2(r)$, then $r=11$ and $L_2(11) \trianglelefteq G < GL_5(3)$, with orbit sizes $1,22,110,110$.
\end{lemma}

\pf Assume $L = L_2(r)$. By Lemma \ref{list}, $r\le 37$, so certainly $\dim V < 250$ and we can use \cite{HM} together with (\ref{liegory}) to identify the possibilities for $L,d,q$:
\[
\begin{array}{|l|l|l|}
\hline
L & d & q \\
\hline
L_2(31) & 15 & 2 \\
L_2(25) & 12 & 2 \\
L_2(23) & 11 & 2,3 \\
L_2(17) & 8 & 2 \\
L_2(13) & 6 & 4,3 \\
             & 7 & 3 \\
L_2(11) & 10 & 2 \\
             & 5 & 4,3,5 \\
             & 6 & 3 \\
L_2(7) & 3 & 9 \\
\hline
\end{array}
\]
The case where $r=11$ and $(d,q) = (5,3)$ gives the example in the conclusion; the orbit sizes are given in \cite{CMS}. (Note that there are two 5-dimensional representations of $L$ over $\F_3$, but they are quasiequivalent, hence have the same orbit sizes.) Also, when $r=13$ and $(d,q) = (6,3)$, $\hat L = SL_2(13)$ is transitive on nonzero vectors (see \cite[Appendix 1]{liebaff}), contrary to our assumption in Theorem \ref{liecase}.

When $r= 11, 13, 17, 23$ or 25 and $(d,q) = (5,5)$, $(6,4)$, $(8,2)$, $(11,2)$ or $(12,2)$, the orbit sizes are given by \cite[Section 5]{CMS} and \cite[Appendix 2]{liebaff}; there is an orbit size divisible by $p$ in all cases.

Now consider the cases where $r=11, 13, 23$ or 31 and $(d,q) = (6,3)$, $(7,3)$, $(11,3)$ or $(15,2)$. For these, we
observe that there is a nonzero vector fixed by a subgroup $H$ of order 11,7,11 or 5 respectively, and $H$ generates $L$ together with any Sylow $p$-subgroup, a contradiction.

This leaves the cases where $r=7$ or 11 and $(d,q) = (3,9),\,(5,4)$ or $(10,2)$. For these cases a \textsc{Magma} computation shows that there is an orbit of size divisible by $p$. \hal

\begin{lemma}
$L$ is not $L_3(4),\,L_3(5)$ or $L_4(3)$.
\end{lemma}

\pf Suppose $L$ is one of these groups. By \cite{HM} and (\ref{liegory}), the possibilities for $L,d,q$ are:
\[
\begin{array}{|l|l|l|}
\hline
L & d & q \\
\hline
L_3(4) & 4 & 9 \\
          & 6 & 3,7 \\
          & 8 & 5 \\
L_4(3) & 26 & 2\\
\hline
\end{array}
\]
Consider the case where $L = L_3(4)$ and $(d,q) = (6,3)$. Here $L < P\Om_6^-(3)$ and we see from \cite[p.52]{Atlas} that $L$ is transitive on the 126 non-isotropic points in $V$, contrary to $p$-exceptionality.

For the remaining cases $(d,q) = (4,9),\,(8,5),\,(6,7),\,(26,2)$, 
there is a nonzero vector fixed by a subgroup $H$ of order 7, 7, 5 or 13 respectively, and $H$ generates $L$ together with any Sylow $p$-subgroup, a contradiction. \hal

\begin{lemma}
$L$ is not $PSp_{2m}(r)$ for $m \ge 2$.
\end{lemma}

\pf Suppose $L = PSp_{2m}(r)$, so that $L$ is as in Lemma \ref{list}. Using \cite{HM} and (\ref{liegory}), we see that the possibilities for $L,d,q$ are:
\[
\begin{array}{|l|l|l|}
\hline
L & d & q \\
\hline
PSp_4(3) & 4 & 25 \\
               & 6 & 5 \\
PSp_4(5) & 12 & 4 \\
PSp_4(7) & 24 & 2 \\
PSp_4(9) & 40 & 2 \\
PSp_6(3) & 13 & 4,7 \\
               & 14 & 7 \\
PSp_8(3) & 40 & 4 \\
Sp_4(4) & 18 & 3 \\
Sp_6(2) & 7 & 3,5,7 \\
             & 8 & 3,5,7 \\
             & 14 & 3 \\
\hline
\end{array}
\]
The three cases $PSp_4(3) < L_6(5)$, $Sp_6(2) < L_7(5)$ and $Sp_6(2) < L_8(7)$ were handled using a \textsc{Magma} computation, and in each case there is an orbit size divisible by $p$. 

Now consider all the other cases in the table with $r$ odd. In each case $r^2$ does not divide $q^d-1$, so there is a subgroup
$H$ fixing a nonzero vector, where $H$ is of index $r$ in a Sylow $r$-subgroup of $L$ (replacing $r$ by 3 when $L=PSp_4(9)$). However, one checks that $H$ generates $L$ with any Sylow $p$-subgroup, contradicting $p$-exceptionality (use \cite{Atlas}, or \cite{odddeg} for $p=2$).

If $L=Sp_4(4)$, an element of order 17 fixes a vector and generates $L$ with any Sylow 3-subgroup. 

Finally let $L=Sp_6(2)$. For the $p=3$ cases, an element of order 7 in $L$ fixes a vector and generates $L$ with any Sylow 3-subgroup. And when $(d,q) = (7,7)$ or $(8,5)$ there are subgroups $\Om_6^-(2)$ or $U_3(3)$, respectively, fixing a vector, and these generate $L$ with any Sylow $p$-subgroup. \hal

\begin{lemma}
$L$ is not $U_n(r)$ for $n \ge 3$.
\end{lemma}

\pf Suppose $L = U_n(r)$ is as in Lemma \ref{list}. By \cite{HM} and (\ref{liegory}), the possibilities are:
\[
\begin{array}{|l|l|l|}
\hline
L & d & q \\
\hline
U_3(3) & 6 & 2 \\
               & 14 & 2 \\
U_3(4) & 12 & 3 \\
U_3(5) & 20 & 2 \\
U_5(2) & 10 & 3,5,11 \\
U_6(2) & 21 & 3 \\
               & 22 & 5 \\
\hline
\end{array}
\]
The group $U_3(3) < L_6(2)$ is transitive on nonzero vectors (see \cite[Appendix 1]{liebaff}), contrary to assumption. 
For the cases $(d,q) = (14,2)$, $(20,2)$, $(10,3)$, $(10,11)$ and $(21,3)$, there is a subgroup $H$ or order $7,7,2^7,2^7,11$, respectively, fixing a vector and generating $L$ with any Sylow $p$-subgroup. In the remaining three cases with $(d,q) = (12,3)$, $(10,5)$ and $(22,5)$, a \textsc{Magma} computation shows that there is an orbit of size divisible by $p$. \hal

\begin{lemma}
$L$ is not an orthogonal group in Lemma $\ref{list}$.
\end{lemma}

\pf Suppose $L = \Om_7(3)$ or $\Om_8^\pm(2)$. Using \cite{HM} and (\ref{liegory}), we see that $L = \Om_8^+(2)$, $d=8$ and $q = 3,5$ or 7. Here $\hat L = 2.L$ and these 8-dimensional representations arise from the action of the Weyl group of $E_8$ on the root lattice. This is transitive on the 240 root vectors, and hence there is an orbit of size 120 on 1-spaces. Hence $\hat L$ is not $p$-exceptional for $p = 3,5$. For $p=7$, the total number of 1-spaces in $V$ is not divisible by 3, so there is a 1-space fixed by a Sylow 3-subgroup of $L$, which generates $L$ with any Sylow 7-subgroup of $L$. \hal

\begin{lemma}
$L$ is not an exceptional group in Lemma $\ref{list}$.
\end{lemma}

\pf Suppose $L$ is an exceptional group.  By \cite{HM} and (\ref{liegory}), the possibilities are:
\[
\begin{array}{|l|l|l|}
\hline
L & d & q \\
\hline
G_2(3) & 14 & 2 \\
G_2(4) & 12 & 3,5,7 \\
^2\!B_2(8) & 8 & 5 \\
\hline
\end{array}
\]
For $(d,q) = (14,2),(12,3),(12,5),(12,7)$, there is a subgroup $H$ of $L$ of order $3^5,2^9,2^{10},2^8$, respectively, fixing a 1-space, and $H$ generates $L$ with any Sylow $p$-subgroup, a contradiction. 
For the remaining case $(d,q) = (8,5)$, a \textsc{Magma} computation gives an orbit of size divisible by 5. \hal

\vspace{4mm}
We now complete the proof of Theorem \ref{liecase} by handling the case where the $p$-exceptional group $G$ in the hypothesis  lies in $\Ga L(V)$ and not in $GL(V)$, where $V = V_d(q)$. Let $G_0 = G\cap GL(V)$. If $p$ divides $|G_0|$ then $G_0$ is $p$-exceptional, hence is given by the linear case of the theorem which we have already proved; but this implies that $GL(V) = GL_5(3)$, in which case $GL(V) = \Ga L(V)$.

Hence $p$ does not divide $|G_0|$, and by Lemma \ref{pexc}, $G$ has a $p$-exceptional normal subgroup $G_0\la \s \ra$ (hence also $\hat L\la \s \ra$), where $\s \in \Ga L(V)\backslash GL(V)$ is a field automorphism of order $p$. Hence $q=p^{kp}$ for some integer $k$. Since $G_0$ is a $p'$-group, we have $p>2$. Moreover, $\s$ induces an automorphism of order $p$ of the simple $p'$-group $L$, which must be a field automorphism, so $L = L(r^p)$ is a group of Lie type over $\F_{r^p}$ for some~$r$. 

Let $\ell$ be the untwisted Lie rank of $L$. Then $d \ge \frac{1}{2}(r^{p\ell}-1)$ by \cite{LS}, and also $|\hat L| < r^{4p\ell^2}$ and 
$|C_V(\s)| = q^{d(1-\frac{1}{p})}$. By Lemma \ref{bound}(ii), $|V:C_V(\s)|<|\s^L|$, which implies
\[
q^{ \frac{1}{2}(r^{p\ell}-1)\cdot (1-\frac{1}{p})} <  r^{4p\ell^2}.
\]
This cannot hold.

Hence the case where $G$ lies in $\Ga L(V)$ and not in $GL(V)$ does not occur. This completes the proof of Theorem \ref{liecase}.

\section{Proof of Theorem \ref{prim}}\label{primpf}

Let $G$ be an irreducible subgroup of $GL_n(p)$ with $p$ prime, and suppose
that $G$ acts primitively on $V_n(p)$. Suppose also that $G$ is $p$-exceptional,
so that $p$ divides $|G|$ and every orbit of $G$ on $V$ has size coprime
to $p$. 

Choose $q = p^k$ maximal such that $G \le \Gamma L_d(q) \le GL_n(p)$, where
$n = dk$. Write $V = V_d(q)$, $G_0 = G \cap GL_d(q)$ and $Z = G_0\cap \F_q^*$,
the group of scalar multiples of the identity in $G_0$. Note that $G_0 \trianglelefteq G$ and $G/G_0$ is cyclic.
Write $K = \F_q$.

If $d=1$ then $G \le \Gamma L_1(q)$, as in Theorem \ref{prim}(ii).
So assume that $d \ge 2$.

\begin{lemma} \label{irred2} $G_0$ is absolutely irreducible on $V=V_d(q)$.
\end{lemma}

\pf We have $G \le N_{GL_n(p)}(K) = \Ga L_d(q)$. So $G_0 = C_H(K) \le GL_d(q)$. 

Let $E = \End_G(V) = \F_r \subseteq K$, and write $q=p^k = r^b$. 
Viewing $V$ as $V_{bd}(r)$, it is an absolutely irreducible $\F_rG$-module.
Now view $V$ as an $\F_qG_0$-module. Then $U:= V \otimes_{\F_r} \F_q$, as an $\F_qG_0$-module,
is the sum of $b$ Frobenius twists of $V$. However $G/G_0$ is cyclic of order
at most $b$, so if $V \downarrow G_0$ were reducible, then $U \downarrow
G$ would be reducible. But $G$ is absolutely irreducible, 
so this is a contradiction. 

Hence $V \downarrow G_0$ is irreducible.
As $C_{\End(V)}(G_0)$ is a field extension of $K$, the choice of $K$ implies
that $C_{\End(V)}(G_0)=K$, and so $V$ is an absolutely irreducible $KG_0$-module. \hal

\vspace{4mm}
If $G$ preserves a tensor product decomposition $V = U \otimes
W$ over $\F_q$, where $\dim U \ge \dim W \ge 2$ (i.e. $G \le {\Ga L(V)}_{U \otimes W}$), 
then Theorems \ref{c4lin} and
\ref{c4semi} (together with Lemma \ref{pexc}) give a contradiction. So assume that $G$ does
not preserve a nontrivial tensor decomposition of $V$. 

\begin{lemma}\label{normal} Let $N$ be a normal subgroup of $G$ such that $N
\le G_0$ and $N \not \le Z$. Then $V \downarrow N$ is absolutely irreducible.
\end{lemma}

\pf By Clifford's Theorem $V \downarrow N$ is a direct sum of homogeneous components; these are permuted by $G$, and hence by the primitivity of $G$, $V\downarrow N$ is homogeneous. Say $V\downarrow N \cong 
W \oplus \cdots \oplus W$, a direct sum of $k$ copies of an irreducible $KN$-module $W$. Let $K_0 = C_{\End(W)}(N)$, a field extension of $K$. By the first few lines of the proof of  \cite[5.7]{asch}, $K_0^*$ can be identified with $Z(C_{GL(V)}(N))$ and $G \le N_{\Ga L(V)}(K_0)$, so $K_0 = K$ by choice of $K$. Hence $W$ is an absolutely irreducible $KN$-module. At this point \cite[3.13]{asch} shows that there is a $K$-space $A$ such that $V$ can be identified with $W \otimes A$ in such a way that $N \le GL(W) \otimes 1_A$, $G_0 \le GL(W)\otimes GL(A)$ and $G \le N_{\Ga L(V)}(GL(W) \otimes GL(A))$. By our assumption that $G$ preserves no nontrivial tensor decomposition, this implies that $W=V$, completing the proof. \hal

\vspace{4mm}
Now let $S$ be the socle of $G/Z$, and write $S = M_1\times \cdots 
\times M_k$ where each $M_i$ is a minimal normal subgroup of $G/Z$. 
Let $R$ be the full preimage of $S$ in $G$, and $P_i$ the preimage 
of $M_i$, so that $R = P_1\ldots P_k$, a commuting product. If some
$P_i$, say $P_k$, is not contained in $G_0$, then $M_k$ is generated 
by a field automorphism $\phi$ of prime order, and so $G/Z \le 
C_{P\Gamma L(V)}(\phi) = PGL_d(q_0)\la \psi \ra$, where $\psi$ 
generates the Galois group of $\F_q/\F_p$ and $\F_{q_0}$ is a proper
subfield of $\F_q$. This contradicts Theorem \ref{sub}.

Hence $R = P_1\ldots P_k\le G_0$. As $P_1 \trianglelefteq G$, 
Lemma \ref{normal} implies that $V \downarrow P_1$ is absolutely
irreducible, hence $C_{G_0}(P_1) = Z$. It follows that $k=1$ and 
$R=P_1$. Also $G$ is not realised (modulo scalars) over a proper subfield of $\F_q$, by Theorem \ref{sub}.

Suppose first that $M_1=R/Z$ is an elementary abelian $r$-group for 
some prime $r$, and replace $R$ by a minimal preimage of $M_1$ in 
$G$. By Lemma \ref{normal} and since $d\ge 2$, any maximal abelian characteristic 
subgroup of $R$ must be contained in $Z$. Hence $R$ is of 
symplectic type, and now we argue as in \cite[11.8]{asch} that $R$ can be taken to be as in Section \ref{c6sec}.
Since $G$ is primitive, Theorem \ref{c6case} now gives a contradiction.

Now suppose $M_1 = R/Z$ is non-abelian, so $M_1 \cong T^l$ for some
non-abelian simple group $T$, and $R = R_1\ldots R_l$ where  
$R_i/Z \cong T$ and the factors are permuted transitively by $G$. 
If $l>1$ then \cite[3.16, 3.17]{asch} implies that $R$ preserves a 
tensor decomposition $V = V_1\otimes \cdots \otimes V_l$ with 
$\dim V_i$ independent of $i$, and $G \le N_{\Gamma L(V)}(\bigotimes 
GL(V_i))$; then Theorems \ref{c7m} and \ref{c10m} give a contradiction.

It remains to consider the case where $l=1$, so that 
$M_1 = {\rm Soc}(G/Z)$ is simple. Then $G/Z$ is almost simple, 
and has socle absolutely irreducible on $V$ and not realisable 
over a proper subfield of $\F_q$. In other words $G$ is in the 
class $\mathcal{C}_9$, and so $G$ is given by the results in 
Sections \ref{c91}-\ref{c94}.

This completes the proof of Theorem \ref{prim}.

\section{Deduction of Theorems \ref{gen}, \ref{tiepnav} and \ref{threehalf}}\label{ded}

\vspace{4mm}
\no {\bf Proof of Theorem \ref{gen} }
Let $G \le GL_d(p) =  GL(V)$ be irreducible and $p$-exceptional, and let $G_0 = O^{p'}(G)$. Then $G_0$ is also 
$p$-exceptional by Lemma \ref{normalsubgroups}. By Clifford's Theorem, $V\downarrow G_0 = V_1\oplus \cdots \oplus V_t$, where the $V_i$ are irreducible $G_0$-modules, conjugate under $G$. Note that $O^{p'}(G_0) = G_0$ and $p$ divides $|G_0^{V_1}|$, so $G_0^{V_1}$ is $p$-exceptional.  
If $G_0^{V_1}$ is primitive, it is given by Theorem \ref{prim}, and if it is imprimitive it is given by Theorem \ref{imprim}. 

We now claim that the $V_i$ are pairwise non-isomorphic $G_0$-modules. Suppose false, and let $W = V_1\oplus
\cdots \oplus V_k$ be a homogeneous component for $G_0$ with $k>1$. Then $G_0^W \le GL(V_1) \otimes 1
\le GL(V_1) \otimes GL_k(p)$. Since $G_0$ is $p$-exceptional on $V$, it is also $p$-exceptional on $W$.
We now apply Theorem \ref{c4lin}, which classifies $p$-exceptional groups which preserve tensor product decompositions. From this theorem, it is clear that such a group cannot act just as scalars on one of the tensor factors, which is what $G_0^W$ does. This is a contradiction, proving the claim.

By the claim, $G$ permutes the summands $V_i$. Finally, the kernel $K$ of the action of $G$ on the set of summands $\{V_1,\ldots ,V_t\}$ contains $G_0$, so $G/K$ is a $p'$-group.

\vspace{4mm}
\no {\bf Proof of Theorem \ref{tiepnav} } 
Let $G$ be a finite group
and let $p > 2$ be a prime. Assume that $G = O^{p'}(G) = O^{p}(G)$, that 
$G$ has abelian Sylow $p$-subgroups, and that
$V$ is a faithful irreducible $\F_pG$-module  
such that every orbit of $G$ on $V$ has length coprime to $p$. These assumptions 
imply that $G=G'$ and also that $p$ divides $|G|$. 

Suppose first that $G$ acts primitively on $V$. Then $G$ is given by Theorem \ref{prim}. As $G$ is insoluble, it is either transitive on $V^\sharp$ or one of the examples in (iii) of the theorem. In the first case we refer to the list of transitive linear groups in \cite[Appendix]{liebaff}: the only examples where $G=G'$ and $G$ has abelian Sylow $p$-subgroups occur in conclusions (i) and (iii) of Theorem \ref{tiepnav}; and the examples in Theorem \ref{prim}(iii) with $p>2$ are also in conclusion (iii).

Now suppose $G$ is imprimitive on $V$. As $G = O^{p'}(G)$, $G$ is given by Theorem \ref{imprim}. Theorem \ref{p-conc} and the assumptions that $p>2$ and $G$ has abelian Sylow $p$-subgroups, now force $G$ to be as in conclusion (ii) of Theorem  \ref{tiepnav}. This completes the proof of the corollary.

\vspace{4mm}
\no {\bf Proof of Theorem \ref{threehalf} }
Let  $G\le GL_d(p)$ be a $\frac{1}{2}$-transitive linear group of order divisible by $p$, and
write $V = V_d(p)$. Let $H = VG \le AGL_d(p)$, the corresponding affine permutation group acting on $V$.
Since $G$ has order divisible by $p$, it does not act semiregularly on $V^\sharp$, and so $H$ is $\frac{3}{2}$-transitive on $V$, and hence is a primitive permutation group on $V$ by \cite[10.4]{W}. This implies that $G$ acts irreducibly on $V$.

Since $G$ is $\frac{1}{2}$-transitive and has order divisible by $p$, it is $p$-exceptional.
So if $G$ acts primitively as a linear group on $V$, then 
it is given by Theorem \ref{prim}. For $G =A_c$ or $S_c$ as in (iii)(a) of the theorem, the orbit sizes are
given in the proof of Lemma \ref{delete}, and we see that the only $\frac{1}{2}$-transitive example is
for $c=6$ with $(d,p) = (4,2)$, in which case $G$ is transitive on $V^\sharp$.
Hence $G$ is as in the conclusion of Theorem \ref{threehalf}.
Finally, if $G$ acts imprimitively on $V$ then it is given by \cite[Theorem 1.1]{passman2} (which determines all imprimitive $\frac{1}{2}$-transitive linear groups). The only example of order divisible by $p$ is $D_{18} < \Gamma L_1(2^6) < GL_6(2)$, as in (ii) of Theorem \ref{threehalf}.


\begin{thebibliography}{99}

\bibitem{asch} M. Aschbacher, On the maximal subgroups of the finite classical groups, 
{\it Invent. Math.} {\bf 76} (1984), 469--514. 

\bibitem{bam} J. Bamberg, M. Giudici, M.W. Liebeck, C.E. Praeger and J. Saxl, 
The classification of almost simple $\frac{3}{2}$-transitive groups, {\it Trans. Amer. Math. Soc.}, \textbf{365} (2013), 4257--4311.

\bibitem{affdtg} J. van Bon, A.A. Ivanov and J. Saxl, Affine distance-transitive graphs with sporadic stabilizer, 
{\it Europ. J. Comb.} {\bf 20} (1999), 163--177.

\bibitem{Mag} W. Bosma, J. Cannon and C. Playoust, The \textsc{Magma} algebra system I: 
The user language, {\it J. Symbolic Comput.} {\bf 24} (1997) 235--265.

\bibitem{CNS} P.J. Cameron, P.M. Neumann and J. Saxl, On groups with no regular orbits on the set of subsets,  
{\it Arch. Math.}  {\bf 43}  (1984), 295--296.

\bibitem{Carter}
R.W. Carter, \emph{Finite groups of Lie type. Conjugacy classes and complex characters}, Wiley, Chichester, 1993. 

\bibitem{CMS} A.M. Cohen, K. Magaard and S. Shpectorov,  Affine distance-transitive graphs: the cross characteristic case, {\it Europ. J. Comb.} {\bf 20} (1999), 351--373.

\bibitem{Atlas}
  J. H. Conway, R. T. Curtis, S. P. Norton, R. A. Parker, and R. A. Wilson,
{\it Atlas of Finite Groups}, Clarendon Press, Oxford, $1985$.

\bibitem{curtis} C.W. Curtis, Modular representations of finite groups with split $(B,\,N)$-pairs, in {\it Seminar on Algebraic Groups and Related Finite Groups}, Springer Lecture Notes {\bf 131},  pp.57--95, Springer, Berlin, 1970.

\bibitem{Dickson}
L. E. Dickson, \emph{Linear Groups With an Exposition of the Galois Field Theory},
 Dover, New York, 1958.

\bibitem{dol} S. Dolfi, R. Guralnick, C.E. Praeger and P. Spiga, Coprime subdegrees for primitive permutation groups and completely reducible linear groups, {\it Israel J. Math.}, to appear.

\bibitem{dorn} L. Dornhoff, {\it Group representation theory, Part B: 
Modular representation theory}, Pure and Applied Mathematics, Vol.7, Marcel Dekker, Inc., New York, 1972.

\bibitem{Fe} W. Feit, {\it The representation theory of finite groups},
  North-Holland, Amsterdam, $1982$.

\bibitem{GW}
A. Gambini Weigel and T. S. Weigel,
On the orders of primitive linear $p'$-groups. 
\emph{Bull. Austral. Math. Soc.} {\bf 48} (1993),  495--521. 

\bibitem{gap} 
The GAP group, `{\it {\sf GAP} - groups, algorithms, and
programming}', Version 4.4, 2004, {\sf http://www.gap-system.org}.

\bibitem{kgv} D. Gluck, K. Magaard, U. Riese and P. Schmid, 
The solution of the $k(GV)$ problem, {\it J. Algebra} {\bf 279} (2004), 694--719. 

\bibitem{GW1}
D. Gluck and T. Wolf, 
Defect groups and character heights in blocks of solvable groups. II, 
{\it J. Algebra} {\bf 87} (1984), 222--246. 

\bibitem{GW2}
D. Gluck and T. Wolf, Brauer's height conjecture for $p$-solvable
groups,
{\it Trans. Amer. Math. Soc.} {\bf 282} (1984), 137--152. 

\bibitem{goodwin} D. Goodwin, Regular orbits of linear groups with an application to the $k(GV)$-problem I, II, 
{\it  J. Algebra}  {\bf 227}  (2000),  395--432 and 433--473.

\bibitem{GorLy3}
D. Gorenstein, R. Lyons and R. Solomon,  \emph{The classification of the finite simple groups. Number 3. Part I. Chapter A. Almost simple $K$-groups.} Mathematical Surveys and Monographs, 40.3. American Mathematical Society, Providence, RI, 1998.

\bibitem{gural} R. Guralnick, Generation of simple groups, {\it J. Algebra} {\bf 103} (1986), 381--401.

\bibitem{gusa} R.M. Guralnick and J. Saxl, Generation of finite almost simple groups by conjugates, {\it J. Algebra}
{\bf 268} (2003), 519--571.

\bibitem{her} C. Hering, Transitive linear groups and linear groups which contain irreducible subgroups of prime order II,  {\it J. Algebra} {\bf 93}  (1985),  151--164.

\bibitem{HM} G. Hiss and G. Malle, Low-dimensional representations of quasi-simple groups (Corrigenda),  {\it LMS J. Comput. Math.} {\bf 5} (2002), 95--126. 

\bibitem{james} G.D. James, On the minimal dimensions of irreducible representations of symmetric groups, {\it Math. Proc. Camb. Phil. Soc.} {\bf 94} (1983), 417--424.

\bibitem{JLPW}
  C. Jansen, K. Lux, R. A. Parker, and R. A. Wilson, {\it An Atlas of
Brauer Characters}, Oxford University Press, Oxford, 1995.

\bibitem{KL}
  P. B. Kleidman and M. W. Liebeck, {\it The Subgroup Structure of the
Finite Classical Groups}, London Math. Soc. Lecture Note Ser. no.
$129$, Cambridge University Press, $1990$.

\bibitem{KP} C. K\"ohler and H. Pahlings, Regular orbits and the $k(GV)$-problem, in {\it Groups and computation, III (Columbus, OH, 1999)}, pp.209--228, Ohio State Univ. Math. Res. Inst. Publ., 8, de Gruyter, Berlin, 2001.

\bibitem{LS} V. Landazuri and G.M. Seitz, On the minimal degrees of projective representations of the finite Chevalley groups, {\it J. Algebra} {\bf 32} (1974), 418--443.
 
\bibitem{liebaff} M.W. Liebeck, The affine permutation groups of rank three,
{\it Proc. London Math. Soc.}, {\bf 54} (1987), 477--516.

\bibitem{regsgps}
M. W. Liebeck, C. E. Praeger, J. Saxl, 
Regular subgroups of primitive permutation groups. 
{\em Mem. Amer. Math. Soc.} {\bf 203} (2010), no. 952, vi+74 pp. 

\bibitem{odddeg} M.W. Liebeck and J. Saxl,  The primitive permutation groups of odd degree, 
{\it  J. London Math. Soc.}  {\bf 31}  (1985), 250--264.

\bibitem{malle} G. Malle, Fast-einfache Gruppen mit langen Bahnen in absolut irreduzibler Operation,  
{\it J. Algebra} {\bf 300}  (2006),  655--672.

\bibitem{Mu}
  M. Murai, On Brauer's height zero conjecture, \emph{Proc. Japan Acad. Ser. A Math. Sci.} \textbf{88} (2012), 38--40.

\bibitem{NS}
  G. Navarro and B. Sp\"ath, A reduction theorem for the Brauer 
height zero conjecture, \emph{J. Eur. Math. Soc.}, (to appear).


\bibitem{NT1}
  G. Navarro and P. H. Tiep, Brauer's height zero conjecture for the 
$2$-blocks of maximal defect, {\it J. Reine Angew. Math.} \textbf{669} (2012), 225--247.


\bibitem{NT2} G. Navarro and P.H. Tiep, Characters of $p'$-degree over normal 
subgroups, \emph{Annals of Math.} \textbf{178} (2013), 1135--1171.

\bibitem{passman1} D. Passman, Solvable $\frac{3}{2}$-transitive group, {\it J. Algebra} {\bf 7} (1967), 192--207.

\bibitem{passman2} D.S. Passman, Exceptional $\frac{3}{2}$-transitive permutation groups, {\it Pacific J. Math.} {\bf 29} (1969), 669--713.

\bibitem{SZ} G.M. Seitz and A.E. Zalesskii, On the minimal degrees of projective representations of the finite Chevalley groups II, {\it J. Algebra} {\bf 158} (1993), 233--243.


\bibitem{Ser}
 A. Seress, Primitive groups with no regular orbits on the set of subsets,
{\it J. London Math. Soc.} {\bf 29} (1997), 697--704.


\bibitem{TZ2}
  P. H. Tiep and A. E. Zalesskii, Unipotent elements of finite groups of 
Lie type and realization fields of their complex representations, 
{\it J. Algebra} {\bf 271} $(2004)$, 327--390.

\bibitem{wag} A. Wagner, An observation on the degrees of projective representations of the symmetric and alternating groups over an arbitrary field, {\it Arch. Math.} {\bf 29} (1977), 583--589.

\bibitem{ward} H.N. Ward, Representations of symplectic groups, 
{\it J. Algebra} {\bf 20} (1972),  182--195.

\bibitem{W} H. Wielandt, {\it  Finite Permutation Groups}, Academic Press, New York, 1964.

\bibitem{atlas} R.A. Wilson, et al., A World-Wide-Web Atlas of finite group 
representations, http://brauer.maths.qmul.ac.uk/Atlas/v3/.

\bibitem{Z} H. Zassenhaus, \"Uber endliche Fastk\"orper, 
{\it Abh. Math. Semin. Hamb. Univ.}  {\bf 11} (1935), pp.187--220.

\end{thebibliography}
\end{document}